\documentclass{amsart}

\usepackage{amsmath,amssymb,amscd,bbm}
\usepackage[mathscr]{eucal}
\usepackage[all]{xy}
\usepackage[utf8]{inputenc}
\usepackage{nicefrac}
\usepackage{mathtools}
\usepackage{color}
\usepackage{upgreek}
\usepackage{stmaryrd}
\usepackage{hyperref}

 \newtheorem{thm}{Theorem}[section]
 \newtheorem{cor}[thm]{Corollary}
 \newtheorem{lem}[thm]{Lemma}
  
 \newtheorem{prop}[thm]{Proposition}
 \theoremstyle{definition}
 \newtheorem{defn}[thm]{Definition}
 \theoremstyle{remark}
 \newtheorem{rem}[thm]{Remark}
 \newtheorem{rems}[thm]{Remarks}
 \newtheorem{exa}[thm]{Example}
 \newtheorem{exas}[thm]{Examples}
 
 \numberwithin{equation}{section}

\usepackage{hapasymb}

\newcommand{\type}{\mathrm{type}}
\newcommand{\cotype}{\mathrm{cotype}}

\newcommand{\dps}{\displaystyle}

\newcommand{\sect}[1]{\mathrm{S}_{#1}}
\newcommand{\cstrip}[1]{\overline{\mathrm{St}}_{#1}}
\newcommand{\csect}[1]{\overline{\mathrm{S}}_{#1}^\ast}

\newcommand{\gauss}{\mathrm{G}}  
\newcommand{\Laplace}{\mathit{\Delta}}
\newcommand{\gammabound}[1]{\llbracket #1 \rrbracket_\gamma}

\renewcommand{\Im}{\im}
\renewcommand{\Re}{\re}

\newcommand{\Ha}{\mathrm{H}^\infty}
\newcommand{\Ho}{\Ha_0}

\newcommand{\Hr}{\mathrm{H}\ddot{\mathrm{o}}\mathrm{r}}

\newcommand{\UCb}{\mathrm{UC}_{\mathrm{b}}}
\newcommand{\Test}{\Cc^\infty}

\newcommand{\emdf}[1]{\textbf{#1}}
\newcommand{\red}[1]{\textcolor{red}{#1}}

\newcommand{\spec}{\usigma}

\newcommand{\resol}{\urho}

\newcommand{\vanish}[1]{\relax}

\newcommand{\beq}{\begin{equation}}
\newcommand{\eeq}{\end{equation}}

\newcommand{\hs}{\hskip-0.1em}
\newcommand{\set}[1]{\hs\left[\,#1\,\right]}

\newcommand{\prfnoi}{\smallskip\noindent}

\newcommand{\Exp}{\mathbb{E}}

\newcommand{\coo}{\mathrm{c}_{00}}

\newcommand{\bignorm}[1]{\bigl\| #1 \bigr\|}
\newcommand{\Bignorm}[1]{\Bigl\| #1 \Bigr\|}

\newcommand{\Id}{\mathrm{I}}

\newcommand{\suchthat}{\,\,|\,\,}

\newcommand{\N}{\mathbb{N}}
\newcommand{\Z}{\mathbb{Z}}

\newcommand{\R}{\mathbb{R}}
\newcommand{\C}{\mathbb{C}}

\newcommand{\We}{\mathrm{W}}

\newcommand{\sector}[1]{\mathrm{S}_{#1}}

\newcommand{\strip}[1]{\mathrm{St}_{#1}}

\newcommand{\res}[1]{|_{#1}}

\newcommand{\Sum}[2][\relax]{%
 \ifx#1\relax \sideset{}{_{#2}}\sum 
 \else \sideset{}{^{#1}_{#2}}\sum
 \fi}

\newcommand{\car}{\mathbf{1}}
\DeclareMathOperator{\re}{Re}
\DeclareMathOperator{\im}{Im}
\newcommand{\konj}[1]{\overline{#1}}

\newcommand{\abs}[1]{\vert #1 \vert}
\newcommand{\bigabs}[1]{\bigl\vert #1 \bigr\vert}

\DeclareMathOperator{\Hol}{Hol}

\newcommand{\Ce}{\mathrm{C}}
\newcommand{\Cc}{\mathrm{C}_\mathrm{c}}
\newcommand{\Co}{\mathrm{C}_{0}}
\newcommand{\Cb}{\mathrm{C}_{\mathrm{b}}}
\newcommand{\BUC}{\mathrm{BUC}}

\newcommand{\loc}{\text{\upshape{\tiny \bfseries loc}}}

\newcommand{\Har}[1]{\mathrm{H}^{#1}}
\newcommand{\Ell}[1]{\mathrm{L}^{#1}}

\newcommand{\ohne}{\setminus}

\newcommand{\leer}{\emptyset}

\newcommand{\dann}{\Rightarrow}

\DeclareMathOperator{\dom}{dom}
\DeclareMathOperator{\ran}{ran}

\newcommand{\cls}[1]{\overline{#1}}
\newcommand{\cl}[1]{\overline{#1}}

\newcommand{\rand}{\partial}
\DeclareMathOperator{\supp}{supp}

\newcommand{\spann}{\mathrm{span}}

\newcommand{\BL}{\mathcal{L}}

\newcommand{\norm}[2][\relax]{%
   \ifx#1\relax \ensuremath{\lVert#2\rVert}
   \else \ensuremath{\left\Vert#2\right\Vert_{#1}}
   \fi}

\makeatletter
\newcommand{\sprod}[2]{\ensuremath{%
  \setbox0=\hbox{\ensuremath{#2}}
  \dimen@\ht0
  \advance\dimen@ by \dp0
  \left[ #1\rule[-\dp0]{0pt}{\dimen@}, #2\hspace{1pt}\right]}}
\newcommand{\bsprod}[2]{\ensuremath{%
  \setbox0=\hbox{\ensuremath{#2}}
  \dimen@\ht0
  \advance\dimen@ by \dp0
  \bigl[ #1\rule[-\dp0]{0pt}{\dimen@}, #2\hspace{1pt}\bigr]}}

 \makeatother
\newcommand{\dprod}[2]{\ensuremath{\langle#1,#2\rangle}}

\newcommand{\Fourier}{\mathcal{F}}
\newcommand{\fourier}[1]{\widehat{#1}}

\makeatletter
\renewcommand{\tocsection}[3]{%
  \indentlabel{\@ifnotempty{#2}{\ignorespaces#1 #2\quad}}#3}
\renewcommand{\tocsubsection}[3]{%
  \indentlabel{\@ifnotempty{#2}{\ignorespaces#1 #2\quad}}#3}

\newcommand\@dotsep{4.5}
\def\@tocline#1#2#3#4#5#6#7{\relax
  \ifnum #1>\c@tocdepth 
  \else
    \par \addpenalty\@secpenalty\addvspace{#2}%
    \begingroup \hyphenpenalty\@M
    \@ifempty{#4}{%
      \@tempdima\csname r@tocindent\number#1\endcsname\relax
    }{%
      \@tempdima#4\relax
    }%
    \parindent\z@ \leftskip#3\relax \advance\leftskip\@tempdima\relax
    \rightskip\@pnumwidth plus1em \parfillskip-\@pnumwidth
    #5\leavevmode\hskip-\@tempdima{#6}\nobreak
    \leaders\hbox{$\m@th\mkern \@dotsep mu\hbox{.}\mkern \@dotsep mu$}\hfill
    \nobreak
    \hbox to\@pnumwidth{\@tocpagenum{\ifnum#1=1\fi#7}}\par
    \nobreak
    \endgroup
  \fi}
\AtBeginDocument{%
\expandafter\renewcommand\csname r@tocindent0\endcsname{0pt}
}
\def\l@subsection{\@tocline{2}{0pt}{2.5pc}{5pc}{}}
\makeatother

\makeatletter
\ifcsname phantomsection\endcsname
    \newcommand*{\qrr@gobblenexttocentry}[5]{}
\else
    \newcommand*{\qrr@gobblenexttocentry}[4]{}
\fi
\newcommand*{\addsubsection}{%
    \addtocontents{toc}{\protect\qrr@gobblenexttocentry}%
    \subsection}
\makeatother

\begin{document}

\title[Holomorphic Hörmander-Type Functional Calculus]{Holomorphic Hörmander-Type Functional Calculus on Sectors and Strips}

\author[Markus Haase]{Markus Haase}

\address{%
Mathematisches Seminar, Kiel University\\
Heinrich-Hecht-Platz 6\\
24118 Kiel, Germany}

\email{haase@math.uni-kiel.de}

\author[Florian Pannasch]{Florian Pannasch}

\email{pannaschflorian@gmail.com}




\date{\today}

\begin{abstract}
In this paper, recent abstract multiplier theorems for $0$-sectorial
and $0$-strip type operators by Kriegler and Weis (2018) are refined
and 
generalized to arbitrary sectorial and strip-type operators. To this end, 
holomorphic Hörmander-type functions
on sectors and strips are introduced
with a finer scale of smoothness than the
classical polynomial scale.  
Moreover, we establish alternative descriptions
of these spaces involving Schwartz and ``holomorphic Schwartz'' functions.  
Finally,  the abstract results are combined with
a recent result by Carbonaro and Dragičević (2017) to obtain an 
improvement---with respect to the smoothness condition---of the 
known Hörmander-type multiplier theorem for general symmetric
contraction semigroups. 
\end{abstract}

\maketitle
\setcounter{tocdepth}{2}
\tableofcontents

\allowdisplaybreaks
\raggedbottom

\section{Introduction}\label{s.int}

In 1960, Hörmander made an important contribution to the theory of
Fourier  multiplier theorems by proving  what is now known as
the  {\em Mikhlin--Hörmander multiplier theorem}, see \cite[Theorem 2.5]{Hoe1960}:

\begin{thm}[Hörmander, 1960]\label{int.t.hor}
Let $m\in \Ell{\infty} (\R_{>0})$, $N,  \, d\in \N$ with $N > \tfrac{d}{2}$, and suppose that
\begin{equation}\label{int.eq.hor}
\sup_{R> 0} \int_{\frac{R}{2}}^{2R} \big| s^k m^{(k)} (s) \big|^2 \, \frac{\ud s}{s}\, <\, \infty
\end{equation} for all $1\leq k \leq N$. Then, for every $1<p<\infty $, the operator
\[
\calS (\R^d) \to \Ell{p} (\R^d), \quad f\mapsto \Fourier^{-1}\big( m(|\bfx|^2) \widehat{f} \big) = m(-\Laplace)f
\] extends to a bounded operator on $\Ell{p}(\R^d)$.
\end{thm} 

Here, $\calS (\R^d)$ denotes the space of Schwartz functions  on
$\R^d$. It has been 
noted subsequently that Hörmander's theorem can be extended 
to all functions $f$ satisfying 
\begin{equation}\label{int.eq.hor-frac}
\sup_{t>0} \| \eta \cdot m(t\bfs)  \|_{\We^{\beta, 2}(\R)}\, <\, \infty,
\end{equation} 
where $\eta\in \Test(\R_{>0})$ is any
non-zero test function with support contained in $\R_{> 0} = (0, \infty)$, 
$\We^{\beta, 2}(\R)$ is the usual fractional
$\Ell{2}$-Sobolev space, and  $\beta > \frac{d}{2}$.

Hörmander's theorem is a result about the functional calculus
for the operator $A = -\Laplace$, where the Laplacian $\Laplace= -A$ is the
generator of the Gauss-Weierstrass semigroup on $\R^d$. 
A classical generalization of this setting is towards (negative)
generators $-A$ of so-called {\em symmetric contraction
  semigroups}. That means, $-A$ generates a $C_0$-semigroup $(T_t)_{t\ge 0}$
of symmetric contractions on $\Ell{2}(\Omega)$ for some
$\sigma$-finite measure space $\Omega$, such that each operator $T_t$ 
extends to a contraction on each $\Ell{p}$-space\footnote{see Section
\ref{s.app} for a slightly more detailed account}. 
As $A$ is positive and self-adjoint on $\Ell{2}$, the spectral theorem provides a 
functional calculus in which the operator $m(A)$
is  (well-defined and) bounded on $\Ell{2}$ whenever $m: \R_{\ge 0} \to
\C$ is bounded and measurable.  One says that $m$
is an {\em $\Ell{p}$-spectral multiplier} if  
$m(A)$ extends to a bounded operator on $\Ell{p}(\Omega)$. For $1\le p < \infty$
this property is equivalent to  an estimate of the form
\[
\|m(A) f\|_p\, \leq \, c \|f\|_p \qquad (f\in \Ell{2}\cap \Ell{p}).
\] 
Generalizing (or adapting) Hörmander's theorem to this setting amounts
to proving that $m$ is an $\Ell{p}$-spectral multiplier if it
satisfies \eqref{int.eq.hor-frac} for some $\beta>0$.  Clearly, as the 
situation is more abstract than before, this requires additional assumptions.
The following
seminal result was obtained by Meda (\cite[Theorem 4]{Me1990}).

\begin{thm}[Meda, 1990]
Let $-A$ be the (injective) generator of a symmetric contraction
semigroup on some $\sigma$-finite measure space, and let $1<p<\infty$.
Suppose, in addition, that there are constants $c, \, \alpha \geq 0$ with
\[
\| A^{-\ui s} f\|_p \, \leq\, c (1+|s|)^{\alpha} \| f \|_p \qquad (f\in \Ell{p}, \, s\in \R).
\] 
Then $A$ has a bounded $\Hr^{\beta, 2}(\R_{>0})$-calculus on
$\Ell{p}$ for all $\beta > \alpha + 1$.
\end{thm}

Here, $\Hr^{\beta, 2}(\R_{>0})$ is 
our notation for the space of all
functions $m:\R_{>0}\to \C$ satisfying (\ref{int.eq.hor-frac}). 
(Whether this condition holds or not is independent of the
chosen test function $\eta$, see also Remark \ref{sec.r.cHr}.)
Note that here we require $A$ to be injective, to the effect that $m$ need not
be defined at $0$. This restriction is inessential, cf. Remark \ref{sec.r.injpart}.

Meda's theorem has been generalized in two directions. The first one is
by García-Cuerva, Mauceri, Meda, Sjögren, and Torrea  
and allows for an additional exponential factor in the assumptions
about the growth of the imaginary powers.

\begin{thm}[Garcia-Cuerva et al., 2001]\label{int.t.GCetal}
Let $-A$ be the (injective) generator of a symmetric contraction semigroup on some $\sigma$-finite measure space, let $1<p<\infty$, and suppose there are constants $c, \, \alpha \geq 0$, and an angle $\omega \in (0, \tfrac{\pi}{2})$ with
\beq\label{int.eq.GCetal}
\| A^{-\ui s} f\|_p \, \leq\, c (1+|s|)^{\alpha}\ue^{\omega |s|} \| f \|_p \qquad (f\in \Ell{p}, \, s\in \R).
\eeq
Then $m(A)$ extends to a bounded operator on $\Ell{p}$ for every
$m\in \Hr^{\beta, 2}(\sector{\omega})$, $\beta > \alpha +1 $.
\end{thm}

Here, $\Hr^{\beta, 2}(\sector{\omega})$ is the space of all functions $m$
which are bounded and holomorphic on the sector $\sector{\omega} := \{
z\in \C\ohne\{0\} \suchthat \abs{\arg z} < \omega\}$  such
that its boundary values on the rays $\ue^{\pm \ui \omega} \cdot \R_{>
  0}$ are
contained in the ``Hörmander space'' $\Hr^{\beta, 2}(\R_{>0})$.  The sector
$\sector{\omega}$ as the (essential) domain of
definition of $m$ is natural here, as the exponential growth condition on the
imaginary powers corresponds to the $\Ell{p}$-spectrum of
$A$ being contained in $\cls{\sector{\omega}}$. 

More recently, Carbonaro and Dragičević in \cite{CarDra2017} could show
that {\em each} negative generator $A$ of a symmetric contraction semigroup 
satisfies the growth estimate \eqref{int.eq.GCetal} with $\alpha = \frac{1}{2}$
and $\omega = \omega_p$, where
\[ \omega_p \coloneqq \arcsin \abs{1- \tfrac{2}{p}},
\]
cf. Theorem \ref{app.t.CD}.
Combining this result with Theorem \ref{int.t.GCetal}
then led them to the following corollary, which up to now has been  
the strongest result about the Hörmander
calculus for arbitrary generators of symmetric contraction semigroups.

\begin{cor}[Carbonaro and Dragičević, 2017]\label{int.c.CD}
Let $-A$ be the (injective) generator of a symmetric contraction
semigroup on some $\sigma$-finite measure space, and let $1<p<\infty$.
Then $m(A)$ extends to a bounded operator on $\Ell{p}$ for every
$m\in \Hr^{\beta,2}(\sector{\omega_p})$ and $\beta > \frac{3}{2}$. 
\end{cor}

A different generalization of Meda's
theorem was achieved by 
Kriegler and Weis in  \cite{Kr2009, KrWe2017, KrWe2018})
building on  the holomorphic functional calculus
for sectorial and strip-type operators on general Banach spaces.
The authors  construct, for a special class of $0$-sectorial and $0$-strip
type operators, an abstract Hörmander functional calculus and infer
abstract multiplier theorems (= boundedness of this Hörmander
calculus) under additional  ``geometrical'' conditions. 
As an application of their theory, they improved Meda's result
as follows. (See \cite[Theorem 6.1.(2)]{KrWe2018} for the even more
general actual formulation.)

\begin{thm}[Kriegler and Weis, 2018]\label{int.t.KrWe}
Let, on some Banach space $X$, $A$  be a $0$-sectorial operator with a bounded $\Ha
(\sect{\theta})$-calculus for some angle $\theta \in (0, \pi)$, and
let $\alpha\ge 0$.  
\begin{aufzi}
\item Suppose that  $\tfrac{1}{2} > \tfrac{1}{\mathrm{type}(X)} -
  \tfrac{1}{\mathrm{cotype}(X)}$ and that the set $\{ (1 +
  \abs{s})^{-\alpha} A^{-\ui s} \suchthat s\in \R\}$ is bounded.
Then $A$ has a bounded $\Hr^{\beta, 2}(\R_{>0})$-calculus for all
$\beta > \alpha + \tfrac{1}{2}$. If, in addition, $X$ has Pisier's
contraction property, then this calculus is $R$-bounded.

\item Suppose that the set  $\bigl\{ (1 +\abs{s})^{-\alpha} A^{-\ui s} \suchthat s\in \R\}$
is semi-$R$-bounded and $X$ has  Pisier's
contraction property, then $A$ has an $R$-bounded
$\Hr^{\beta, 2}(\R_{> 0})$-calculus for each $\beta > \alpha + \frac{1}{2}$.
\end{aufzi}
\end{thm}

Here, $\mathrm{type}(X)$ and $\mathrm{cotype}(X)$ denote
the Rademacher type and cotype, respectively, of the Banach space
$X$.  For the other notions see Appendix \ref{s.alp}. 

Meda's result is contained in Theorem \ref{int.t.KrWe}. Indeed, 
it has been proved already by 
 Cowling in \cite[Theorem 2]{Co1983} that  if $-A$ is the (injective) 
generator of a symmetric contraction semigroup then $A$ 
has a bounded $H^\infty$-calculus on some sector. (This is
also a consequence of Corollary \ref{int.c.CD}.) Moreover, if
$X = \Ell{p}(\Omega)$ for $1 < p < \infty$ then 
$\type(X) = \min\{2, p\}$, $\cotype(X) = \max\{2, p\}$,  and $X$ has Pisier's
contraction property. Hence, part a) of Theorem \ref{int.t.KrWe} can
be applied and one obtains
an $R$-bounded $\Hr^{\beta, 2}(\R_{> 0})$-calculus for each $\beta > \alpha
+ \frac{1}{2}$ (and not just  a bounded  calculus for $\beta >
\alpha +1$ as in Meda's theorem).

\subsection*{Main results}

In the present paper we generalize the Kriegler--Weis results towards
general sectorial (i.e., not necessarily $0$-sectorial) operators. 
Moreover, following a suggestion by
Strichartz from \cite{Str1971}  we refine the ``smoothness'' 
conditions in the description of the Hörmander and underlying Sobolev
spaces employing  more general (not necessarily polynomial) weights, the so-called
``strongly admissible'' functions (see Section \ref{s.adm}).

\begin{thm}\label{int.t.main}
Let $A$ be a densely defined operator with dense range and a bounded
$\Ha(\sector{\theta})$-calculus  for some $\theta > 0$
on some Banach space $X$. Suppose, in addition,  that 
\[    \norm{A^{-\ui s}} \le \tilde{v}(s) \ue^{\omega \abs{s}}\qquad
(s\in \R)
\]
for some $\omega \ge 0$ and some measurable function $\tilde{v}: \R
\to \R_{>0}$. Let, furthermore, $v: \R \to [1, \infty)$ be strongly
admissible.
Then the following assertions hold.
\begin{aufzi}
\item $A$ has a bounded
$\Hr^2_v(\sector{\omega})$-calculus whenever $\frac{\tilde{v}}{v} \in
\Ell{1}(\R)$.
 
\item If $r\in [1,2]$ is such that $\frac{1}{r} > \frac{1}{\type(X)}-
  \frac{1}{\cotype(X)}$ 
  then $A$ has a bounded
$\Hr^2_v(\sector{\omega})$-calculus whenever $\frac{\tilde{v}}{v} \in \Ell{r}(\R)$.
If, in addition, $X$ has Pisier's 
contraction property, then this calculus is $R$-bounded.

\item If $\bigl\{ \tilde{v}(s)^{-1} \ue^{-\omega
    \abs{s}} A^{-\ui s} \suchthat s\in \R\}$ is semi-$R$-bounded and 
$X$ has Pisier's  contraction property,
then $A$ has an $R$-bounded
$\Hr^2_v(\sector{\omega})$-calculus whenever $\frac{\tilde{v}}{v} \in \Ell{2}(\R)$.
\end{aufzi}
\end{thm}

Parts b) and c) generalize Theorem \ref{int.t.KrWe} 
(take $\omega = 0$, $\tilde{v}(s) = (1 + \abs{s})^\alpha$, $v(s) = (1
+ \abs{s})^\beta$). Part a)
has no analogue in \cite{KrWe2018} but is somehow the  most direct
generalization of (the proof of) Theorem \ref{int.t.GCetal} to general Banach
spaces. 

If we combine Theorem \ref{int.t.main} with the result of Carbonaro and
Dragičević  (Theorem \ref{app.t.CD}), we obtain the following improvement 
of Corollary \ref{int.c.CD}.

\begin{cor}\label{int.c.CD-improved}
Let $1<p<\infty$  and let $-A_p$ be the $\Ell{p}$-generator of a symmetric
contraction semigroup. Then  (the injective part of) $A_p$  
has an $R$-bounded
$\Hr^{2}_v(\sect{\omega_p})$-calculus for each admissible function
$v:\R \to [1, \infty)$ such that 
\[
\frac{(1+ |\bfs|)^{\frac{1}{2}}}{v}\, \in \, \Ell{2}(\R).
\]
In particular, $A_p$ has an $R$-bounded $\Hr^{\beta, 2}(\sect{\omega_p})$-calculus for each $\beta > 1$.  
\end{cor}

Parts a) and b) of Theorem \ref{int.t.main} (as well as the corresponding
results about strip-type operators and parts of the underlying function theory) and Corollary
\ref{int.c.CD-improved} have first been established in the second author's
Ph.D. thesis \cite{Pan2019}.

\subsection*{Synopsis}

Having sketched the main achievements of the paper, let us
give a short synopsis. Roughly the paper is divided into four parts:
function theory (Sections \ref{s.adm}--\ref{s.hrs}), operator theory
(Sections \ref{s.fcc}--\ref{s.sec}), Applications (Section
\ref{s.app}) and Appendices (Appendix \ref{s.aux} and \ref{s.alp}). 

\medskip

\subsubsection*{Function Theory}

In Section \ref{s.adm} we introduce (strongly) admissible
functions, following an idea of Strichartz
\cite{Str1971}.  These functions 
shall then replace (viz. generalize) the polynomial weights in describing
$\Ell{2}$-Sobolev spaces. 

In Section \ref{s.sob} we introduce certain function spaces on
vertical strips $\strip{\omega}$, e.g., 
the space of  holomorphic Schwartz functions $\Ha_0(\strip{\omega})$, the
generalized Fourier algebra $\uA_v(\cstrip{\omega})$ and the
generalized Sobolev space $\We_v^2(\strip{\omega})$.
In order to deal with boundary values we employ the classical Paley--Wiener theorem
for the Hardy space $\Har{2}(\strip{\omega})$, which is sufficient for
our purposes. The function theory on sectors is  reduced to 
the strip case via the $\exp/\log$-correspondence and is part of
Section \ref{s.sec}. (This allows us to  work with Fourier 
instead of Mellin transforms.)

Section \ref{s.hrr} is devoted to the spaces $\Hr^2_v(\R)$  of (generalized) Hörmander  functions
on the real line. The main novelties here are that we replace test
functions by general Schwartz functions in the definition and use
admissible functions in contrast to just polynomial weights. 
The latter
forbids using the standard complex interpolation technique  (like
in \cite[Proof of Lemma 6]{CarDra2017}), which is why the 
quite detailed exposition in the preparatory Section \ref{s.sob} is necessary.

The subsequent Section \ref{s.hrs} deals with 
the space $\Hr^2_v(\strip{\omega})$ of
generalized Hör\-man\-der
functions on a vertical strip.  Theorem \ref{hrs.t.indep} is the
central technical result, and the   ``Calderón type reproducing
formula'' \eqref{hrs.eq.CRF}
\[   \vphi f = \int_\R (\tau_t \psi^* \cdot \vphi) (\tau_t \psi \cdot f)
\, \ud{t}
\]
established in its proof, is the main auxiliary result. This formula
needs ``holomorphic cut-off functions'' $\vphi, \psi$, which is the reason why our 
definition of the holomorphic Hörmander spaces differs from the
classical one by virtue of boundary values. (Note 
that there are no non-trivial holomorphic test functions.) This also 
motivates in hindsight the generalization
to Schwartz functions in the definition of $\Hr^2_v(\R)$ in Section \ref{s.hrr}.
At the end of the section, in Theorem \ref{hrs.t.strongadm}, we confine to Hörmander spaces with respect to
strongly admissible functions and show, in particular, that they
coincide with what one would get by using the classical approach via
boundary values.

\subsubsection*{Operator Theory}

In Section \ref{s.fcc} we review the theory of holomorphic functional
calculus for strip type operators and give conditions when an operator
admits an (possibly unbounded) Hörmander calculus (Proposition
\ref{fcc.p.ubdhc}). This generalizes the $0$-sectorial case treated by 
Kriegler and Weis (Remark \ref{fcc.r.KrWe-ubdhc}).  

Section \ref{s.bsc} is devoted to operators with a bounded Sobolev
calculus. Section \ref{s.bhc} then contains the principal  results of
the paper, 
Theorem \ref{bhc.t.meda}, \ref{bhc.t.main} and \ref{bhc.t.main-semi},
about bounded Hörmander calculus on strips. 

In Section \ref{s.sec} we briefly transfer the function theory from strips to
sectors and introduce sectorial operators and their functional
calculus. Finally, we transfer the chief results from the previous
section to sectorial operators, leading to Theorems
\ref{sec.t.meda}--\ref{sec.t.main-semi} (which are summarized in
Theorem \ref{int.t.main} above).

\subsubsection*{Applications}

This short Section \ref{s.app} specializes the results of the previous sections
to generators of $C_0$-groups on Hilbert spaces, sectorial operators
on $\Ell{p}$-spaces, generators of symmetric contraction semigroups
(taking up the discussion in the Introduction) and, as a special
case, the Ornstein--Uhlenbeck operator.

\subsubsection*{Appendices}

In Appendix \ref{s.aux} we recall some results from classical analysis
related to Young's inequality  and important for establishing the
Calderón type reproducing formula for holomorphic Hörmander
functions. In Appendix \ref{s.alp} we review some
``Abstract Littlewood--Paley'' theory, i.e., the Kalton--Weis theory of
abstract square function estimates 
and its connection with bounded $\Ha$-calculus, and related notions from 
Banach  space geometry (spaces of finite cotype, Pisier's contraction
property, $R$-boundedness).


\subsection*{Notation and Terminology}


For a Banach space $X$ the dual space is denoted by $X'$ and the dual
pairing by $\langle \cdot, \cdot \rangle_{X, X'}$. Also, $\BL (X)$ is the Banach algebra of
all bounded operators on $X$, with  $\Id\in \BL (X)$ being  the identity operator. 

A  possibly unbounded linear operator $A$ between two Banach
spaces $X$ and $Y$ is  identified with its graph in $X\oplus Y$. By
$\dom A$, $\ran A$, $\ker A$, we denote the
domain, the range, and the null space of $A$, respectively. If $X=Y$
then  $\spec(A)$ and $\resol(A)$ denote the spectrum and the
resolvent set of $A$, respectively.  Whenever $\lambda \in \resol(A)$,
we write $R(\lambda, A)\coloneqq (\lambda - A)^{-1}$ for the resolvent
operator of $A$ at $\lambda$.

\medskip

We employ the notation 
\[   F(x, f, \dots) \lesssim G(x,f, \dots) \qquad (x\in \cdots.\, f\in
\dots, \dots)
\]
as an abbreviation for ``there is a real number $c\ge 0$ such that 
$F(x, f, \dots) \le c G(x,f, \cdots)$ for all $x\in \dots$, $f\in
\dots$ and $\dots$.'' That means, $c$ may depend on other objects but
not on $x, f, \dots$. 

\medskip

By $\N$ we denote the positive integers, i.e., $\N = \{1, 2,
\ldots\}$, and 
$\R_{>0}\coloneqq (0, \infty)$. For $\omega > 0$ we use the notation 
\[
\strip{\omega} \coloneqq \big\{ z\in \C\, \big|\, |\Im z| < \omega \big\}
\] 
for the strip of width $2\omega$ symmetric about the real  axis, and
$\strip{\omega} := \R$ for $\omega = 0$. Similarly, for $0 < \omega
\le \upi$
\[
\sector{\omega}\coloneqq \exp(\strip{\omega}) = \big\{ z\in \C\setminus \{ 0\}  \, \big|\, |\mathrm{arg}\, z| < \omega  \big\}
\] 
and $\sector{\omega} = (0, \infty)$ for $\omega = 0$.

On any set $M$, $\car$ denotes the function which is constantly equal to
$1$. Also we sometimes write $\bfz$ or $\bfs$ to denote the the
identity mapping on $M$. Usually we use $\bfz$ whenever $M= O$ is a
domain in the complex plane, and $\bfs$ or $\bft$ when $M \subseteq
\R$ is a real interval.

Whenever $f$ is a function defined on a domain $M \subseteq \C$ and
$z\in \C$, we let
\[ \tau_z f := f(\cdot - z)
\]
with domain $M + z$.  Moreover, $f^*$ is the function
\[ f^*(z) := \konj{ f(\konj{z})}\quad (\konj{z} \in M).
\]
If $r\in \R$ is such that $\R + \ui r\subseteq M$, then $f_{|r} :=
f(\bfs + \ui r)$ is a function on $\R$.

For functions $f, \, g: \R \to \C$,
\[
(f\ast g)(t)\coloneqq \int_\R f(t-s)g(s)\, \ud s\qquad (t\in \R)
\] 
is  the usual convolution of $f$ and $g$ (when well-defined). 
The Fourier transform of a  function $f\in \Ell{1}(\R)$ is
\[
(\Fourier f)(t)\coloneqq \widehat{f}(t)\coloneqq \int_\R f(s)\ue^{-\ui st}\, \ud s\qquad (t\in \R)
\] 
and the inverse Fourier transform of $f$ is
\[
(\Fourier^{-1} f)(t)\coloneqq f^\vee(t)\coloneqq \frac{1}{2\pi }\int_\R f(s)\ue^{\ui st}\, \ud s\qquad (t\in \R).
\] 
We shall freely use standard results of elementary Fourier
analysis.

\subsection*{Glossary of Function Spaces}

For a nonempty open set $O\subseteq \C$ we let
\begin{align*}
\Hol (O) & \coloneqq \big\{ f: O \to \C\, \big|\, f \text{  is holomorphic}\big\},\\[0.1cm]
\Ha (O) & \coloneqq \big\{ f\in \Hol (O)\, \big|\,f \text{  is
          bounded}  \big\}\, =\, \Hol (O) \cap \Cb (O),
\end{align*} 
the latter space being equipped with the supremum norm
$\| f\|_{O, \infty}\coloneqq \sup_{z\in O}|f(z)|$. The space of test
functions on  $\R$ is  $\Test(\R)$, and $\calS(\R)$ denotes the space
of Schwartz functions. For $\omega > 0$ and $1 \le p < \infty$, the
Hardy space of order $p$ on $\strip{\omega}$ is
\[  \Har{p}(\strip{\omega}) = \{ f\in \Hol(\strip{\omega}) \suchthat
\norm{f}_{\Har{p}(\strip{\omega})} < \infty\},
\]
where 
\[ \norm{f}_{\Har{p}(\strip{\omega})}^p := 
\sup_{\abs{r} < \omega} \int_\R  \abs{f(s +\ui r)}^p \, \ud{s}. 
\]
The following function spaces are defined in the 
course of the paper (page number in brackets):
\begin{itemize}
\item $\Ell{p}_{v, \omega}(\R),\,  \Ell{p}_v(\R)$: weighted
  $\Ell{p}$-spaces (\pageref{f.weightedLp}) \smallskip
\item $\Ha_0(\strip{\omega})$, $\Ho[\cstrip{\omega}]$,
$\Ho(\sector{\omega})$: 
holomorphic Schwartz functions (\pageref{f.holSchwartz-strip} and \pageref{f.holSchwartz-sector}) \smallskip

\item $\uA(\cstrip{\omega})$:  generalized Fourier algebra (\pageref{f.Fourieralg}) \smallskip

\item $\We^2_v(\strip{\omega})$, $\We^2_v(\sector{\omega})$ generalized
$\Ell{2}$-Sobolev spaces (\pageref{f.Sob-strip} and \pageref{f.Sob-sector}) \smallskip

\item  $\We^{\alpha, 2}(\R)$, 
$\We^{\alpha, 2}(\sector{\omega})$: classical $\Ell{2}$-Sobolev
  spaces  (\pageref{f.Sob-strip-class} and  \pageref{f.Sob-sector-class}) \smallskip

\item $\Har{p}(\sector{\omega})$:   Hardy space of order $p$ on a
  sector (\pageref{f.Hardy-sector}) \smallskip

\item $\Hr^2_v(\R; \psi)$,  $\Hr^2_v(\R)$:  generalized Hörmander
  functions  on $\R$ (\pageref{f.Hrr-psi} and \pageref{f.Hrr}) \smallskip

\item $\Hr^2_v(\strip; \psi)$,  $\Hr^2_v(\strip{\omega})$:  generalized Hörmander
functions on $\strip{\omega}$ (\pageref{f.Hrst-psi} and \pageref{f.Hrst})


\item $\Hr^2_v(\sector{\omega})$, $\Hr^{\beta, 2}(\sector{\omega})$:
  generalized and classical Hörmander
functions  on a sector, \pageref{f.Hrse-class}

\item $\calE(\strip{\theta})$, $\calE[\cstrip{\omega}]$, $\calE(\sector{\theta})$, $\calE[\csect{\omega}]$ elementary
functions on
strips and sectors, \pageref{f.elem-strip} and \pageref{f.elem-sector}.

\end{itemize}

\section{Admissible Weight Functions}\label{s.adm}

Classical fractional Sobolev spaces are defined via the (inverse) Fourier transform:
\begin{equation}\label{adm.eq.FracSob}
\We^{\alpha, 2} (\R) \, =\, \big\{ f\in \Ell{2}(\R)\, \big|\, (1+|\bfs|)^\alpha f^\vee\in \Ell{2} (\R)\big\},
\end{equation} where $\alpha > \tfrac{1}{2}$, see \cite[Definition
6.2.2.]{Gra2009} for example. The first step to a more general class
of Sobolev spaces is to replace the polynomial weights
$(1+|\bfs|)^\alpha$  in (\ref{adm.eq.FracSob}) by more general
functions. This goes back to Strichartz \cite{Str1971}.

\begin{defn}\label{adm.d.weight}
A measurable function $v: \R \to (0, \infty)$ is called \emdf{admissible} if $v$ has the following properties:
\begin{aufziii}
\item $\forall\, s\in \R\colon \quad v(s)\geq 1$;
\item $\dps M_v \coloneqq \sup_{s, t\in \R} \frac{v(s+t)}{v(s)+v(t)} <
  \infty.$
\end{aufziii} 
An admissible function $v$ is called 
\emdf{strongly admissible} if, in addition,
\begin{aufziii}
\item[$\mathrm{3)}$] $\tfrac{1}{v}\in \Ell{2}(\R )$.
\end{aufziii}
\end{defn}

Note that it is immediate from 1) and 2) that 
\[   v(s+t) \le 2 M_v v(s) v(t)  \quad \text{and}\quad 
v(2t) \le 2M_v \, v(t) \qquad (s,t\in \R). 
\]
Let us call a function $v: \R \to \R_{> 0}$ \emdf{quasi-monotonic} if it is decreasing on
$\R_{\le 0}$ and increasing on $\R_{\ge 0}$. 
The following lemma helps to recognize admissible functions.

\begin{lem}\label{adm.l.critsubadd} Let $v:\R \to (0, \infty)$ be
be quasi-monotonic
with $v(0) \ge 1$. Then $v$ is admissible if and only if 
\[ \sup_{t\in \R} \frac{v(2t)}{v(t)} < \infty.
\]
\end{lem}

\begin{proof}
As seen above,  each admissible function has the required property. 
Suppose $c \ge 0$ is such that $v(2t)\le c v(t)$ for all $t\in \R$,
and fix $s, t\in \R$ with $s\leq t$. Then $2s \leq s+t\leq 2t$. Hence, if $s+t \geq 0$, then 
$v(s+t)\, \leq\, v(2t)\, \leq \, cv(t) 
\le \, c\, (v(t) + v(s))$  by monotonicity of $v$ on
$\R_{\ge 0}$. Likewise, if $s+t< 0$, then 
$v(s+t)\, \leq\, v(2s)\, \leq \, cv(s) \le c (v(s) + v(t))$. 
\end{proof}

We collect some properties of admissible functions.

\begin{lem}\label{adm.l.ExGen}\label{adm.l.localbdd}
Let $v,\,  \tilde{v}: \R \to [1, \infty )$ be admissible functions. Then the following statements hold:
\begin{aufzi}
\item $v$ is locally bounded.
 
\item $v^\alpha$ is admissible for each $\alpha \in \R_{\ge 0}$.
\item If both $v$ and $\tilde{v}$ are quasi-monotonic, then 
$v\tilde{v}$ is admissible and quasi-monotonic as well.
\item For all $\theta, \, \tilde{\theta} \in \R_{\ge 0}$ with
  $\theta+\tilde{\theta}= 1$ the function $\theta v +
  \tilde{\theta}\tilde{v}$ is admissible with $M_{\theta v +
    \tilde{\theta}\tilde{v}} \le \max\{ M_v, M_{\tilde{v}}\}$.
\end{aufzi}
\end{lem}

\begin{proof}
a)\ It suffices to show that $v$ is locally bounded at $t=0$, i.e., that there is $\veps > 0$ such that $v$ is bounded on $(-\veps , \veps )$. Choose $R>0$ such that
\[
A_R \coloneqq \big\{ s\in \R\, \big|\, v(s) + v(-s)\leq R \big\}
\] has Lebesgue measure greater than zero. By a theorem of Steinhaus \cite{Str1972}, there is $\veps > 0$ with
\[
(-\veps, \veps)\, \subseteq \, A_R - A_R.
\] 
Now, property {\rm 2)} of Definition \ref{adm.d.weight} yields the claim.

\prfnoi
b) and d) are immediate. For c) use Lemma \ref{adm.l.critsubadd}.
\end{proof}

Let us list some examples.

\begin{exas}\label{adm.ex.ExamplesForWeights}
\begin{aufziii}
\item (Polynomial weights) $\quad$ Let $v\coloneqq 1 + |\bfs|$. Then
  $v$ is admissible with $M_v = 1$. Consequently, by Lemma
  \ref{adm.l.ExGen}.b), $v_\alpha \coloneqq v^{\alpha}$ is
  admissible for each $\alpha \geq 0$, and strongly admissible
for $\alpha > \tfrac{1}{2}$. Note that $v_\alpha$ is even quasi-monotonic.

\item (Mixed weights) $\quad$ Let $v_{\log} \coloneqq \ln (\ue +
  |\bfs|)$. Then $v_{\log}$ is quasi-monotonic.
For each $s\in \R$
\[
\frac{v_{\log} (2s)}{v_{\log} (s)}\, =\, \frac{\ln\big(\frac{\ue}{2} + |s| \big) + \ln 2}{\ln (\ue + |s|)}\, \leq\, 1 + \frac{\ln 2}{\ln (\ue + |s|)}\, \leq\, 1 + \ln 2.
\] Hence, by Lemma \ref{adm.l.critsubadd}, $v_{\log}$ is admissible. By Lemma \ref{adm.l.ExGen},
\[
v_{\alpha, \beta} \coloneqq (1+|\bfs|)^{\alpha} v_{\log}^\beta\, =\, (1+|\bfs|)^{\alpha} \big(\ln (\ue + |\bfs|)\big)^\beta 
\] is a quasi-monotonic admissible function for each $\alpha, \, \beta > 0$. Whenever either $\alpha > \tfrac{1}{2}$ and $\beta \geq
0$, or $\alpha = \tfrac{1}{2}$ and $\beta > 1$, the function
$v_{\alpha, \beta}$ is strongly admissible.
\item (Counterexample)$\quad$ The function
\[
v\coloneqq \exp \big(\sqrt{|\bfs|}\big).
\] is quasi-monotonic.
However,
\[
\frac{v(2s)}{v(s)}\, =\, \exp\Big( (\sqrt{2}-1)\sqrt{|s|}\Big)
\, \xrightarrow{s\to \pm \infty} \, \infty.
\] Hence, $v$ is not admissible.
\end{aufziii}
\end{exas}

The  counterexample from above becomes clear also from the next 
lemma.

\begin{lem}\label{adm.l.polynomialgrowth}
Every admissible function grows at most polynomially. More precisely: Let $v:\R \to [0, \infty )$ be a function which is locally bounded at zero, and suppose that there is a constant $c> 0$ with
\[
v(2s) \, \leq \, c v(s) \qquad (s\in \R).
\] Then there is $\alpha \geq 0$ with 
\[
\sup_{s\in \R} \frac{v(s)}{(1+|s|)^\alpha}\, <\, \infty.
\]
\end{lem}

\begin{proof}
By passing, if necessary, to the function $\max\{ v(s), v(-s)\}$ we
may suppose that $v$ is even.  The hypotheses imply that $v$ is
bounded on $[0,2]$. Choose $\alpha > 0$ large enough so that $c\leq 2^\alpha$. Then
\begin{align*}
\frac{v\big( 2^n 2^r \big)}{\big( 2^n 2^r \big)^\alpha}
\,  \leq\,   c^n  \frac{v\big( 2^r \big)}{2^{\alpha n}2^{\alpha r}}
 \, =\,  \left(\frac{c}{2^\alpha}\right)^n \frac{v\big( 2^{r}\big)}{2^{\alpha r}}
 \, \leq\,   \sup_{1\leq s\leq 2} v(s)
\end{align*} for each $0\leq r\leq 1$ and $n\in \N_0$. As each $s\geq 1$ has a representation $s= 2^{n+r}$ for some $n\in \N_0$ and $r\in [0,1]$, we obtain
\[
\sup_{s\geq 1} \frac{v(s)}{s^\alpha} \, \leq \sup_{1\leq s\leq 2} v(s). 
\] Consequently,
\[
\sup_{s\in \R} \frac{v(s)}{(1+|s|)^\alpha} \, =\, \sup_{s\geq 0} \frac{v(s)}{(1+s)^\alpha} \, \leq\,   \sup_{0 \leq s\leq 2} v(s)\, <\, \infty,
\] which yields the claim.
\end{proof}

Two $\R_{> 0}$-valued functions $v,\tilde{v}$ are called \emdf{equivalent}
if there is $c > 0$ such that 
\[
\frac{1}{c} \tilde{v} \, \leq\, v\, \leq\, c \tilde{v}.
\]
Evidently, if $v,\tilde{v}$ are equivalent and $v$ is 
strongly admissible, then so is $\tilde{v}$.

\begin{lem}\label{adm.l.AWalwaysSmooth}
Let $v: \R \to [1, \infty )$ be an admissible function. Then there is
an admissible function $\tilde{v}\in \uC^\infty (\R)$ equivalent to
$v$. If $v$ is even then $\tilde{v}$ can be
chosen to be even as well. 
\end{lem} 

\begin{proof}
Let $\eta \in \Test(\R)$ be a positive, even function with 
\beq\label{adm.eq.AWalwaysSmooth}
\int_\R \eta (r)\, \ud r\, =\, 1,
\eeq
and such that $\eta\vert_{[0, \infty )}$ is decreasing. Then
$\eta (2t)\, \leq\, \eta (t)$ for all $t\in \R$.
Set $\tilde{v}\coloneqq v\ast \eta$. Then $1\leq \tilde{v}\in
\uC^{\infty} (\R)$, and one has, for all $s,\, t \in \R$,
\begin{align*}
\tilde{v}(s+t) =&\, \int_\R v(s+t-r)\eta (r)\, \ud r
\leq\, M_v \Bigg( \int_\R v(s-\tfrac{r}{2})\eta (r)\, \ud r + \int_\R v(t-\tfrac{r}{2})\eta (r)\, \ud r \Bigg)\\
=&\,  2 M_v \Bigg(\int_\R v(s-r)\eta (2r)\, \ud r + \int_\R v(t-r)\eta (2r)\, \ud r\Bigg)
   \leq\, 2M_v \bigl( \tilde{v} (s) + \tilde{v}(t) \bigr)
\end{align*} 
Hence, $\tilde{v}$ is an admissible function. As $v$ is
locally bounded, choose $C> 0$ such that $v(t) \le C$ for $t\in
\supp(\eta)$. Then, by \eqref{adm.eq.AWalwaysSmooth},
\begin{align*}
v(s) =&\, \int_\R v(s)\eta (r)\, \ud r\, \leq\,  2M_v \int_\R
        v(r)v(s-r)\eta (r)\, \ud r\, \leq\, 2 M_v C \,\tilde{v}(s),
\end{align*} and, as $\eta$ is even,
\begin{align*}
\tilde{v}(s) =&\, \int_\R v(s-r)\eta (r)\, \ud r\, =\, \int_\R v(s+r) \eta (r)\, \ud r\\
\leq&\, 2M_v \int_\R v(r)\eta (r)\, \ud r\,  v(s)\, \leq\, 
2M_v C\, v(s).
\end{align*} 
If $v$ is even then so is $\tilde{v}$, hence the claim is proved.
\end{proof}

\begin{rem}
We shall see later that equivalent admissible functions lead to the same
function spaces with equivalent norms 
(see, e.g., Propositions \ref{adm.p.Lpom}, \ref{sob.p.We}. \ref{hrs.p.Hr}).
Thanks to  Lemma \ref{adm.l.AWalwaysSmooth}, any appearing
admissible function $v$  can always be assumed to be smooth. 

E.g., instead of the non-smooth polynomial function $v_\alpha(t)
= (1 + \abs{t})^\alpha$, $\alpha \ge 0$, we may work with the
equivalent smooth function
$(1+ t^2)^\frac{\alpha}{2}$.
\end{rem}

Let $v$ be admissible, $\omega \in \R_{\ge 0}$  and $p \in [1,\infty]$. 
Then  we may form the space \label{f.weightedLp}
\begin{align*}
 \Ell{p}_{v,\omega}(\R) & := \{ f\in \Ell{p}(\R) \suchthat 
  \ue^{\omega \abs{\bfs}} vf \in \Ell{p}(\R)\},
\end{align*} 
where  $\Ell{p}_v(\R) := \Ell{p}_{v,0}(\R)$, for simplicity. Note that 
$\ue^{\omega \abs{s}}$ and $\cosh(\omega \bfs)$ are equivalent, and
hence to describe the space $\Ell{p}_{v, \omega}(\R)$ we may replace
the former by the latter whenever it is convenient, cf. part b)
below. 

Let us collect some more or less well-known or obvious facts.

\begin{prop}\label{adm.p.Lpom}
Let $v$ be admissible, $\omega\in \R_{\ge 0}$  and $p \ge 1$. Then the following assertions
hold:
\begin{aufzi}
\item $\Ell{p}_{v,\omega}(\R)$ is a Banach space with respect to the norm
$\norm{f}_{\Ell{p}_{v,\omega}} = \norm{v \ue^{\omega \abs{\bfs}}f}_{\Ell{p}}$, continuously embedded
into $\Ell{p}(\R)$ and isomorphic to $\Ell{p}(\R)$ via the
mapping $f \mapsto v \ue^{\omega\abs{\bfs}} f$. 

\item If $\tilde{v}$ is admissible with  $\tilde{v}\lesssim v$, 
$\Ell{p}_{v,\omega}(\R) \subseteq \Ell{p}_{\tilde{v},\omega}(\R)$, continuously.
In
particular, if $\tilde{v}$ and $v$ are equivalent, then 
$\Ell{p}_{v,\omega}(\R) = \Ell{p}_{\tilde{v},\omega}(\R)$ with equivalent norms.

\item If $0 \le \tilde{\omega}< \omega$ and $\tilde{v}$ is admissible,
then $\Ell{p}_{v, \omega}(\R) \subseteq \Ell{1}_{\tilde{v}, \tilde{\omega}}(\R)
\cap \Ell{p}_{\tilde{v},\tilde{\omega}}(\R)$, continuously.

\item If $p < \infty$ then $\Test(\R)$ is dense in $\Ell{p}_{v,\omega}(\R)$.

\item If $p < \infty$ then translation  $(\tau_t)_{t\in \R}$ is a  strongly continuous
  group  on $\Ell{p}_{v,\omega}(\R)$ with $\norm{\tau_t} \le 2 M_v \ue^{\omega\abs{t}}v(t)$ for
  $t\in \R$.

\item 
If  $f \in \Ell{1}_{v,\omega}(\R)$ and $g\in \Ell{p}_{v, \omega}(\R)$
then $f \ast g \in
  \Ell{p}_{v,\omega}(\R)$  with
\[   \norm{f \ast g}_{\Ell{p}_{v,\omega}(\R)} \le 2M_v \norm{f}_{\Ell{1}_{v,\omega}(\R)}
\norm{g}_{\Ell{p}_{v,\omega}(\R)}.
\]

\prfnoi
If $f \in \Ell{p'}_{v,\omega}(\R)$ and $g\in \Ell{p}_{v,
  \omega}(\R)$ then $f \ast g\in \Ce(\R) \cap \Ell{\infty}_{v,
  \omega}(\R)$ 
with 
\[    \norm{f\ast g}_{\Ell{\infty}_{v, \omega}(\R)} \le 2M_v  
\norm{f}_{\Ell{p'}_{v,\omega}(\R)} \norm{g}_{\Ell{p}_{v,\omega}(\R)}.
\]

\item If $v$ is strongly admissible, then 
$\Ell{2}_{v,\omega}(\R)$
is (after scaling the norm) a Banach algebra with respect to convolution,
continuously embedded into $\Ell{1}_{\car,\omega}(\R)$.

\item Let  $\tilde{v}$ be admissible and quasi-monotonic with 
$v \lesssim \tilde{v}$, and let $\eta \in \Ell{1}_{\tilde{v},
  \omega}(\R)$. Define $\eta_n(\bfs) := n
  \eta(n\bfs)$ for $n \in \N$. Then 
$\sup_{n\in \N}
 \norm{\eta_n}_{\Ell{1}_{v, \omega}(\R)} < \infty$ and 
for each $f\in \Ell{p}_{v, \omega}(\R)$
one
  has  $\dps (\ue^{\ui t \bfs} \eta_n) \ast f \stackrel{n \to \infty}{\longrightarrow}  \fourier{\eta}(0) \cdot
  f$ in $\Ell{p}_{v,\omega}(\R)$ uniformly  in $t$ from compact
  subsets of $\R$. 
\end{aufzi}
\end{prop}

\begin{proof}
a) and b) are straightforward. c) follows from the
computation
\[   \tilde{v} \ue^{\tilde{\omega}\abs{\bfs}} f
=   \frac{\tilde{v}}{v} \ue^{-( \omega - \tilde{\omega})\abs{\bfs}} 
\, v \ue^{\omega\abs{\bfs}} f
\]
and the fact that $\tilde{v}$ is dominated by some $v_\alpha$.  

\prfnoi
d) follows from a) and b)
and the fact that we may replace $v$ by a smooth equivalent
and $\ue^{\omega\abs{\bfs}}$ by the (smooth) function $\cosh(\omega\bfs)$.

\prfnoi
e) For $f\in \Ell{p}_{v,\omega}(\R)$ one has 
\[ \norm{\tau_t f}_{\Ell{p}_{v,\omega}} = \norm{v \ue^{\omega{\bfs}}\tau_t f}_{\Ell{p}}
= \norm{v(\bfs + t) \ue^{\omega\abs{\bfs + t}} f}_{\Ell{p}} \le 2 M_v
v(t) \ue^{\omega\abs{t}}
\norm{f}_{\Ell{p}_{v,\omega}}
\]
as $v(x + t) \le 2 M_v v(x) v(t)$ for all $x,t\in \R$. Hence, each
translation operator is bounded and the translation group is
locally bounded. The strong continuity now follows from a density
argument
employing c).

\prfnoi
f) Define  $F := v \ue^{\omega
  \abs{\bfs}}\abs{f}$ and $G:= v \ue^{\omega
  \abs{\bfs}}\abs{g}$.
\begin{align*}
v(s) & \ue^{\omega \abs{s}} |(f \ast g)(s)|
\leq  v(s) \int_\R \ue^{\omega\abs{s-t}}|f(s{-}t)| \, \ue^{\omega
  \abs{t}}|g(t)|\, \ud t \\
& \le 2M_v \int_\R v(s-t)\ue^{\omega\abs{s-t}}|f(s{-}t)| \, v(t)\ue^{\omega
  \abs{t}}|g(t)|\, \ud t 
= 2M_v (F\ast G)(s)
\end{align*}
In the first case, $F\in \Ell{1}(\R)$ and $G\in \Ell{p}(\R)$, hence 
$F\ast G \in \Ell{p}(\R)$ again. This yields $f \ast g \in \Ell{p}_{v,
  \omega}(\R)$ and the claimed norm estimate. The same argument
applies
in the second case, where $F\in \Ell{p'}(\R)$ and $G \in
\Ell{p}(\R)$. Note that here $f\ast g \in \BUC(\R)$ as $f \in \Ell{p'}(\R)$
and $g\in \Ell{p}(\R)$.

\prfnoi
g) Let $f,g \in \Ell{2}_{v,\omega}(\R)$. Since $\frac{1}{v} \in
\Ell{2}(\R)$, one has 
\[     \ue^{\omega \abs{\bfs}} f = \frac{1}{v} (v \ue^{\omega
  \abs{\bfs}} f) \in \Ell{2}(\R) \cdot \Ell{2}(\R) \subseteq
\Ell{1}(\R)
\]
with (hence) 
\[    \norm{f}_{\Ell{1}_{\car,\omega}(\R)} \le 
\norm{\tfrac{1}{v}}_{\Ell{2}(\R)} \norm{f}_{\Ell{2}_{v,\omega}(\R)}.
\]
This proves the continuous embedding $\Ell{2}_{v,\omega}(\R) \subseteq
\Ell{1}_{\car,\omega}(\R)$.  

As $f,g \in \Ell{2}(\R)$, the convolution $f\ast g$ is defined
pointwise and a function in  $\Co(\R)$. Set $F := \ue^{\omega
  \abs{\bfs}}\abs{f}$ and $G:= \ue^{\omega
  \abs{\bfs}}\abs{g}$. Then 
\begin{align*}
& v(s) \ue^{\omega \abs{s}} |(f \ast g)(s)|
\leq  v(s) \int_\R \ue^{\omega\abs{s-t}}|f(s{-}t)| \, \ue^{\omega \abs{t}}|g(t)|\, \ud t \\
&\leq M_v \int_\R \big( v(s{-}t) + v(t) \big) F(s{-}t) \, G(t)\, \ud t \\
&=\, M_v (vF \ast G) (s) + M_v (F\ast vG)(s).
\end{align*} 
By the first part, $F, G \in \Ell{1}(\R)$. It follows that 
\begin{align*}
 \norm{f\ast g}_{\Ell{2}_{v,\omega}(\R)}
&  \leq\, M_v \big( \|vF \ast G\|_{\Ell{2}(\R)} + \| F\ast vG\|_{\Ell{2}(\R)} \big)\\[0.1cm]
  & \leq\, M_v\big( \|v F\|_{\Ell{2}(\R)} \, \|G\|_{\Ell{1}(\R)} + \|F\|_{\Ell{1}(\R)} \, \|vG\|_{\Ell{2}(\R)} \big)\\[0.1cm]
  &=\, M_v \big( \| f\|_{\Ell{2}_{v,\omega}(\R)}\,  \| g\|_{\Ell{1}_{\car,\omega}(\R)} + \|f\|_{\Ell{1}_{\car,\omega}(\R)}\,  \|g\|_{\Ell{2}_{v,\omega}(\R)} \big)\\[0.1cm]
  &\leq \, \bigl( 2M_v \| \tfrac{1}{v}\|_{\Ell{2}(\R)} \bigr)\,  
\| f\|_{\Ell{2}_{v,\omega}(\R)}\,  \| g\|_{\Ell{2}_{v,\omega}(\R)}.
\end{align*} 

\prfnoi
h)\ Observe that 
\[   v(s) \ue^{\omega \abs{s}}\eta_n(s) = \frac{v(s)}{\tilde{v}(s)}
\frac{\tilde{v}(s)}{\tilde{v}(ns)} \frac{\ue^{\omega
    \abs{s}}}{\ue^{n\omega \abs{s}}} \cdot  n \cdot ( \tilde{v} \ue^{\omega \abs{\cdot}}\eta)(ns) \qquad (s\in \R).
\]
This implies $\norm{\eta_n}_{\Ell{1}_{v, \omega}(\R)} \le
\norm{v/\tilde{v}}_\infty \norm{\eta}_{\Ell{1}_{\tilde{v},
    \omega}(\R)}$
for each $n\in \N$.

For the second assertion fix $f\in \Ell{p}_{v, \omega}(\R)$. It
suffices to  show that 
\[  (\ue^{\ui t \bfs} \eta_n) \ast f  - \fourier{\eta}(\tfrac{-t}{n})
f= 
\int_\R \ue^{\ui t s} \eta_n(s) (\tau_s f- f)\, \ud{s}
\]
converges to $0$ in $\Ell{p}_{v,\omega}(\R)$ as $n\to \infty$
uniformly in $t\in \R$. Taking norms leads to the estimate
\[ \int_\R \abs{\eta_n(s)} \norm{\tau_s f- f}_{\Ell{p}_{v, \omega}(\R)}\, \ud{s}
\lesssim \sup_{\abs{s}\le \delta} \norm{\tau_s f- f}_{\Ell{p}_{v, \omega}(\R)} 
+ \int_{\abs{s}\ge n\delta} \abs{\eta(s)} \norm{\tau_{\frac{s}{n}}f - f}_{\Ell{p}_{v,\omega}(\R)}
\]
for each $\delta > 0$. For fixed $\delta > 0$ the second summand
converges to $0$ as $n \to \infty$. This follows from 
\[ \norm{\tau_{s/n}f - f}_{\Ell{p}_{v,\omega}(\R)} \lesssim  \ue^{\omega
  \abs{\tfrac{s}{n}}}v(\tfrac{s}{n}) 
\lesssim \ue^{\omega \abs{\tfrac{s}{n}}} \tilde{v}(\tfrac{s}{n}) 
\le
\ue^{\omega \abs{s}} \tilde{v}(s)\qquad (s\in \R, n\in \N)
\]
and the fact that $\eta \tilde{v} \ue^{\omega \abs{\bfs}} \in \Ell{1}(\R)$. By e), the first summand
converges to $0$ as $\delta \to 0$. This proves the claim.
\end{proof}

\begin{rems}\label{adm.r.Lpom}
1.) The second part of Assertion f) can be strengthened in case $1< p < \infty$. Indeed,
one then has (in the proof) $F \ast G \in \Co(\R)$, and hence
$v \ue^{\omega \abs{\bfs}} (f \ast g) \in \Co(\R)$ whenever
$v$ is continuous.

\prfnoi
2.) Assertion g) is essentially due to Strichartz. It is the only point
where the subadditivity of $v$ plays a role. For the convolution
results only the submultiplicativity is needed. 
\end{rems}

\section{Generalized Sobolev Functions on Strips}\label{s.sob}

For this section, we fix an admissible function $v: \R \to \R_{\ge 1}$
and a number $\omega \in \R_{\ge 0}$. The \emdf{horizontal strip}
of \emdf{height} $\omega\ge 0$ is
\[ \strip{\omega} := \begin{cases}   \{ z\in \C \suchthat \abs{\im z} <
  \omega\} & \text{if $\omega > 0$},\\
\qquad\R & \text{if $\omega = 0$}.
\end{cases}
\]
Suppose, for the time being, $\omega > 0$.
Then  $\Hol(\strip{\omega})$ is  the space
of all holomorphic functions on $\strip{\omega}$. 
If convenient, we  identify 
a holomorphic function on some strip with its restriction to
a smaller strip, so that we have
\[  0 < \omega'< \omega \quad \dann\quad \Hol(\strip{\omega})
\subseteq \Hol(\strip{\omega'}) \subseteq \Ce^\infty(\R).
\]
For $f\in \Hol(\strip{\omega})$ and $\abs{y} < \omega$ we denote by 
\[  f_{|y} : \R \to \C,\qquad f_{|y}(x) := f(x + \ui y)
\]
the restriction of $f$ to the ordinate $y$. The space of 
bounded holomorphic functions on $\strip{\omega}$ is
\[ \Ha(\strip{\omega}) := \{ f\in \Hol(\strip{\omega}) \suchthat 
 \norm{f}_{\infty, \strip{\omega}} := \sup_{z\in \strip{\omega}}
 \abs{f(z)} < \infty\}.
\]
Besides this space, we also shall consider the space
\[  \Ho(\strip{\omega}) := \{ f\in \Hol(\strip{\omega}) \suchthat
\forall \alpha > 0 : \abs{f(z)} \lesssim \frac{1}{(1 +  \abs{\re z})^\alpha}\}
\]
of \emdf{holomorphic Schwartz functions}\label{f.holSchwartz-strip}
and, for $\omega \ge 0$, 
\[ \Ho[\cstrip{\omega}] := \bigcup_{\theta > \omega}
\Ho(\strip{\theta}).
\]
Note that by a classical theorem of Paley
and Wiener \cite[Theorem 7.22]{RudinFA}, 
$\fourier{\eta}
\in \Ho(\strip{\omega})$ for each test function
$\eta \in \Test(\R)$ and $\omega > 0$. 
The following lemma tells that the name 
``holomorphic Schwartz function'' is justified.

\begin{lem}\label{hrr.l.Ho}
Let $\omega \ge 0$ and 
$f\in \Ha_0[\cstrip{\omega}]$. Then $f'\in \Ha_0[\cstrip{\omega}]$
as well. More precisely, if $\alpha > 0$ and 
$f$ is holomorphic with 
$\abs{f} \lesssim (1 + \abs{\re \bfz})^{-\alpha}$
on $\strip{\theta}$, for some $\theta > 0$,  then $\abs{f'} \lesssim (\abs{1 + \re \bfz})^{-\alpha}$
on each smaller strip. 

In particular, each function $f\in \Ho[\R]$ restricts to a Schwartz function  $\R$. 
\end{lem}

\begin{proof}
Suppose that $\theta > 0$ and 
$f$ is holomorphic on $\strip{\theta}$ with $\abs{f}\lesssim
(1+ \abs{\re \bfz})^{-\alpha}$. 
Fix $0 < \veps < \theta$ and consider $z\in
\strip{\theta-\veps}$. Then
\[   \abs{f'(z)} \le \frac{1}{2\upi} \int_{\abs{w-z}= \veps}
\frac{\abs{f(w)}}{\abs{w-z}^2}\, \abs{\ud{w}}
\lesssim \sup_{\abs{t}\le \veps} (1 + \abs{\re z - t})^{-\alpha} \le
\frac{(1 + \abs{\veps})^\alpha}{(1 + \abs{\re z})^\alpha}
\]
as claimed. 
\end{proof}

\vanish{
\[  \Ho(\strip{\omega}) := \{ f\in \Hol(\strip{\omega}) \suchthat
\exists\, c,d > 0 : \abs{f(\bfz)} \le c \ue^{-d \abs{\re \bfz}}\}
\]
and, for $\omega \ge 0$, 
\[ \Ho[\cstrip{\omega}] := \bigcup_{\theta > \omega}
\Ho(\strip{\theta}).
\]
}

We now return to the more general setting with $\omega \ge 0$, where
$\omega = 0$ is also allowed.  
In the remainder of this section we shall be concerned with two spaces
of functions on $\strip{\omega}$, $\uA_v(\strip{\omega})$ and
$\We^2_v(\strip{\omega})$. The former is a generalization of the
classical Fourier algebra on $\R$ and plays an auxiliary role. The
latter is a generalization of the classical fractional
$\Ell{2}$-Sobolev space and the main object of interest.  We start
with the former.

\medskip
The \emdf{generalized Fourier algebra} on $\cstrip{\omega}$ with \label{f.Fourieralg}
respect to $v$ is
\[    \uA_v(\cstrip{\omega}) := \{ \fourier{g} \suchthat 
g\in \Ell{1}_{v, \omega}(\R)\}
\]
endowed with the norm
\[   \norm{\fourier{g}}_{\uA_v(\cstrip{\omega})} := 
\norm{g}_{\Ell{1}_{v, \omega}(\R)} \qquad (g\in  \Ell{1}_{v,
  \omega}(\R)).
\]
If $v = \car$, we simply write $\uA(\cstrip{\omega})$. Note that if,
in addition, $\omega =0$, then $\uA(\cstrip{0}) = \uA(\R) =
\Fourier(\Ell{1}(\R))$ is the classical Fourier algebra. (This
justifies our terminology.)

Recall that if $g \in \Ell{1}_{v, \omega}(\R)$ 
then its Fourier transform $\fourier{g}$ is given by 
\[ \fourier{g}(z) = \int_\R \ue^{-\ui s z} g(s)\, \ud{s} = 
\int_\R \ue^{-\ui s z} \frac{\ue^{-\omega\abs{s}}}{v(s)} \,
(v(s)\ue^{\omega\abs{s}} g(s) )\, \ud{s} \qquad (z\in \R).
\]    
If $\omega = 0$, then $\fourier{g}\in \Co(\R)$
(Riemann-Lebesgue). If   $\omega > 0$ then the formula above 
is meaningful for $\abs{\im z} \le \omega$ and hence $\fourier{g}$
can be considered to be an element of $\Hol(\strip{\omega}) \cap 
\Cb(\cstrip{\omega})$. Moreover, one then has
\[    (\ue^{y \bfs} g)^\wedge = (\fourier{g})_{|y} \qquad (\abs{y} \le
\omega).
\]
The following proposition collects the relevant properties, some of
them rather obvious.

\begin{prop}\label{sob.p.A}
Let $v, \tilde{v}: \R \to \R_{\ge 1}$ be admissible functions
and $\omega, \tilde{\omega} \ge 0$. 
\begin{aufzi}

\item The space $\uA_v(\cstrip{\omega})$ is a separable Banach space
isometrically isomorphic to $\Ell{1}_{v, \omega}(\R)$ and  to $\Ell{1}(\R)$.

\item If  $\tilde{v}\lesssim
  v$, then $\uA_{v}(\cstrip{\omega}) \subseteq
  \uA_{\tilde{v}}(\cstrip{\omega})$, continuously. In particular, 
if $v, \tilde{v}$ are equivalent, then 
$\uA_v(\cstrip{\omega}) = \uA_{\tilde{v}}(\cstrip{\omega})$ with equivalent norms.

\item Restriction yields a continuous embedding 
\[ \uA_v(\cstrip{\omega}) \subseteq \uA_{\tilde{v}}(\cstrip{\tilde{\omega}})\qquad
(0 \le \tilde{\omega} < \omega). 
\]
\item The space $\Fourier( \Test(\R))$ is dense in
  $\uA_v(\cstrip{\omega})$.

\item One has $\uA_v(\cstrip{\omega}) \subseteq \Co(\cstrip{\omega})$,
continuously.

\item With respect to pointwise multiplication,
  $\uA_v(\cstrip{\omega})$
is a Banach algebra (after a scaling of the norm). 

\vanish{
\item \red{For $\omega > 0$, \quad $\dps
\uA_v(\cstrip{\omega}) = \{ f\in \Ha(\strip{\omega}) \cap 
\Cb(\cstrip{\omega}) \suchthat f_{|\pm \omega} \in \uA_v(\R)\}$
with equivalence of norms
\[  \norm{f}_{\uA_v(\cstrip{\omega})} \eqsim 
\norm{f_{|\omega}}_{\uA_v(\R)} + \norm{f_{|-\omega}}_{\uA_v(\R)}.
\]
}

\item \red{If  $N\in \N$ and $v = (1 + \abs{\bfs})^N$, then
\begin{align*}
 \uA_v(\R) & = \{ f \in \Ce^k(\R) \suchthat f, f', \dots, f^{(N)} \in
\uA(\R)\} \quad \text{and} \\
\uA_v(\cstrip{\omega}) & = 
\{ f \in \Hol(\strip{\omega}) \suchthat f,
                         f', \dots, f^{(N)} \in
                         \uA(\cstrip{\omega})\}\quad \text{for $\omega > 0$.}
\end{align*}
}

}

\item The translation group 
acts strongly continuously and isometrically on $\uA_v(\cstrip{\omega})$. 

\end{aufzi}
\end{prop}

\begin{proof}
a)---d) follow, by construction of $\uA_v(\cstrip{\omega})$, directly  from 
the respective assertions of Proposition \ref{adm.p.Lpom}.

\prfnoi
e)\ By construction and a trivial estimate, $\uA_v(\cstrip{\omega})
\subseteq
\ell^\infty(\cstrip{\omega})$ continuously. By d) it suffices to prove
that $\Fourier(\Test(\R)) \subseteq \Co(\cstrip{\omega})$. But this is
true, since  $\Fourier(\Test(\R)) \subseteq 
\Ho[\cstrip{\omega}]$ as remarked above. 

\prfnoi
f)\ This follows from Proposition \ref{adm.p.Lpom}.f with  $p=1$.

\vanish{
\prfnoi
g)\   

\prfnoi
h)\ 
}

\prfnoi
g)\ follows  readily from  the identity
$\tau_t f = (\ue^{\ui t \bfs} g)^\wedge$ for $t\in \R$
and $f = \fourier{g}  \in \uA_v(\cstrip{\omega})$, $g\in \Ell{1}_{v,\omega}(\R)$.
\end{proof}

\bigskip

\noindent
We now turn to the space of principal interest, which is the
\emdf{generalized Sobolev space} \label{f.Sob-strip}
\[ \We_v^2(\strip{\omega}) := \{ \fourier{g} \suchthat g \in
\Ell{2}_{v,\omega}(\R)\}
\]
and endow it with the norm
\[ \norm{\fourier{g}}_{\We^2_v(\strip{\omega})} :=
\norm{g}_{\Ell{2}_{v, \omega}(\R)} \qquad (g\in \Ell{2}_{v,\omega}(\R)).
\]
By definition of the Fourier transform, an element 
$f\in \We^2_v(\strip{\omega})$ is---a priori---just a tempered
distribution on $\R$. However, since $\Ell{2}_{v, \omega}(\R) \subseteq
\Ell{2}(\R)$,
Plancherel's theorem tells that $f\in \Ell{2}(\R)$. Furthermore, 
for $\omega > 0$ the next result shows that $f$ may be considered to be a
holomorphic function on $\strip{\omega}$.

\begin{prop}\label{sob.p.We}
Let $v, \tilde{v}: \R \to \R_{\ge 1}$ be admissible functions
and $\omega, \tilde{\omega} \ge 0$. 
\begin{aufzi}

\item The space $\We_v^2(\strip{\omega})$ is a separable 
Hilbert space. If $v= \car$ and $\omega=0$ then $\We^2_v(\strip{\omega}) =
    \We^2_{\car}(\strip{0}) = \Ell{2}(\R)$.


\item If $\tilde{v}\lesssim
  v$ then $\We_{v}^2(\strip{\omega}) \subseteq
  \We_{\tilde{v}}^2(\strip{\omega})$ continuously.  In particular, 
if $v, \tilde{v}$ are equivalent, then 
$\We_v^2(\strip{\omega}) = \We_{\tilde{v}}^2(\strip{\omega})$ with equivalent norms.

\item Restriction yields continuous embeddings  
\[ \We^2_{v}(\strip{\omega}) \subseteq \uA_{\tilde{v}}(\strip{\tilde{\omega}}) 
\cap \We^2_{\tilde{v}}(\strip{\tilde{\omega}}) \subseteq 
\Co(\cstrip{\tilde{\omega}})  \qquad(0\le \tilde{\omega} < \omega).
\]
In particular, if $\omega > 0$ then
$\We^2_v(\strip{\omega})$ embeds continuously into 
$\Hol(\strip{\omega})$ endowed with the topology of compact
  convergence through
\beq\label{sob.eq.We-rep}
   f(z) = \int_\R \ue^{-\ui z s}\,f^\vee(s)\, \ud{s} \qquad 
(\abs{\im z} < \omega)
\eeq
for $f\in \We^2_v(\strip{\omega})$. 

\item If  $\frac{\tilde{v}}{v} \in
  \Ell{2}(\R)$, then $\We^2_{v}(\strip{\omega}) \subseteq
  \uA_{\tilde{v}}(\cstrip{\omega})$ continuously.  In particular, 
if  $v$ is strongly admissible, then 
\[ \We^2_v(\strip{\omega}) \subseteq \uA(\cstrip{\omega}) 
\subseteq  
\Co(\cstrip{\omega})
\]
continuously. Moreover, still with $v$ being strongly
admissible, 
$\We^2_v(\strip{\omega})$ is (after a scaling of
norm)  a Banach algebra with respect to pointwise multiplication.

\item For $v = v_\alpha = (1 + |\bfs|)^\alpha$ and $\omega = 0$ the space
$\We_{v_\alpha}^2(\strip{\omega})$ coincides with
the classical fractional Sobolev space 
$\We^{\alpha, 2}(\R)$. \label{f.Sob-strip-class}

In particular, for $\alpha = N \in \N$, the space $\We^2_N(\R)$
coincides with the classical $\Ell{2}$-Sobolev space of all
$\Ell{2}(\R)$-functions with (distributional) derivatives
up to order $N$ being in $\Ell{2}$.



\item Pointwise multiplication is a bounded bilinear mapping
\[  \uA_v(\cstrip{\omega}) \times \We^2_v(\strip{\omega}) \to 
\We^2_v(\strip{\omega})
\]

\item The translation group $(\tau_t)_{t\in \R}$ 
is a strongly continuous group of isometries on
$\We_v^2(\strip{\omega})$. 

\end{aufzi}
\end{prop}

\begin{proof}
a)\ is clear, b) follows  from Proposition \ref{adm.p.Lpom}.b, and 
c) follows from Proposition \ref{adm.p.Lpom}.c and Proposition \ref{sob.p.A}.e.

\prfnoi
d) If $\frac{\tilde{v}}{v} \in \Ell{2}(\R)$ then 
\[  \Ell{2}_{v, \omega}(\R) = \frac{\tilde{v}}{v}
\frac{\ue^{-\omega\abs{\bfs}}}{\tilde{v}} \Ell{2}(\R) \subseteq
\frac{\ue^{-\omega\abs{\bfs}}}{\tilde{v}} \Ell{1}(\R) =
\Ell{1}_{\tilde{v},\omega}(\R).
\]
This proves the first assertion. Suppose that $v$ is strongly
admissible. Then one can take $\tilde{v} := \car$ in the first
assertion and obtains the second. The third follows directly
from Proposition \ref{adm.p.Lpom}.g.

\prfnoi
e)\ is classical.

\prfnoi
f) follows from Proposition \ref{adm.p.Lpom}, f) (with
$p=2$) and g) follows from the identity $(\tau_t f)^\vee = \ue^{\ui t
  \bf s} f^\vee$ for $f\in \We^2_v(\strip{\omega})$ and $t\in \R$. 
\end{proof}

Fix  $\omega > 0$. 
One (other)  classical theorem of Paley and Wiener 
\cite[p.188/198]{Katznelson2004} 
states that 
\[  
\We^2_\car(\strip{\omega}) = \Har{2}(\strip{\omega}),
\]
where the latter is the \emdf{Hardy space} 
of all functions $f\in \Hol(\strip{\omega})$ such that 
\[  \norm{f}_{\Har{2}}^2 := \sup_{\abs{y} < \omega}
\norm{f_{|y}}_{\Ell{2}(\R)}^2 
= \sup_{\abs{y}< \omega} \int_{\R} \abs{f(x + \ui y)}^2\, \ud{x} <
\infty
\]
with equivalence of norms
\[  \norm{f}_{\We^2_\car(\strip{\omega})}  \eqsim 
\norm{f}_{\Har{2}(\strip{\omega})}
\]
From \eqref{sob.eq.We-rep} it follows readily that 
\[       f_{|y} = ( \ue^{y \bfs} f^\vee )^\wedge \qquad (\abs{y} <
\omega),
\]
and, in particular,  $f^\vee = f_{|0}^\vee$. The two $\Ell{2}$-functions 
\[  f_{|\pm\omega} \coloneqq ( \ue^{\pm \omega \bfs} f^\vee )^\wedge
\in \Ell{2}(\R)
\]
are the \emdf{boundary values} of $f$ on $\rand\strip{\omega}$. 
It is easily seen that 
the   mapping
\[   [-\omega, \omega] \to \Ell{2}(\R),\quad y \mapsto f_{|y}
\]
is continuous. In particular,
\[    \lim_{r \nearrow \omega} f_{|r} = f_{|\omega}\quad
\text{and}\quad \lim_{r \searrow {-}\omega} f_{|r} = f_{|_\omega}
\]
in $\Ell{2}(\R)$ (which justifies the term ``boundary value'').

Employing these concepts we arrive at the following characterization
of elements of $\We^2_v(\strip{\omega})$.

\begin{cor}\label{sob.c.We-char}
Let $\omega > 0$ and $v: \R \to \R_{\ge 1}$ admissible. 
Then for a holomorphic function $f\in \Hol(\strip{\omega})$ the
following statements are equivalent:
\begin{aufzii}
\item $f\in \We^2_v(\strip{\omega})$;
\item $f\in \Har{2}(\strip{\omega})$ and $f_{|\pm \omega} \in
  \We^2_v(\R)$;
\item $f_{|0}\in \Ell{2}(\R)$ and $f_{|0}^{\,\,\vee} \in
  \Ell{2}_{v,\omega}(\R)$;
\end{aufzii}
Moreover, there is equivalence of norms
\[    \norm{f_{|\omega}}_{\We^2_v(\R)}^2 +
  \norm{f_{|-\omega}}_{\We^2_v(\R)}^2  \,\, \eqsim \,\,
\norm{f}_{\We^2_v(\strip{\omega})}^2.
\]
\end{cor}

\begin{proof}
(i)$\dann$(ii):  Let $f\in \We^2_v(\strip{\omega})$. Since $\car \le
v$ we have $f\in \We^2_\car(\strip{\omega}) =
\Har{2}(\strip{\omega})$. By definition,   
\[ v (f_{|\pm \omega})^{\,\vee} 
=  v \ue^{\pm \omega \bfs } f^\vee 
=  \frac{\ue^{\pm \omega \bfs}}{\ue^{\omega \abs{\bfs}}} ( v
\ue^{\omega\abs{\bfs}} f^\vee) \in \Ell{\infty}(\R) \cdot \Ell{2}(\R)
\subseteq \Ell{2}(\R). 
\]
This concludes the proof of (ii) and shows
$\norm{f_{|\pm \omega}}_{\We^2_v(\R)} \lesssim \norm{f}_{\We^2_v(\strip{\omega})}$.

\prfnoi
(ii)$\dann$(iii): Suppose $f\in \Har{2}(\strip{\omega})$ and 
$f_{|\pm \omega} \in \We^2_v(\R)$. By definition of $\Har{2}$,
$f_{|0}\in \Ell{2}(\R)$.  Let $g := f^\vee = f_{|0}^{\, \vee}$. 
Then 
\[ v \ue^{\pm \omega \bfs} g = v (f_{|\pm \omega})^\vee \in
\Ell{2}(\R)
\]
and hence 
\[   v \ue^{\abs{\omega} \bfs}g 
= v\ue^{\omega \bfs} g \car_{\R_{\ge 0}} + 
v \ue^{-\omega \bfs} g \car_{\R_{< 0}}  \in \Ell{2}(\R)
\]
 as well. This proves (iii) and shows
 $\norm{f}_{\We^2_v(\strip{\omega})}
\lesssim \norm{f_{|\omega}}_{\We^2_v(\R)} + \norm{f_{|-\omega}}_{\We^2_v(\R)}$.

\prfnoi
(iii)$\dann$(i): Suppose that (iii) holds and let $g := f_{|0}^{\,
  \vee}$. Then $g\in \Ell{2}_{v, \omega}(\R)$ and hence, by definition,
$\fourier{g} \in \We^2_v(\strip{\omega})$. By Fourier inversion
$\fourier{g}_{|0} = f_{|0}$, and hence $f=\fourier{g}$ on
$\strip{\omega}$ by the identity theorem for holomorphic functions.
\end{proof}

Observe that $\Ho(\strip{\omega}) \subseteq \Har{2}(\strip{\omega})$
and hence, by the Paley--Wiener theorem, 
$\Ho(\strip{\omega}) \subseteq \We^2_\car(\strip{\omega})$. By
Proposition \ref{sob.p.We}.c we therefore obtain
\[   \Ho(\strip{\theta}) \subseteq \We^2_v(\strip{\omega}) \cap \uA_v(\cstrip{\omega})\qquad
(0\le \omega < \theta). 
\]
We shall need the following approximation results.

\begin{lem}\label{sob.l.den}
Let $v$ be an admissible function and $\omega \ge 0$. 
\begin{aufzi}
\item  The space $\Fourier(\Test(\R))$
is dense in $\We^2_v(\strip{\omega})$.
\item The space $\calS(\R)$ is contained and $\Test(\R)$ 
is dense in $\We^2_v(\R)$   and $\uA_v(\R)$. 

\item Let $\theta > \omega$. Then the space
\[   \{ f\in \Hol(\strip{\theta}) \suchthat \exists\, a > 0:   \abs{f}
\lesssim \ue^{-a \abs{\Re \bfz}^2}\}
\]
is dense in $\We_v^2(\strip{\omega})$.

\end{aufzi}
\end{lem}

\begin{proof}
a)\ 
By Proposition \ref{adm.p.Lpom}.d, $\Test(\R)$ is dense in $\Ell{2}_{v,
  \omega}(\R)$.  Hence, $\Fourier(\Test(\R))$ is dense in 
$\We^2_v(\strip{\omega})$. (This is analogous to Proposition \ref{sob.p.A}.d.)

\prfnoi
b)\ If $\eta \in \calS(\R)$ then $\eta^\vee \in \calS(\R)$ as well,
and since $v$ grows at most polynomially, $v \eta^\vee \in \Ell{1}(\R)
\cap \Ell{2}(\R)$. It follows that $\calS(\R) \subseteq \uA_v(\R) \cap 
\We^2_v(\R)$.

Next, we have to show that  the space
\[  E  := \{   v \eta^\vee \suchthat \eta \in \Test(\R)\} 
\]
is dense in $\Ell{2}(\R)$ and in $\Ell{1}(\R)$. We 
treat the $\Ell{2}$-case first,  the $\Ell{1}$-case being  similar.
Take  $g\in \Ell{2}(\R)$ such that 
\[   \int_\R  v \eta^\vee g = 0\quad \text{for all $\eta \in
  \Test(\R)$}.
\]
Replacing $\eta$ by $\eta \ast \vphi$ for $\eta, \vphi \in \Test(\R)$
we find
\[  \int_\R  v \eta^\vee \vphi^\vee g = 0\quad \text{for all $\vphi,
\, \eta \in
  \Test(\R)$}.
\]
Since $\{ \vphi^\vee \suchthat \vphi \in \Test(\R)\}$ is dense in
$\Ell{2}(\R)$, it follows that 
\[    v \eta^\vee g = 0 \quad \text{for all $\eta \in \Test(\R)$}.
\]
This implies $g = 0$ and hence b) is proved for $\We^2_v(\R)$.
The proof for $\uA_v(\R)$ is analogous.

\vanish{
\prfnoi
c)\ Write $\psi := \vphi^\vee$. Then $\psi_n(s)= n \psi(ns) = \vphi_n^\vee(s)$.
Moreover
\[ (\tau_t \vphi_n \cdot f)^\vee = (\ue^{\ui t \bf s} \psi_n) \ast
f^\vee 
\to f^\vee
\]
in $\Ell{2}_v(\R)$ by Proposition \ref{adm.p.Lpom}.h. 
}

\prfnoi
c)\  Let $\gauss_\veps = \ue^{-\veps \bfz^2}$, $\veps > 0$,  be a scaled Gaussian
function and let $f = \fourier{\vphi}$ for some 
$\vphi \in \Test(\R)$. Then $f\in \Ho(\strip{\theta}) \subseteq \Har{2}(\strip{\theta})$.
Hence,  $\Ho(\strip{\theta}) \ni \gauss_\veps f \to f$ 
in $\Har{2}(\strip{\theta})$ and, a fortiori, in
$\We^2_v(\strip{\omega})$
(Proposition \ref{sob.p.We}.c with $v=\car$). Now c) follows from a). 
\end{proof}

\vanish{
\begin{rem}
The proof shows more, namely:  {\em the space of
all entire functions satisfying
\[ \exists \, a,b>
0\,  \forall\, N\in \N\, \exists \, \gamma_N > 0 : 
\abs{f(\bfz)} \le   \gamma_n\,  \ue^{b \abs{\im \bfz}^2}\, \ue^{- a \abs{ \re
    \bfz}^2} 
\]
is dense in $\We^2_v(\strip{\omega})$.} 
\end{rem}
}

Also the following result shall be needed later.

\begin{lem}\label{sob.l.int-estimates}
Let $v: \R \to \R_{\ge 1}$ be admissible,  $0 \le \omega \le \theta$, 
and $\psi,\,  \varphi \in \Ho[\cstrip{\theta}]$. Then 
\[
\sup_{z\in \strip{\theta-\omega}} \left(\int_\R \| \tau_t \psi \cdot \tau_z \varphi\|_{\We_v^2(\strip{\omega})} \, \ud t\right)
+
\sup_{z\in \strip{\theta-\omega}} \left(\int_\R \| \tau_t \psi \cdot \tau_z \varphi\|_{\uA_v(\cstrip{\omega})} \, \ud t\right) <\, \infty.
\]
\end{lem}

\begin{proof}
Fix  $\theta' > \theta$ with $\psi, \, \vphi \in \Ho
(\strip{\theta'})$ and define $\omega' := \omega + (\theta' -  \theta) >
\omega$. Since $\Har{2}(\strip{\omega'})  =
\We^2_\car(\strip{\omega'}) \subseteq \We^2_v(\strip{\omega}) \cap
\uA_v(\cstrip{\omega})$
continuously, it suffices to prove
  \[
\sup_{z\in \strip{\theta-\omega}} \left(\int_\R \| \tau_t \psi \cdot \tau_z \varphi\|_{\Har{2}(\strip{\omega'})} \, \ud t\right) <\, \infty.
\] 
Choose $\alpha,\beta > 1$. Then $\abs{\psi} \lesssim (1 + \abs{\re
  \bfz})^{-\alpha}$ and $\abs{\vphi} \lesssim (1 + \abs{\re
  \bfz})^{-(\alpha + \beta)}$ on $\strip{\theta}$.
Hence, for $t\in \R$ and $z = x + \ui y\in
\strip{\theta-\omega}$,
\begin{align*}
 \norm{\tau_t \psi \cdot \tau_z \vphi}_{\Har{2}(\strip{\omega'})}
& = \norm{\tau_{t-x} \psi \cdot \tau_{\ui y}
  \vphi}_{\Har{2}(\strip{\omega'})}
\\ & = \sup_{\abs{y'} < \omega} \norm{\tau_{t-x} \psi_{|y'} \cdot
  \vphi_{|y'- y}}_{\Ell{2}(\R)}
\lesssim (1 + \abs{t-x})^{-\alpha}
\end{align*}
by Lemma \ref{aux.l.exp} (and its proof). The claim follows.
\end{proof}

\section{Generalized Hörmander Functions on the Real Line}\label{s.hrr}

In this section we generalize the classical Hörmander spaces on the
real line. Unless otherwise specified, $v$ denotes an arbitrary admissible
function and  $0 \neq \psi \in \Ce^\infty(\R)$.

A distribution $f$ is called a \emdf{generalized
  Hörmander
function} with respect to the \emdf{localizing function} $\psi$, if $\tau_t\psi\cdot f \in \We^2_v(\R)$ for
all $t \in \R$ and 
\[ \norm{f}_{\Hr^2_v(\R; \psi)} := \sup_{t\in \R} \norm{\tau_t\psi
  \cdot f}_{\We^2_v(\R)} < \infty.
\]
In this case, $f \in \Ell{2}_{\loc}(\R)$ and hence the use of the
word ``function'' is justified. We define \label{f.Hrr-psi}
\[  \Hr^2_v(\R; \psi) := \{ f\in \Ell{2}_{\loc}(\R) \suchthat
\text{$f$ is generalized Hörmander w.r.t. $\psi$}\}
\]
and endow it with the norm(!) from above. 
Let us collect some immediate properties.

\begin{prop}\label{hrr.p.Hr}
Let $v,\tilde{v}$ be admissible and    $0 \neq \psi \in
\Ce^\infty(\R)$.
Then the following assertions hold.
\begin{aufzi}
\item $\Hr^2_v(\R; \psi)$ is a Banach space continuously
included in $\Ell{2}_{\loc}(\R)$.

\item If $\tilde{v}\lesssim
  v$, then $\Hr^2_{v}(\R;\psi) \subseteq \Hr^2_{\tilde{v}}(\R;\psi)$
continuously. In particular, 
if $v, \tilde{v}$ are equivalent then 
$\Hr_v(\R;\psi) = \Hr_{\tilde{v}}(\R;\psi)$ with equivalent norms.

\item If $v$ is strongly admissible then $\Hr^2_v(\R;\psi) \subseteq \Cb(\R)$
  continuously.

\item Suppose $\theta > 0$ and  $\psi \in \Har{2}(\strip{\theta})$. 
Then  restriction to the real line yields a continuous inclusion
$\Ha(\strip{\theta})\subseteq \Hr^2_v(\R;\psi)$.
\end{aufzi}
\end{prop}

\begin{proof}
a) and b) are  straightforward.

\prfnoi
c)\ By Proposition \ref{sob.p.We}.d,
$\tau_t \psi \cdot f\in \We^2_v(\R) \subseteq \uA(\R) \subseteq
\Co(\R)$ for each $t\in \R$. As $\psi \neq 0$, it follows that $f\in
\Ce(\R)$. 
Now pick $x_0\in \R$ with  $\abs{\psi(x_0)} >0$;  
then for each $x\in \R$ 
\[ f(x) = \frac{1}{\psi(x_0)} (\tau_{x-x_0}\psi \cdot f)(x)
\] 
and hence
$\abs{f(x)} \lesssim \norm{\tau_{x{-}x_0} \psi \cdot f}_{\infty}
\lesssim  \norm{\tau_{x{-}x_0} \psi \cdot f}_{\We^2_v(\R)}
\le \norm{f}_{\Hr^2_v(\R;\psi)}$.

\prfnoi
d) Let  $f\in
  \Ha(\strip{\theta})$.  Then, by Proposition \ref{sob.p.We}.c, 
\[   \norm{\tau_t \psi\cdot f}_{\We^2_v(\R)} \lesssim
 \norm{\tau_t \psi\cdot f}_{\Har{2}(\strip{\theta})} \le
 \norm{\psi}_{\Har{2}(\strip{\theta})}  \norm{f}_{\Ha(\strip{\theta})}.\qedhere
\]
\end{proof}

Our aim is to show that 
the space $\Hr^2_v(\R;\psi)$ and its topology does not depend on
$\psi$ when $\psi$ is a Schwartz function. We shall need the
following lemma.

\begin{lem}\label{hrr.l.mea}
Let $0 \neq \psi \in \Ce^\infty(\R)$ and $f$ a function on $\R$. 
\begin{aufzi}
\item Suppose that  $f$ is \emdf{locally in $\We^2_v(\R)$}, i.e. 
$\eta f\in \We^2_v(\R)$ for each $\eta \in \Test(\R)$ 
Then for each $\eta \in
  \Test(\R)$ the mapping  
\[ \R \to \We^2_v(\R),\qquad t\mapsto \tau_t\eta \cdot f
\]
is continuous. 

\item If $f$ is such that $\tau_t \psi \cdot f \in \We^2_v(\R)$ for
  each $t\in \R$  then  $f$ is locally in $\We^2_v(\R)$ and 
the mapping 
\[    \R \to \We^2_v(\R),\qquad t \mapsto \tau_t \psi \cdot f
\]
is strongly measurable.
\end{aufzi}
\end{lem}

\begin{proof}
a) Fix $a > 0$ and 
choose a test function $\vphi$ such that $\vphi \equiv 1$ on $[-a,a]+
\supp(\eta)$. Then 
$\tau_t \eta \cdot f = \tau_t \eta \cdot (\vphi f)$ for $\abs{t} \le
a$.  Since $\vphi f\in \We_v^2(\R)$, it suffices to show that 
$\tau_{\bft}\eta \in \Ce(\R; \uA_v(\R))$. But this follows from 
Lemma \ref{sob.l.den}.b and  Proposition \ref{sob.p.A}.g.

\prfnoi
b)\ For the first assertion, use a smooth
partition of unity and the fact that $\Test(\R) \subseteq \uA_v(\R)$. 

Now fix $\vphi \in \Test(\R)$ with $\vphi(0) = 1$ and let $\vphi_n :=
\vphi(\bfs/n)$ for $n \in \N$. 
As $f$ is locally  in $\We^2_v(\R)$, by a) 
the function  $\tau_\bft(\psi \vphi_n) \cdot f$
 is continuous for each $n\in \N$. Therefore,
it suffices to show that $\tau_t(\vphi_n \psi)\cdot f \to \tau_t \psi
\cdot f$ in $\We^2_v(\R)$ for each $t\in \R$.

Let $\eta := \vphi^\vee$ and $\eta_n := n \eta(n \bfs) = \vphi_n^\vee$ for $n \in
\N$.  Then $\eta$ is a Schwartz function and hence satisfies the 
assumptions made in Proposition \ref{adm.p.Lpom}, part h) with $\omega = 0$
and $\tilde{v}= v_\alpha = (1+ \abs{\bfs})^\alpha$ for $\alpha > 0$
large enough.  Hence, with  $t\in \R$ fixed, 
\[  (\tau_t(\vphi_n \psi)\cdot f)^\vee = (\ue^{\ui \bfs t} \eta_n)
\ast (\tau_t\psi \cdot f)^\vee \to \vphi(0) \cdot (\tau_t\psi \cdot f)^\vee 
= (\tau_t\psi \cdot f)^\vee 
\]
in $\Ell{2}_{v,0}(\R)$ as $n \to \infty$.
Taking Fourier
transforms yields $\tau_t(\vphi_n \psi)\cdot f \to \tau_t\psi \cdot f$
in $\We^2_v(\R)$, and the proof is complete. 
\end{proof}

\vanish{
\begin{lem}
Let $f: \R \to \C$. 
\begin{aufzi}

\item For a function $f$ on $\strip{\omega}$ the following assertions
are equivalent:
\begin{aufzii}
\item $f$ is quasi-locally in $\We^2_v(\strip{\omega})$;
\item For some $0 \neq \psi \in \Ho[\cstrip{\omega}]$ and each $t\in \R$:
  $\tau_t\psi \cdot f \in \We^2_v(\R)$.
\end{aufzii}

\item If there is $\Psi \subseteq \Ce^\infty(\R)$ such that 
$\psi f\in \We^2_v(\R)$ for all $\psi \in \Psi$ and $\bigcap_{\psi\in
  \Psi}
\set{\psi = 0} = \leer$, then $f$ is locally in $\We^2_v(\R)$.

\end{aufzi}
\end{lem}
}

The  classical Hörmander spaces are defined as above (with $v=
v_\alpha = (1 + \abs{\bfs})^\alpha$) with  test functions
$\psi \in \Test(\R)$ as localizers.
We shall show  that one can, more generally,
use  Schwartz functions.


\begin{thm}\label{hrr.t.indep}
Let $\eta, \psi \in \calS(\R)\ohne \{0\}$. If
$f\in \Hr^2_v(\R; \eta)$,  then $f\in \Hr^2_v(\R;\psi)$ with
\[   \norm{f}_{\Hr^2_v(\R;\psi)} \lesssim \norm{f}_{\Hr^2_v(\R; \eta)}.
\]
\end{thm}

\begin{proof}
As $\eta \neq 0$ there is no loss of generality to suppose that 
$\int_\R \abs{\eta}^2  = 1$. 
We claim that 
\beq\label{hrr.eq.intrep}
 \psi f = \int_\R  (\tau_t\konj{\eta}\cdot \psi) \, 
(\tau_t \eta \cdot f) \, \ud{t}
\eeq
as a Bochner integral in $\We^2_v(\R)$. Note that the function
$\tau_\bft \konj{\eta} \cdot \psi$ is continuous $\R \to \uA_v(\R)$
and, by Lemma \ref{hrr.l.mea}, the function $\tau_\bft \eta \cdot f$ is
strongly measurable $\R \to \We^2_v(\R)$. Hence, the function
$(\tau_\bft\konj{\eta}\cdot \psi) \, 
(\tau_\bft \eta \cdot f)$ is strongly measurable $\R \to \We^2_v(\R)$.
Moreover, 
\[ \int_{\R} \norm{ (\tau_t\konj{\eta}\cdot \psi) \, 
(\tau_t \eta \cdot f)}_{\We^2_v(\R)} \, \ud{t}
\lesssim \int_\R \norm{\tau_t\konj{\eta}\cdot \psi}_{\uA_v(\R)} \,
\ud{t}\, 
\norm{f}_{\Hr^2_v(\R;\eta)}.
\]
In order to see that 
\[
 \int_\R \norm{\tau_t\konj{\eta}\cdot \psi}_{\uA_v(\R)} \,
\ud{t} < \infty,
\]
we choose $N\in \N$ so
large
that $v/ v_N \in \Ell{2}(\R)$. Then $\We^{N,2}(\R) = \We^2_{v_N}(\R)
\subseteq \uA_v(\R)$ continuously. Since the norm in $\We^{N,2}(\R)$
effectively looks at the $\Ell{2}$-norms of all derivatives up to order
$N$, it suffices to have
\beq\label{hrr.eq.Bochnerint}
 \int_\R \norm{ \tau_t \konj{\eta}^{(j)} \cdot
  \psi^{(N{-}j)}}_{\Ell{2}(\R)} \, \ud{t} < \infty
\eeq
for all $0 \le j \le N$. As $\konj{\eta}^{(j)}$ and
$\psi^{(N{-}j)}$
are both Schwartz functions, \eqref{hrr.eq.Bochnerint} 
follows from  Lemma \ref{aux.l.exp}.

Now that we know that the integral \eqref{hrr.eq.intrep} 
exists in $\We^2_v(\R)$ in the Bochner sense, 
it remains to determine its value. 
To this aim, let $\vphi \in \Test(\R)$ be arbitrary. Then 
\[  \vphi \int_\R (\tau_t\konj{\eta}\cdot \psi) \, 
(\tau_t \eta \cdot f) \, \ud{t} =  \int_\R \tau_t \abs{\eta}^2 \cdot 
(\vphi \psi f) =  \Bigl(\int_\R \abs{\eta}^2 \Bigr) \, \vphi \psi f
= \vphi \psi f
\]
in $\Ell{2}(\R)$ by Lemma \ref{aux.l.Lp}, 
since $\vphi \psi f \in \Ell{2}(\R)$.
This implies \eqref{hrr.eq.Bochnerint}.

Finally, take $s\in \R$ and replace $\psi$ above by $\tau_s \psi$. 
Then we find
\[   \norm{\tau_s \psi \cdot f}_{\We^2_v(\R)} \le 
\int_\R \norm{\tau_t \konj{\eta} \cdot \tau_s\psi}_{\uA_v(\R)} \,
\ud{t}\,\, \norm{f}_{\Hr^2_v(\R)}
=  \int_\R \norm{\tau_t \konj{\eta} \cdot \psi}_{\uA_v(\R)} \,
\ud{t}\,\, \norm{f}_{\Hr^2_v(\R;\eta)}
\]
since $\norm{\tau_t \konj{\eta} \cdot \tau_s\psi}_{\uA_v(\R)}
= \norm{\tau_{t-s} \konj{\eta} \cdot \psi}_{\uA_v(\R)}$ for all 
$t, s\in \R$. This concludes the proof. 
\end{proof}

With the last result at hand we are in the position to define
the \emdf{generalized Hörmander space} on $\R$ associated with the
admissible function $v$ as \label{f.Hrr}
\[ \Hr^2_v(\R) := \Hr^2_v(\R; \psi),
\]
with $\psi$ being {\em any}  non-zero Schwartz function.

\vanish{
\begin{exa}\label{hrr.exa.exp} 
Let $v: \R \to [1, \infty )$ be admissible. 
Then $\ue^{-\ui s \bfz} \in \Hr^2_v(\R)$ for any 
$s\in \R$, and
\[  \norm{ \ue^{\ui s \bfz}}_{\Hr^2_v(\R)}\lesssim v(s)  \qquad (s\in \R).
\]
Indeed: fix $0 \neq \psi \in \calS(\R)$ and $s\in \R$. Then 
\[
\| \psi \ue^{-\ui s \bfz}\|_{\We_v^2(\R)}\, =\, 
\big\| v \cdot \tau_{s}\psi^\vee \big\|_2\, =\, \big\| \tau_{-s} v \cdot \psi^\vee \big\|_2\, \lesssim_v \, v(s)\|\psi\|_{\We_v^2(\R)}.
\] 
Hence,
\[
\sup_{t\in \R} \| \tau_t \psi \cdot \ue^{-\ui s\bfz} \|_{\We_v^2
  \R)}\, \lesssim_{v, \psi}\, v(s) \, <\, \infty
\] for each $s\in \R$.
\end{exa}
}

\section{Generalized Hörmander Functions on Strips}
\label{s.hrs}

We now want to pass to holomorphic functions on strips. Again,
$v$ denotes an arbitrary admissible function and  
$\omega \ge 0$ is an arbitrary  non-negative real number.

\medskip

Let $0 \neq \psi \in \Ho[\cstrip{\omega}]$. 
A function $f$ on $\strip{\omega}$ 
is called \emdf{generalized
  Hörmander} with respect to $\psi$, 
if $\tau_t\psi\cdot f \in \We^2_v(\strip{\omega})$ for
all $t \in \R$ and 
\[ \norm{f}_{\Hr^2_v(\strip{\omega}; \psi)} := \sup_{t\in \R} \norm{\tau_t\psi
  \cdot f}_{\We^2_v(\strip{\omega})} < \infty.
\]
We let \label{f.Hrst-psi}
 \[  \Hr^2_v(\strip{\omega}; \psi) := \{ f\suchthat
\text{$f$ is generalized Hörmander w.r.t. $\psi$}\}.
\]
Note that if $\omega = 0$ then this definition is the same
as the one in the previous section. More generally, 
 restriction to $\R$ induces an 
embedding $\Hr^2_v(\strip{\omega};\psi) \subseteq \Hr^2_v(\R)$ 
(since it also induces an embedding 
$\We^2_v(\strip{\omega}) \subseteq \We^2_v(\R)$ by 
Proposition \ref{sob.p.We}.c).

Note that each non-zero
 $\psi \in \Ho[\R]$ restricts to a non-zero 
Schwartz function on $\R$ (Lemma
\ref{hrr.l.Ho}), hence  $\Hr^2_v(\strip{0}; \psi) = 
\Hr^2 _v(\R)$, as defined at the end of Section \ref{s.hrr}.


\begin{prop}\label{hrs.p.Hr}
Let $v, \tilde{v}$ be admissible, $\omega>0$ and $0 \neq \psi \in
\Ho[\cstrip{\omega}]$. 
Then the following assertions hold:
\begin{aufzi}

\item $\Hr^2_v(\strip{\omega};\psi)$ is a Banach space with respect 
to the norm $\norm{\cdot}_{\Hr^2_v(\strip{\omega});\psi)}$ with 
$\We^2_v(\strip{\omega}) \subseteq \Hr^2_v(\strip{\omega};\psi)$,
continuously.

\item Restriction yields continuous embeddings 
\begin{align*} \Hr^2_v(\strip{\omega}; \psi) \subseteq  
\Ha(\strip{\theta}) \quad & \text{if $0 < \theta < \omega$\quad and}\\
\Ha(\strip{\theta}) \subseteq \Hr^2_v(\strip{\omega};\psi) \quad & \text{if
                                          $\theta > \omega$}. 
\end{align*}
In particular, $\Hr^2_v(\strip{\omega};\psi)$ embeds continuously into
$\Hol(\strip{\omega})$.

\item If $\tilde{v}\lesssim v$ then $\Hr^2_v(\strip{\omega};\psi) \subseteq
 \Hr^2_{\tilde{v}}(\strip{\omega};\psi)$ continuously. In
 particular,
if $v$ and $\tilde{v}$ are equivalent, then $\Hr^2_v(\strip{\omega};\psi) =
 \Hr^2_{\tilde{v}}(\strip{\omega};\psi)$ with equivalent
 norms.

\item Pointwise multiplication is a bounded bilinear mapping
\[   \uA_v(\cstrip{\omega}) \times \Hr^2_v(\strip{\omega};\psi) 
\to \Hr^2_v(\strip{\omega};\psi). 
\]
\end{aufzi}
\end{prop}

\begin{proof}
a) and b)\ 
It is clear that $\Hr^2_v(\strip{\omega};\psi)$ is a linear space
and $\norm{\cdot}_{\Hr^2_v(\strip{\omega};\psi)}$ is a seminorm on
it. The continuous inclusion 
$\We^2_v(\strip{\omega}) \subseteq \Hr^2_v(\strip{\omega};\psi)$
follows
since $\psi \in \uA_v(\cstrip{\omega})$ whereon  translation is a
strongly continuous and isometric group.

Now let $f\in \Hr^2_v(\strip{\omega};\psi)$. Then  for each $t \in \R$ 
\[ f = \frac{\tau_t \psi \cdot f}{\tau_t \psi} \quad
\text{on}\quad \set{\tau_t \psi \neq 0}. 
\]
As the zeroes of $\psi$ form a discrete set, we can vary $t$ and find  $f\in
  \Hol(\strip{\omega})$. Moreover, it follows that
  $\norm{\cdot}_{\Hr^2_v(\strip{\omega};\psi)}$ is a norm.

Next,  fix $0 < \theta < \omega$. Then $\We^2_v(\strip{\omega}) \subseteq
\Ha(\strip{\theta})$ continuously by Proposition \ref{sob.p.We}.c.  It follows that for each $t\in \R$
and $\delta > 0$ 
\[  \abs{f}  = \frac{ \abs{\tau_t \psi \cdot f}}{\abs{\tau_t \psi}}
\le \frac{1}{\delta} \norm{\tau_t \psi \cdot f}_{\Ha(\strip{\theta})}
\lesssim \frac{1}{\delta} \norm{f}_{\Hr^2_v(\strip{\omega};\psi)}
\]
on the set $\set{\abs{\tau_t \psi} > \delta} \cap \strip{\theta}$.
Hence, by choosing $\delta > 0$ sufficiently small and varying $t\in
\R$ we find
\[  \norm{f}_{\Ha(\strip{\theta})} \lesssim
\norm{f}_{\Hr^2_v(\strip{\omega};\psi)}.
\]
This proves the first assertion of b). 
  
To show completeness, note that the mapping
\[  \Phi: \Hr^2_v(\strip{\omega};\psi) \to \ell^\infty(\R;
\We^2_v(\strip{\omega})), \qquad
\Phi f := (\tau_t \psi \cdot f)_{t\in \R}
\]
is an isometric embedding. Hence, it suffices to show that $\ran\Phi$
is closed. To this end, let  $(f_n)_n$ be a sequence in
$\Hr^2_v(\strip{\omega};\psi)$ and $g := (g_t)_t \in \ell^\infty(\R;
\We^2_v(\strip{\omega}))$ with $\Phi f_n \to g$. Then $(f_n)_n$ is
Cauchy in  $\Hr^2_v(\strip{\omega};\psi)$ and hence there is $f\in
\Hol(\strip{\omega})$ such that $f_n \to f$ uniformly on each 
$\strip{\theta}$, $0< \theta < \omega$. It follows that 
$g_t = \tau_t\psi f$ for each $t\in \R$ and hence $f\in
\Hr^2_v(\strip{\omega};\psi)$ with $\Phi f =g$. This concludes the
proof of a).  

It remains to show the second assertion in b). Fix $\theta > \omega$
and $f\in \Ha(\strip{\theta})$. Without loss of generality we may
suppose that $\psi \in \Ho(\strip{\theta})$. Then 
$\tau_t\psi \cdot f \in \Har{2}(\strip{\theta})$ with 
\[  \norm{\tau_t \psi \cdot f}_{\Har{2}(\strip{\theta})} 
\le \norm{f}_{\Ha(\strip{\theta}} \norm{\tau_t\psi}_{\Har{2}(\strip{\theta})}
= \norm{f}_{\Ha(\strip{\theta}} \norm{\psi}_{\Har{2}(\strip{\theta})}.
\]
As $\Har{2}(\strip{\theta}) = \We^1_\car(\strip{\theta}) \subseteq
\We^2_v(\strip{\omega})$ continuously, the assertion follows.

\prfnoi
c) is obvious and e) follows from Proposition \ref{sob.p.We}.f. 
\end{proof}

Next, we aim at showing that the 
generalized Hörmander space $\Hr^2_v(\strip{\omega};\psi)$ 
does not depend on the particular choice of $0\neq \psi \in
\Ho[\cstrip{\omega}]$. For $\omega = 0$ this has been done in the
previous  section (Theorem \ref{hrr.t.indep}), 
as $\Ho[\R]$-functions restrict to Schwartz functions on $\R$
(Lemma \ref{hrr.l.Ho}). 
 
For $\omega > 0$ the matter is
a little more delicate since we do not have the means of a (finite) 
partition of unity at hand, which played an important role in the
proof of Lemma \ref{hrr.l.mea}.

Recall that for $\omega\ge 0$, $\theta>0$, 
and $\psi: \strip{\theta{+}\omega}\to \C$ we write
\[
(\tau_z \psi)(w)\, =\, \psi(w-z)\quad (z\in \strip{\theta},\, w\in \strip{\omega}).
\]
Still, unless explicitly noted otherwise, $v$ is an arbitrary
admissible function.

\begin{lem}\label{hrs.l.mea}
Let $\theta >0$, $\omega\geq 0$, $0\neq \psi \in \Ho (\strip{\theta{+}\omega})$, and
$f\in \Hr_v^2(\strip{\omega};\psi)$. Then the following statements hold:
\begin{aufzi}
\item The mapping
\[
\R \to \We_v^2 (\strip{\omega}), \quad t\mapsto \tau_t\psi \cdot f,
\] 
is weakly continuous.
\item If $\tau_z \psi\cdot f\in \We_v^2 (\strip{\omega})$ for all
  $z\in \strip{\theta - \omega}$, and  
\[
\sup_{z\in \strip{\theta}} \| \tau_z \psi \cdot f\|_{\We_v^2 (\strip{\omega})} \, <\, \infty,
\]
then the mapping
\[
\strip{\theta}\to \We_v^2 (\strip{\omega}), \quad z\mapsto \tau_z \psi \cdot f,
\] is holomorphic.
\end{aufzi}
\end{lem}

\begin{proof}
For {\rm a)} let $\rho\coloneqq 0$ and for {\rm b)}, let
$\rho\coloneqq \theta$. Define
\[ H: \strip{\rho} \to \Ell{2}(\R),\qquad H(z) := \tilde{v}
\cosh(\omega \bfs) (\tau_z \psi \cdot f)^\vee,
\]
where $\tilde{v}$ is equivalent to $v$ and smooth (Lemma \ref{adm.l.AWalwaysSmooth}). 
By the assumptions made in a) and b), $H$ is bounded. 

Next, fix $\eta\in \Test(\R)$ and define $\vphi := \tilde{v}
\cosh(\omega\bfs) \eta \in \Test(\R)$. Then 
\[  I(z) :=  \int_\R H(z) \eta = \int_\R \tilde{v} \cosh(\omega \bfs) (\tau_z \psi \cdot f)^\vee \, \eta
= \int_\R (\tau_z \psi \cdot f)^\vee\, \vphi 
= \int_\R \tau_z \psi \cdot f \cdot \vphi^\vee 
\]
for each $z\in \strip{\rho}$. As $f_{|0} \in \Hr^2_v(\R)$ and
$\vphi^\vee \in \calS(\R)$, $f\vphi^\vee \in \We^2_v(\R)\subseteq
\Ell{2}(\R)$. Hence,   $I: \strip{\rho} \to \C$ is continuous.  

Since $\Test(\R)$ is dense in
$\Ell{2}(\R)$, this already completes the proof of a). In case b)
one employs Morera's theorem in combination with a Fubini argument to
see that $I$ is holomorphic. 
\end{proof}

Lemma \ref{hrs.l.mea} is preliminary, as the next
result entails  much stronger statements.

\begin{thm}\label{hrs.t.indep}
Let $\omega \geq 0$ and  $0\neq \psi \in \Ho [\cstrip{\omega}]$.
Furthermore, let $\theta > 0$ and $\vphi \in \Ho(\strip{\theta{+}\omega})$.
Then the following assertions hold:
\begin{aufzi}
\item  For each $f\in \Hr_v^2 (\strip{\omega};\psi)$ 
\[   \tau_\bfz\vphi \cdot f\in \Ha( \strip{\theta};
\We_v^2(\strip{\omega})) 
\]
and 
\[  \sup_{z\in \strip{\theta}} \norm{\tau_z \vphi \cdot
    f}_{\We^2_v(\strip{\omega})} \lesssim 
\norm{f}_{\Hr^2_v(\strip{\omega};\psi)}
\]
(with the constant not depending on $f$). 

\item For each $f\in \Hr_v^2(\strip{\omega};\psi)$  and $z\in
  \strip{\theta}$ one has $(\tau_z\vphi \cdot f)^\vee \in
  \Ce(\R)$ with 
\[  \sup_{z\in \strip{\theta}, s\in \R} v(s)  \ue^{\omega \abs{s}}
  \abs{(\tau_z \vphi \cdot f)^\vee(s)} \lesssim
\norm{f}_{\Hr^2_v(\strip{\omega}; \psi)}
\]
(with the constant not depending on $f$).

\item  If $\omega > 0$ one has the representation formula
\beq\label{hrs.eq.rep}
 \vphi(0)\cdot  f(z) = \int_\R (\tau_z\vphi \cdot f)^\vee(s) \ue^{-\ui s z} \,
\ud{s} \qquad (z\in \strip{\theta}\cap \strip{\omega})
\eeq
whenever $f\in \Hr^2_v(\strip{\omega};\psi)$.
\end{aufzi}
\end{thm}

\begin{proof}
Define $\psi^* := \konj{\psi(\konj{\bfz})}
\in  \Ho[\cstrip{\omega}]$ and  fix $f\in \Hr^2_v(\strip{\omega};\psi)$ 
and $z\in \strip{\theta}$. Then $\tau_z \vphi \in \Ho[\strip{\omega}]$.

\prfnoi
a)\ By   Lemma \ref{hrs.l.mea}, the
function $F_z :=  \tau_\bft(\psi \psi^*) \cdot (\tau_z\vphi  \cdot f)
= (\tau_\bft \psi^*\cdot \tau_z\vphi) \cdot (\tau_\bft\psi \cdot f)$ 
is bounded and weakly continuous $\R \to \We_v^2(\strip{\omega})$. 
Moreover, 
\begin{align*}
 \sup_{z\in \strip{\theta}} \int_\R & \norm{\tau_t(\psi \psi^*) \cdot (\tau_z\vphi\cdot
  f)}_{\We_v^2(\strip{\omega})} \, \ud{t} 
\\ & \lesssim \sup_{z\in \strip{\theta}} \int_\R \norm{\tau_t \psi^* \cdot \tau_z\vphi}_{\uA_v(\cstrip{\omega})}
\, \ud{t} \,\,\norm{f}_{\Hr^2_v(\strip{\omega};\psi)} < \infty
\end{align*}
by Lemma \ref{sob.l.int-estimates}. 
Since $\We^2_v(\strip{\omega})$ is reflexive, the 
integral
\[ g_z := \int_\R \tau_t(\psi^* \psi)\cdot  \tau_z\vphi \cdot f\, \ud{t}
\]
exists in $\We^2_v(\strip{\omega})$ in the weak sense
(as defined, e.g., in \cite[Them 3.26]{RudinFA}).  Actually, one can
say more here: since $\We^2_v(\strip{\omega})$ is separable, Pettis'
measurability theorem yields that $F_z$ is strongly measurable 
$\R \to \We^2_v(\strip{\omega})$ and hence, by the estimate above,
$F_z \in \Ell{1}(\R; \We_v^2(\strip{\omega}))$. 

We claim that $g_z = \norm{\psi}_{\Ell{2}(\R)}^2 \cdot \tau_z\vphi
\cdot f$,  i.e.,
\[
  \norm{\psi}_{\Ell{2}(\R)}^2 \cdot \tau_z\vphi
\cdot f  =    \int_\R \tau_t(\psi^* \psi)\cdot  \tau_z\vphi \cdot f\, \ud{t}.
\]
For $\omega= 0$, this has already been established in the proof of Theorem
\ref{hrr.t.indep}. In the case $\omega > 0$ evaluations at points $x\in \R$
are bounded linear functionals on $\We^2_v(\strip{\omega})$ (by Proposition
\ref{sob.p.We}.b). Inserting $x\in \R$ yields
\[  g_z(x) = \int_\R \abs{\psi(x-t)}^2\, \ud{t}\,\,  \vphi(x-z) f(x)
=   \norm{\psi}_{\Ell{2}(\R)}^2\cdot (\tau_z\vphi \cdot f)(x).
\]
Since both $g_z$ and $\tau_z\vphi \cdot f$ are holomorphic functions on
$\strip{\omega}$, our claim is proved. 

Now, since $\psi \neq 0$ also $\norm{\psi}_{\Ell{2}(\R)}^2 \neq
0$. Therefore,  $\tau_z\vphi \cdot f \in \We^2_v(\strip{\omega})$ with 
\[  \norm{\tau_z\vphi \cdot f}_{\We^2(\strip{\omega})} \lesssim 
\norm{f}_{\Hr^2(\strip{\omega};\psi)},
\]
where the hidden constant depends on $\vphi, \psi, v, \omega$, but neither
on $f$ nor on $z\in \strip{\theta}$.  

Finally, apply b) of Lemma \ref{hrs.l.mea} to find that $\tau_\bfz \vphi
\cdot f$ is holomorphic $\strip{\theta} \to
\We^2_v(\strip{\omega})$.  This concludes the proof of a). 

\prfnoi
b)\  We again employ the representation formula
\beq\label{hrs.eq.CRF}    \vphi f = \int_\R (\tau_t\psi^* \cdot \vphi) \cdot
(\tau_t \psi \cdot f)\, \ud{t}.
\eeq
By a) applied with $\vphi = \psi$ we know
that  $\tau_\bft \psi \cdot f \in \Cb(\R; \We^2_v(\strip{\omega}))$,
and hence \[ (\tau_\bft \psi \cdot f)^\vee \in 
\Cb(\R; \Ell{2}_{v, \omega}(\R)).
\]
On the other hand, by Lemma \ref{sob.l.int-estimates},
\[   \tau_\bft \psi^* \cdot \vphi \in 
\Ell{1}(\R; \uA_v(\cstrip{\omega}) \cap \We^2_v(\strip{\omega}))
\]
and hence
\[ (\tau_\bft \psi^* \cdot \vphi)^\vee \in 
\Ell{1}(\R; \Ell{1}_{v, \omega}(\R) \cap \Ell{2}_{v,\omega}(\R)).
\]   
Taking the convolution, we obtain 
(by Proposition \ref{adm.p.Lpom}.f)
\[ (\tau_\bft \psi^* \cdot \vphi)^\vee  \ast (\tau_\bft \psi \cdot
f)^\vee \in \Ell{1}(\R; \Ell{2}_{v, \omega}(\R) \cap \Ell{\infty}_{v,
  \omega}(\R)\cap \Ce(\R) )
\]
and hence
\[  h :=  \int_\R (\tau_t \psi^* \cdot \vphi)^\vee  \ast (\tau_t \psi \cdot
f)^\vee \, \ud{t} \in \Ell{2}_{v, \omega}(\R) \cap \Ell{\infty}_{v,
  \omega}(\R) \cap \Ce(\R). 
\]
Applying the Fourier transform (which is bounded $\Ell{2}_{v,\omega}(\R)
\to \We^2_v(\strip{\omega})$ and transforms  the convolution into the 
pointwise product) yields
\[  \fourier{h} = \int_\R  
(\tau_t \psi^* \cdot \vphi) \cdot (\tau_t \psi \cdot
f) \, \ud{t} = \vphi\cdot f
\]
and hence $h = (\vphi \cdot f)^\vee$.  Replacing $\vphi$ by 
$\tau_z \vphi$ in this argument proves the first assertion of b). The
second follows by tracing through the estimates hidden in 
this qualitative argument.

\prfnoi
c)\  For $z\in \strip{\omega}$ and
$w\in \strip{\theta}$ one has
\[ \vphi(z-w)f(z) = (\tau_w \vphi\cdot f)(z)
= \int_\R (\tau_w\vphi \cdot f)^\vee(s) \ue^{-\ui s z}.
\]
Putting $w=z$ yields \eqref{hrs.eq.rep}.
\end{proof}

\begin{rems}
1.)\ Taking into account Remark \ref{adm.r.Lpom}.1, we see that in part b) of
Theorem \ref{hrs.t.indep} one even has
$\tilde{v} \ue^{\omega\abs{\bfs}} (\tau_z \vphi \cdot f)^\vee\in \Co(\R)$
whenever $\tilde{v}$ is continuous and equivalent to $v$.

\prfnoi
2.)\ In the case $\omega = 0$, assertion b) from above holds
for $\theta = 0$ and an arbitrary  Schwartz function $\vphi$ on $\R$. This follows by making appropriate
changes in the given proof.  

\prfnoi
3.) In general, the representation formula
\eqref{hrs.eq.rep} need not hold if $\omega = 0$, as
in this case $z\in \R$ and $(\tau_z\vphi \cdot f)^\vee$ may only be
in $\Ell{2}(\R)$. However, if $v$ is strongly admissible then
$(\tau_z\vphi f)^\vee \in \Ell{1}(\R)$ so that the
\eqref{hrs.eq.rep} still holds. More generally, if $v$ is strongly
admissible then \eqref{hrs.eq.rep} holds for $z\in
\strip{\theta} 
\cap \cstrip{\omega}$. 
\end{rems}

As a consequence of Theorem  \ref{hrs.t.indep} 
we may now define \label{f.Hrst}
\[ \Hr^2_v(\strip{\omega}) := \Hr^2_v(\strip{\omega};\psi)
\]
where $0 \neq \psi \in \Ho[\strip{\omega}]$ is arbitrary. 
In order to have a canonical norm, we let
\[  \norm{f}_{\Hr^2_v(\strip{\omega})} := 
\norm{f}_{\Hr^2_v(\strip{\omega}; G)} = 
\sup_{t \in \R} \norm{\tau_t\gauss\cdot f}_{\We^2_v(\strip{\omega})}
\]
where $\gauss  = \ue^{-\bfz^2}$ is the Gaussian function.

\begin{exa}\label{hrs.exa.exp} 
Let $v: \R \to [1, \infty )$ be admissible and $\omega \ge 0$. 
Then 
$\ue^{-\ui s \bfz} \in \Hr^2_v(\strip{\omega})$ for each
$s\in \R$, and
\[  \norm{ \ue^{\ui s \bfz}}_{\Hr^2_v(\R)}\lesssim v(s) \ue^{\omega \abs{s}}  \qquad (s\in \R).
\]
Indeed: fix $0 \neq \psi \in \Ho[\strip{\omega}]$ and $s\in \R$. Then 
\[
\| \psi \ue^{-\ui s \bfz}\|_{\We_v^2(\strip{\omega})}\, =\, 
\big\| v \ue^{\omega \abs{\bfr}} \tau_{s}\psi^\vee \big\|_2\, =\,
  \big\| \tau_{-s} v \ue^{\omega \abs{ \bfr + s}} \psi^\vee
  \big\|_2\, \lesssim_v \, v(s) \ue^{\omega \abs{s}} \|\psi\|_{\We_v^2(\strip{\omega})}.
\] 
Hence,
\[
\sup_{t\in \R} \| \tau_t \psi \cdot \ue^{-\ui s\bfz} \|_{\We_v^2(\strip{\omega})}\, \lesssim_{v, \psi}\, v(s) \ue^{\omega \abs{s}}\, <\, \infty.
\]
\end{exa}

\bigskip

In the next theorem we collect the properties of Hörmander
functions in the strongly admissible case.

\begin{thm}\label{hrs.t.strongadm}
Let $\omega \ge 0$ and $v: \R \to \R_{\ge 1}$ be strongly admissible. 
Then the following assertions hold:
\begin{aufzi}
\item The space $\Hr^2_v(\strip{\omega})$ embeds continuously into
  $\UCb(\cstrip{\omega})$.

\item For $\omega > 0$, each 
 $f\in
  \Hr^2_v(\strip{\omega})$ extends continuously to $\cstrip{\omega}$,
  and its boundary functions $f_{|\pm \omega}$ are in $\Hr^2_v(\R)$.
  Conversely, if $f\in \Ha(\strip{\omega})$ is such that its 
  $\Ell{\infty}$-boundary values $f_{|\pm\omega}$ are contained in 
$\Hr^2_v(\R)$, then $f\in \Hr^2_v(\strip{\omega})$. Moreover
one has equivalence of norms
\[ \norm{f}_{\Hr^2_v(\strip{\omega})} \simeq 
\norm{f_{|\omega}}_{\Hr^2_v(\R)} + \norm{f_{|-\omega}}_{\Hr^2_v(\R)}.
\]

\item The space $\Hr^2_v(\strip{\omega})$ is (after a scaling of the
  norm)
a Banach algebra with respect to pointwise multiplication.

\item Let $\theta > 0$ and $\vphi\in \Ho(\strip{\theta{+}\omega})$. 
Then for each $f\in \Hr^2_v(\strip{\omega})$ one has
\[  \vphi(0)\cdot f(z) = \int_\R (\tau_z\vphi\cdot f)^\vee(s) \,
\ue^{-\ui s z}\, \ud{s} \qquad (z\in \strip{\theta} \cap
\cstrip{\omega}).
\]
\end{aufzi}
\end{thm}

\begin{proof}
We fix, once and for all, a function $0\neq \psi \in \Ho[\strip{2\omega}]$.

\prfnoi
a) Let $\theta >0$ and $\vphi\in \Ho(\strip{\theta{+}\omega})$ and
$f\in \Hr^2_v(\strip{\omega})$. Since $v$ is strongly admissible, for
each 
$w\in \strip{\theta}$ we have 
$\ue^{\omega\abs{\bfs}}(\tau_w\vphi \cdot f)^\vee \in \Ell{1}(\R)$. Hence, 
\[ \vphi(z-w) f(z) = (\tau_w\vphi \cdot f)(z) = \int_\R (\tau_w\vphi \cdot f)^\vee(s)
\cdot \ue^{-\ui s z}\, \ud{s} \quad
\]
for all $z\in \strip{\omega}$. Observe that the right-hand side has a
continuous extension to $\cstrip{\omega}$. Hence, if we specialize
$\vphi = \gauss = \ue^{-\bfz^2}$ and $w=0$ we find that
\[ f(z) = \ue^{z^2} \int_\R  (\gauss f)^\vee(s)\ue^{-\ui s z}\, \ud{s} 
\]
extends continuously to $\cstrip{\omega}$. (Compare the proof of
Proposition \ref{hrr.p.Hr}.c for the case $\omega = 0$.) 

Furthermore, by Theorem  \ref{hrs.t.indep}.a,  
the function $v \ue^{\omega \abs{\bfs}} (\tau_\bfz \gauss \cdot f)^\vee$ is holomorphic
$\C \to \Ell{2}(\R)$
and bounded on each strip. In particular, it is uniformly continuous on $\cstrip{\omega}$.
Hence, to prove a) it suffices to show that the function
\[ H: \cstrip{\omega} \to \Ell{2}(\R),\qquad H(z) = \frac{1}{v}
\ue^{-\omega\abs{\bfs}}
\ue^{-\ui \bfs z}
\]
is bounded and uniformly continuous. Since $v$ is strongly admissible,
$\abs{H(z)} \le \frac{1}{v}$. This shows the boundedness. Moreover,
\begin{align*}
 \norm{H(z) - H(w)}_{\Ell{2}(\R)}^2 & = 
\int_\R \frac{\ue^{-2\omega \abs{s}}}{v(s)^2} \abs{\ue^{-\ui sz} -
  \ue^{-\ui sw}}^2\, \ud{s} 
\\ & \le \int_{\abs{s}\ge N} \frac{4}{v(s)^2}\, \ud{s} + \int_{-N}^N \abs{1 - \ue^{-\ui
    s (w-z)}}^2\, \ud{s} 
\end{align*}
for all $z,w\in \cstrip{\omega}$ and $N \in \N$. From this, the
uniform continuity follows readily.

\prfnoi
d)\ Now, coming back to the more general situation we find that 
\[ \vphi(z-w) f(z) = \int_\R (\tau_w\vphi \cdot f)^\vee(s)
\cdot \ue^{-\ui s z}\, \ud{s} 
\]
holds whenever $z\in \cstrip{\omega}$ and $z-w \in
\strip{\omega{+}\theta}$. This implies d).

\prfnoi
b)\ Fix $\omega > 0$ and $f\in \Hr^2_v(\strip{\omega})$. Then, as
shown above, $f$ extends continuously to $\cstrip{\omega}$. Obviously,
$(\tau_t\psi \cdot f)_{|\pm\omega} = \tau_t \psi_{|\pm \omega} \cdot
f_{|\pm \omega}$ for each $t\in \R$. 
By Corollary \ref{sob.c.We-char}, (i)$\dann$(ii), 
$\psi_{|\pm \omega} \in \Hr^2_v(\R)$ with 
\[  \norm{f}_{\Hr^2_v(\strip{\omega})} \eqsim
\norm{f_{|\omega}}_{\Hr^2_v(\R)} + \norm{f_{|-\omega}}_{\Hr^2_v(\R)}.
\]
Conversely, suppose that $f\in \Ha(\strip{\omega})$ and
its (distributional or almost everywhere) boundary functions $f_{|\pm
  \omega} = \lim_{y \nearrow \omega} f_{|\pm y}$ are
in $\Hr^2_v(\R)$. Then  
\[ (\tau_t \psi \cdot f)_{|y} = 
\tau_t \psi_{|y} \cdot f_{|y}
\to \tau_t \psi_{|\pm \omega} \cdot f_{|\pm \omega}
\]
in $\Ell{2}(\R)$ (cf. the discussion of the Paley--Wiener
theorem in Section \ref{s.sob}). By asumption, the latter
functions are  $\We^2_v(\R)$ uniformly in $t\in \R$.
As, obviously, $\tau_t\psi \cdot f\in \Har{2}(\strip{\omega})$,
Corollary \ref{sob.c.We-char}, (ii)$\dann$(i), applies
and yields $f\in \Hr^2_v(\strip{\omega})$ as desired.

\prfnoi
c)\ Take $f,g\in \Hr^2_v(\strip{\omega})$. 
Since $v$ is strongly admissible,
$\We_v^2(\strip{\omega})$ is, after scaling of norm,
a Banach algebra with respect to pointwise multiplication
(Proposition \ref{sob.p.We}.d). It follows that  
for each 
$t\in \R$, 
\[  \tau_t\psi^2 \cdot fg = (\tau_t\psi \cdot f)
\cdot (\tau_t \psi \cdot g) \in \We_v^2(\strip{\omega})
\]
with 
\[  \norm{\tau_t\psi^2 \cdot fg}_{\We_v^2(\strip{\omega})} \lesssim  
\norm{\tau_t\psi \cdot f}_{\We_v^2(\strip{\omega})}\,
\norm{\tau_t\psi \cdot g}_{\We_v^2(\strip{\omega})}.
\]
Taking the supremum over $t\in \R$ yields
$fg\in \Hr^2_v(\strip{\omega})$ with 
\[  \norm{fg}_{\Hr_v^2(\strip{\omega})} \lesssim  
\norm{f}_{\Hr_v^2(\strip{\omega})} \,
\norm{f}_{\Hr_v^2(\strip{\omega})}. 
\]
(Note that we employ Theorem \ref{hrs.t.indep}.a here.)
This proves the claim. 
\end{proof}

\vanish{
\begin{rem}
Theorem \ref{hrs.t.strongadm} shows in particular, that for $\alpha >
\frac{1}{2}$ and $v = (1 + \abs{\bfs})^\alpha$
the space $\Hr^2_v(\strip{\omega})$ coincides
with the classical 
\end{rem}
}

\section{Functional Calculus}\label{s.fcc}

In this section, $X$ always denotes a Banach space.  We fix $\omega
\ge 0$ and a {\em strongly} admissible
function $v: \R \to [1, \infty)$. Then $\We^2_v(\strip{\omega})$ and $\Hr^2_v(\strip{\omega})$
are Banach algebras continuously embedded into
$\Ha(\strip{\omega}) \cap \UCb(\cstrip{\omega})$
(Proposition \ref{sob.p.We}.d and Theorem  \ref{hrs.t.strongadm}).  

\medskip

A closed operator $A$ on $X$ is called \emdf{strip type} of
\emdf{height} $\omega$ if $\spec(A) \subseteq \cstrip{\omega}$ and
\[    \sup_{\abs{\im \lambda} > \omega'} \norm{R(\lambda,A)} < \infty
\]
for each $\omega'> \omega$. Such a  strip type operator $A$ admits a natural
(unbounded) holomorphic functional calculus 
\[  \Phi: \Ha[\cstrip{\omega}] \to \{ \text{closed operators on $X$}\}
\]
satisfying the axioms of an ``abstract functional calculus''
in the sense of \cite{Ha2020}. 
We briefly recall the construction, cf. also 
\cite[Ch. 4]{Ha2006} and \cite{Ha2018,HaHa2013}.

\medskip
The algebra of \emdf{elementary functions} is \label{f.elem-strip}
\[ \calE[\cstrip{\omega}] = \bigcup_{\theta > \omega}
\calE(\strip{\theta})
\]
where 
$\calE(\strip{\theta}) := \Ha(\strip{\theta}) \cap
\mathrm{H}^1(\strip{\theta})$ for $\theta > 0$. The 
\emdf{elementary calculus} 
\[ \Phi_0: \calE[\cstrip{\omega}] \to \BL(X)
\]
is given by 
\[ f(A) := \Phi_0(f) := \frac{1}{2\upi \ui}
\int_{\rand{\strip{\omega'}}} f(w) R(w,A)\, \ud{w}
\qquad (\omega < \omega'< \theta,\, f\in \calE(\strip{\theta})).
\]
The elementary calculus 
extends canonically to an (in general unbounded) functional calculus
 $\Phi$ with  $\Ha[\cstrip{\omega}] \subseteq \dom(\Phi)$.

\begin{defn}\label{fcc.d.ample} 
In the situation above, a set $\calD \subseteq \calE[\cstrip{\omega}]$
is called \emdf{ample} if 
\[   \cls{\spann} \bigcup_{f\in \calD} \ran(f(A)) =X.
\]
\end{defn}

Later, we shall need the following result.

\begin{lem}\label{fcc.l.ample}
Let $A$ be an operator of strip type $\omega \geq 0$ on a Banach space
$X$. Then the following statements hold:
\begin{aufzi}
\item For each $0 \neq \psi \in \calE[\cstrip{\omega}]$\quad   
$\dps
\bigcap_{t\in \R} \ker (\tau_t \psi)(A)\, =\, \{ 0\}$.

\item $\ker (\ue^{-\bfz^2})(A) = \{0\}$.

\item If $A$ is densely defined and $\theta > \omega$ 
then the sets
$\{ \vphi^2 \suchthat \vphi\in \Ho(\strip{\theta})\}$ 
and $\{ \tau_t \psi \suchthat t\in \R\}$ for given 
$0 \neq \psi \in \calE(\strip{\theta})$ are ample.
\end{aufzi} 
\end{lem}

\begin{proof}
Fix $\theta > \omega$ and $0 \neq \psi \in \calE(\strip{\theta})$. Without loss of generality we may assume that
\[
\|\psi\vert_\R \|_2^2\, =\, \int_\R \psi^\ast (t) \psi (t)\,  \ud t\, =\, \int_\R \psi^\ast (s-t) \psi (s-t)\, \ud t\, =\, 1 \qquad (s\in \R). 
\] Then, by holomorphy,
\begin{equation}\label{prel.id.IntIsOne}
\int_\R \psi^\ast (z-t) \psi (z-t)\, \ud t\, =\, 1 
\end{equation} for all $z\in \strip{\theta}$.
Set 
\[
X_\psi\coloneqq \cls{\spann}\bigcup_{t\in \R} \ran (\tau_t \psi)(A),
\] and let $|\Im \lambda|>\omega$. Then, for a suitable $\delta > 0$,
\begin{align*}
R(\lambda , A)^2 &=\, \frac{1}{2\pi \ui} \int_{\partial \strip{\delta}} \frac{1}{(\lambda - z)^2} R(z, A) \, \ud z\\
&=\, \frac{1}{2\pi \ui} \int_{\partial \strip{\delta}} \left(\int_\R \frac{\psi^\ast (z-t) \psi (z-t)}{(\lambda - z)^2}\, \ud t\right) R(z, A)\, \ud z\\
&=\, \int_\R \left(\frac{1}{2\pi \ui}\int_{\partial \strip{\delta}} \frac{\psi^\ast (z-t) \psi (z-t)}{(\lambda - z)^2} R(z, A)\, \ud z \right)\ud t\\
&=\, \int_\R (\tau_t \psi^\ast) (A)\,\,  (\tau_t \psi) (A)\,\,  R(\lambda, A)^2\, \ud t.
\end{align*} Here we have used identity (\ref{prel.id.IntIsOne}) in the second line and Fubini's theorem in the third line. It follows that
\[
\bigcap_{t\in \R}\ker  (\tau_t \psi) (A)\, \subseteq\, \ker R(\lambda, A)^2\, =\, \{ 0\}, 
\] and 
$ \dom(A^2) = \ran R(\lambda, A)^2 \subseteq X_\psi$.
Hence we obtain {\rm a)} and the second part of {\rm c)}.

\prfnoi
Now we specialize $\psi := \ue^{-\bfz^2}$. Then, obviously,
\[   \{ \tau_t \psi \suchthat t\in \R\} \subseteq 
\{ \vphi^2 \suchthat \vphi \in \Ho(\strip{\theta})\}.
\]
This yields the second part of {\rm c)}. 
For {\rm b)} write $\tau_t \psi = f_t \psi$ with $f_t = \ue^{-t^2}
\ue^{2t\bfz}$. Note that $\tau_s \psi \cdot f_t \in \calE[\cstrip{\omega}]$
for each $s\in \R$. It follows that 
\[ (\tau_s \psi)(A) (\tau_t \psi)(A) = 
(\tau_s \psi \cdot \tau_t \psi)(A)= 
(\tau_s \psi \cdot f_t)(A) \psi(A) \qquad (s, t\in \R). 
\]
By a), this implies $\ker \psi(A) = \{0\}$.
\end{proof}

Returning to our main object of study, 
the question arises whether the holomorphic
calculus can be extended to
a functional calculus defined on $\Hr^2_v(\strip{\omega})$.
Observe that  the  ``algebraic extension'' method will not work
because the functions
in the algebra $\Ha[\cstrip{\omega}]$ are all defined on open strips 
containing $\cstrip{\omega}$ and hence a quotient of such functions
is necessarily meromorphic on such a larger strip. In contrast, a
function in $\Hr^2_v(\strip{\omega})$ in general 
is defined on $\cstrip{\omega}$
(sharp) and may not have a meromorphic extension to a larger  strip. 

Therefore, in order to obtain a $\Hr^2_v(\strip{\omega})$-calculus
for $A$, one looks for a suitable {\em topological extension} of the
natural holomorphic calculus. 
To this end, note  that restriction
induces an embedding 
\[    
\calE[\cstrip{\omega}] 
\hookrightarrow \bigcup_{\theta > \omega} \Har{2}(\strip{\theta})
\subseteq \We^2_v(\strip{\omega})
\]
and hence we may consider $\Phi_0$ as a non-fully defined operator 
$\We^2_v(\strip{\omega}) \supseteq \dom(\Phi_0) 
\to \BL(X)$. Suppose that  if $\Phi_0$ is closable. Then
its closure 
\[ \Phi_{\We^2_v} := \cls{\Phi_0} : \dom(\cls{\Phi_0}) \to \BL(X)
\]
is an algebra homomorphism, called  the \emdf{elementary Sobolev
  calculus}. This is a ``topological extension'' in the sense of \cite{Ha2020}.

As usual, one  may extend $\Phi_{\We^2_v}$ by the 
algebraic method to a still larger algebra of functions on $\cstrip{\omega}$. By general
theory (see e.g. \cite[Cor. 9.2]{Ha2020}) 
this extension constitutes an extension
the natural holomorphic calculus $\Phi$. In effect, there is no danger in keeping
the symbol $\Phi$ also for this second extension and always
write $f(A)$ instead of $\Phi(f)$.

We are now heading for a condition on $A$ guaranteeing that $\Phi_0$ is closable
and that $\Hr^2_v(\strip{\omega})$ is contained in 
the domain of the resulting calculus. (In this case, we say that 
$A$ \emdf{admits an (unbounded) $\Hr^2_v(\strip{\omega})$-calculus}.

\begin{defn}
Let $A$ be a strip-type operator of height $\omega$. A function $h\in
\calE[\cstrip{\omega}]$ is called a \emdf{$\We^2_v(\strip{\omega})$-regularizer} for
$A$ if
there is a constant $c = c(h)$ such that 
\[     \norm{(\psi h)(A)} \le c \norm{\psi}_{\We^2_v(\strip{\omega})}
\qquad (\psi \in \Ho[\cstrip{\omega}]).
\]
The operator $A$ is called \emdf{almost
  $\We^2_v(\strip{\omega})$-regular}, if  
\[   \bigcap \{ \ker h(A) \suchthat \text{$h$ is 
a $\We^2_v(\strip{\omega})$-regularizer for $A$} \} = \{0\}.
\]
\end{defn}

In the terminology of \cite{Ha2020}, $A$ is almost
  $\We^2_v(\strip{\omega})$-regular iff the set of 
$\We^2_v(\strip{\omega})$-regularizers for $A$ is an {\em anchor set} 
within the natural holomorphic calculus.

\begin{prop}\label{fcc.p.W2vreg-char} 
Let $A$ be a strip-type operator of height $\omega$. Then, for a
function
$h\in \calE[\strip{\omega}]$ the following statements are equivalent:
\begin{aufzii}
\item $h$ is $\We^2_v(\strip{\omega})$-regularizer for $A$\textup{;}

\item $\exists\, c\geq 0, \, \theta> \omega\colon\, \, \forall \, \psi \in \Ho (\strip{\theta})\colon \quad \| (\psi h)(A)\|\, \leq\, c\| \psi\|_{\We_v^2 (\cstrip{\omega})}$\textup{;}

\item $\forall\, x\in X, \, x'\in X'\colon \quad s\mapsto \tfrac{\ue^{-\omega
      |s|}}{v} \langle (\ue^{-\ui s \bfz} h)(A)x, x' \rangle \, \in \, \Ell{2}(\R)$.

\end{aufzii}
\end{prop}

\begin{proof}
Clearly {\rm (i)} implies {\rm (ii)}. 
For the proof of the remaining implications we abbreviate
\[  U_s^h :=   (\ue^{-\ui s \bfz}h)(A) \in \BL(X) \qquad (s\in \R)
\]
and  note that
\begin{equation}\label{fcc.id.UnbddSFC}
(fh)(A)\, =\, \int_\R f^\vee (s) U_s^h\, \ud s
\end{equation} 
for all $f\in \calE [\strip{\omega}]$, the integrand being 
continuous and integrable in $\BL(X)$. (Write out $U_s^h = (\ue^{-\ui s
  \bfz}h)(A)$ as a  Cauchy-integral over $\rand\strip{\omega'}$, where
$\omega'> \omega$ is sufficiently close to $\omega$.)

Now suppose that (iii) holds. By the closed graph theorem there is a constant $c\geq 0$ with
\[ 
\Big\| \frac{\ue^{-\omega |\bfs|}}{v} \langle U_\bfs^h x, x' \rangle
\Big\|_2 \, \leq\, c \| x\| \, \|x'\| \qquad (x\in X,\, x'\in X').
\]
Let $f\in \calE [\strip{\omega}]$,
$x\in X$, and $x'\in X'$. Then, by \eqref{fcc.id.UnbddSFC}
and Cauchy--Schwarz
\begin{align*}
\left| \langle (fh)(A)x, x'\rangle \right| =&\, \left| \int_\R f^\vee (s) \langle U_s^h x, x'\rangle\, \ud s \right|\\
=&\, \left| \int_\R \big( v(s) \ue^{\omega |s|} f^\vee (s)\big) \cdot \Big( \frac{\ue^{-\omega |s|}}{v(s)} \langle U_s^h x, x' \rangle\Big)\, \ud s \right|\\[0.2cm]
\leq&\, c \| f\|_{\We_v^2(\strip{\omega})}\| x\|\, \|x'\|.
\end{align*} 
This yields (i).

Now suppose that (ii) holds. By Lemma \ref{adm.l.AWalwaysSmooth} we may suppose that
$v\in \uC^\infty (\R)$. 
Let $\eta\in \Cc^\infty (\R)$ and let $x\in X$, and $x'\in
X'$. Abbreviate $f = \Fourier\bigl( \frac{\eta}{\cosh(\omega \bfs)
  v}\bigr) \in \Fourier(\Test(\R)) \subseteq
\Ho(\strip{\theta})$. Then 
\begin{align*}
\left| \int_\R \eta(s) \frac{\langle U_s^h x, x'\rangle}{\cosh(\omega s) v(s)} \, \ud s \right| =&\, \left| \int_\R f^\vee(s)\langle U_s^h x, x' \rangle \, \ud s\right|
=\, \big| \langle (fh)(A)x, x' \rangle \big|\\
\lesssim &\, \norm{f}_{\We_v^2 (\strip{\omega})} \, \| x\|\, \| x'\|
\lesssim \, \| \eta \|_2 \, \| x\|\, \| x'\|,
\end{align*} 
As $\Test(\R)$ is dense in $\Ell{2}(\R)$, (iii) follows. 
\end{proof}

Next, we show that $A$ admits an unbounded
$\Hr^2_v(\strip{\omega})$-calculus
if  $A$ is almost $\We^2_v(\strip{\omega})$-regular. 

\begin{prop}\label{fcc.p.ubdhc} 
Let $A$ be a strip type operator of height $\omega\ge 0$
on a Banach space $X$ and let $v: \R \to [1, \infty )$ be strongly
admissible.  Suppose that $A$ is almost
$\We^2_v(\strip{\omega})$-regular. 
Then the following statements hold:
\begin{aufzi}
\item $A$ admits an (possibly unbounded) $\Hr^2_v(\strip{\omega})$-calculus. 
\item For each $f\in \We_v^2(\strip{\omega})$ and each
  $\We^2_v(\strip{\omega})$-regularizer
$h\in \calE[\cstrip{\omega}]$ one has $(fh)(A) \in \BL(X)$ and
\[ 
\langle (fh)(A)x, x' \rangle\, =\, \int_\R f^\vee (s) \langle U^h_sx, x' \rangle \, \ud s
\]
for all $x\in X $ and $x' \in X'$, where $U^h_s= (\ue^{-\ui s \bfz} h)(A)$ for $s\in \R$.
\end{aufzi}
\end{prop}

\begin{proof}
{\rm a)}, first part:\  We start with showing that the elementary calculus $\Phi_0$
is $(\We^2_v(\strip{\omega}), \BL(X))$-closable. 
Let $\psi_n \in \calE[\cstrip{\omega}]$ with $\psi_n \to 0$
in $\We_v^2 (\strip{\omega})$ and $\psi_n (A) \to T$ in $\BL (X)$. We
need to show that $T=0$. Fix a $\We^2_v(\strip{\omega})$-regularizer
$h\in \calE[\cstrip{\omega}]$. Then
\[    \norm{h(A)\psi_n(T)} = \norm{(\psi_n h)(A)} \lesssim
\norm{\psi_n}_{\We^2_v(\strip{\omega})} \to 0 \qquad (n \to \infty).
\]
This implies $h(A)T = 0$. Since $\bigcap_h \ker(h(A)) = \{0\}$ it follows
that $T= 0$ as desired.

\prfnoi
{\rm b)}\  Fix $f\in \We_v^2 (\strip{\omega})$ and $h\in
\calE[\cstrip{\omega}]$ as before. 
Choose $\psi_n\in \calE[\cstrip{\omega}]$ with $\psi_n \to f$ in
$\We_v^2(\strip{\omega})$. 
Then, $\psi_n h\to fh$ in $\We_v^2 (\strip{\omega})$. Since $h$ is a
$\We^2_v(\strip{\omega})$-regularizer, the sequence $((\psi_nh)(A))_n$ is operator norm
Cauchy and therefore converges to some operator $T\in \BL(X)$. It
follows that $f\in \dom(\Phi_{\We^2_v})$ and $f(A) = T$. Furthermore,
\[
(\psi_nh)(A)\, =\, \int_\R \psi_n^\vee (s) U_s^h \, \ud s
\qquad (n \in \N)
\] 
by (\ref{fcc.id.UnbddSFC}). The claim may now follows from 
Proposition \ref{fcc.p.W2vreg-char}.(iii).

\prfnoi
{\rm a)}, second part:\ Let $f\in \Hr^2_v(\strip{\omega})$. Then
$\ue^{-\bfz^2} f \in \We^2_v(\strip{\omega})$. Since
$(\ue^{-\bfz^2})(A)$ is injective, we may suppose without loss of
  generality that $f\in \We_v^2(\strip{\omega})$. By b) one has
$hf \in \dom(\Phi_{\We^2_v})$ for each
$\We_v^2(\strip{\omega})$-regularizer $h$, and since these functions
$h$ form an anchor set (by hypothesis), $f(A)$ is defined via the
algebraically extended  calculus.
\end{proof}

\begin{rem}\label{fcc.r.KrWe-ubdhc}
The notion of almost $\We^2_v(\strip{\omega})$-regularity is coined in
order to generalize the setting of 
Kriegler and Weis in \cite{Kr2009}, \cite{KrWe2017}, and
\cite{KrWe2018}. There, $0$-strip type operators $A$ and strongly
admissible functions $v =  (1 + \abs{\bfs})^\beta$ with
$\beta > \tfrac{1}{2}$ are considered; in addition it is required that
there is $h\in \Ho[\R]$ such that
$h(A)$ is injective, and
$\tfrac{1}{v_\beta}U_\bfs h(A)$ has weakly square integrable orbits.
(See also \cite[(3.9)]{KrWe2017} for a sectorial formulation.) 
\end{rem}

We shall see in the following section that 
if $-\ui A$ generates a $C_0$-group then
almost $\We_v^2(\strip{\omega})$-regularity, and in fact
a stronger property, is automatic.  The following example
shows that  there are almost $\We_v^2(\strip{\omega})$-regular
operators $A$ such that $-\ui A$ does not generate a group

\begin{exa}
Let $p\in (1, \infty)$, $p\neq 2$, and
let $-\Laplace$ denote 
the negative Laplacian on $\Ell{p}(\R)$. Then, $-\Laplace$ is of strip
type $0$ (see \cite[p. 236]{Ha2006}) 
but the operators $\ue^{\ui s\Laplace}$ are unbounded for all $s\neq
0$ (\cite[Proposition 8.3.8]{Ha2006}). In particular, 
$-\ui \Laplace$ does not generate a $C_0$-group on $\Ell{p}(\R)$. 
However, by the Mikhlin multiplier theorem (\cite[Theorem E.6.2.b)]{Ha2006}),
\begin{align*}
&\big\| \ue^{\ui t\Laplace} R(\lambda, -\Laplace) \big\| \,\lesssim \, \max_{k=0, 1} \Big\|\,  \bfs^k \Big( \frac{\ue^{\ui t\bfs^2}}{\lambda + \bfs^2} \Big)^{(k)} \Big\|_{\infty, \R\setminus \{ 0\}}\\[0.2cm]
& \leq\, \Big\|\frac{\ue^{\ui t\bfs^2}}{\lambda + \bfs^2}\Big\|_{\infty, \R\setminus\{ 0\}} + \Big\| \, \ui t\cdot \frac{2\bfs^2 \ue^{\ui t\bfs^2}}{\lambda + \bfs^2} - \frac{2\bfs^2 \ue^{\ui t \bfs^2}}{(\lambda + \bfs^2)^2}  \Big\|_{\infty, \R\setminus\{ 0\}}\, \lesssim \, 1 + |t|
\end{align*} 
for all $t\in \R$ and $|\Im \lambda| > 0$. Hence, $h\coloneqq (\lambda
- \bfz)^{-2}$ is a $\We^2_v(\R)$-regularizer whenever 
$v:\R \to [1, \infty )$ is an admissible function with
$\tfrac{1+|\bfs|}{v}\in \Ell{2}(\R)$. A fortiori, $-\Laplace$ is
almost $\We^2_v(\R)$-regular for such $v$.  
\end{exa}

\section{Bounded Sobolev Calculus}\label{s.bsc}

Still,  $X$ is a Banach space, $\omega\ge 0$ and 
$v: \R \to [1, \infty)$ is strongly admissible. 

\begin{lem}\label{bsc.l.res}
Let $\lambda \in \C$ with 
$\abs{\im \lambda} > \omega$. Then  
\[    (\lambda - \bfz)^{-1}  \in \We^2_v(\strip{\omega}).
\]
The mapping 
\[ \C \ohne \cstrip{\omega} \to \We^2_v(\strip{\omega}),\qquad \lambda \mapsto
r_\lambda :=  (\lambda -\bfz)^{-1} 
\]
is holomorphic and  for each $\omega'> \omega$ one has
\[    \norm{(\lambda - \bfz)^{-1}}_{\We^2_v(\strip{\omega})} 
\lesssim 
(\abs{\im \lambda} - \omega')^{-\frac{1}{2}} \qquad (\abs{\im \lambda}
> \omega').
\]
\end{lem}

\begin{proof}
Fix $\omega'> \omega$. Then $\Har{2}(\strip{\omega'}) \subseteq
\We^2_v(\strip{\omega})$ continuously. It is elementary to
verify
\[    \norm{ (\lambda - \bfz)^{-1}}_{\Har{2}(\strip{\omega'})} =
\sqrt{\upi}  (\abs{\im \lambda} - \omega')^{-\frac{1}{2}}.
\]
This yields the first assertion and the claimed norm
estimate. To prove holomorphy of the mapping 
$\lambda \mapsto r_\lambda$
it suffices, by the Arendt--Nikolski theorem,  
to  test this against a point-separating 
set $W$ of bounded linear functionals on
$\We^2_v(\strip{\omega}$). Obviously, one 
can take  $W = \{ \delta_w  \suchthat w\in \cstrip{\omega}\}$. 
\end{proof}

\begin{defn}\label{fcc.def.DefBddSobCalc}
A possibly unbounded operator  $A$  on $X$ 
is said to \emdf{ have a bounded $\We^2_v(\strip{\omega})$-calculus}
if there exists a bounded algebra
homomorphism
\[ \Psi: \We^2_v(\strip{\omega}) \to \BL(X)
\]
such that $\Psi(r_\lambda) = R(\lambda, A)$ for 
some/all $\lambda \in \C\ohne \cstrip{\omega}$.
\end{defn}

Observe that if $\Psi$ is as in the definition, then 
$\lambda \mapsto \Psi(r_\lambda)$ is a pseudo-resolvent
and hence there is a linear (possibly multivalued)
closed linear operator $A$ on $X$ with $R(\lambda, A) = 
\Psi(r_\lambda)$.  This operator is uniquely determined
by each single operator $\Psi(r_\lambda)$
\cite[Prop. A.2.4]{Ha2006}.

\begin{prop}\label{bsc.p.bsc-char}
For an operator $A$ on a Banach space $X$ the following
assertions are equivalent:
\begin{aufzii}
\item $A$ has a bounded $\We^2_v(\strip{\omega})$-calculus $\Psi$;
\item $A$ is of strip type $\omega$ and
there is $\theta > \omega$ and a constant $c\geq 0$ with
\[
\| \psi(A)\| \, \leq \, c \| \psi\|_{\We_v^2 (\strip{\omega })}\qquad \text{for all $\psi\in\Ho(\strip{\theta})$}.
\]
\end{aufzii}  
In this case, $A$ is almost $\We^2_v(\strip{\omega})$-regular and  
$\Psi$ coincides 
with the topological extension $\Phi_{\We^2_v}$
of the natural holomorphic calculus for $A$. 
\end{prop}

\begin{proof}
(ii)$\dann$(i):  Suppose (ii) and 
pick $h\in \calE[\cstrip{\omega}]$. Then 
\[ \norm{(\psi h)(A)} = \norm{h(A)} \norm{\psi(A)}\le c \norm{h(A)}
\norm{\psi}_{\We^2_v(\strip{\omega})} \quad (\psi \in \Ho(\strip{\theta})).
\]
Hence, by Proposition \ref{fcc.p.W2vreg-char}, 
$h$ is a $\We_v^2(\strip{\omega})$-regularizer for $A$. 

Next, let $f\in \We^2_v(\strip{\omega})$ and let $(\psi_n)_n$ be a
sequence in $\Ho(\strip{\theta})$ with $\psi_n \to f$ in $\We^2_v(\strip{\omega})$.
By (ii), $(\psi_n(A))_n$ is a Cauchy sequence in $\BL(X)$
and hence there is $T\in \BL(X)$ with $\psi_n(A) \to T$. It follows that
$f\in \dom(\Phi_{\We^2_v})$, $f(A) = T$, and  
$\norm{f(A)}\le c \norm{f}_{\We^2_v(\strip{\omega})}$. 
This gives (i) with $\Psi(f) = f(A)$.

\prfnoi
(i)$\dann$(ii): Suppose (i) holds. Then it follows from Lemma
\ref{bsc.l.res} that $A$ is of strip type $\omega$. Fix
$\theta > \omega$ and $f\in \calE(\strip{\theta})=
\Har{1}(\strip{\theta}) \cap \Har{\infty}(\strip{\theta})$.
Then, with $\omega < \omega' < \theta$
\[ f  = \frac{1}{2\upi \ui} \int_{\rand{\strip{\omega'}}} 
f(w) r_w \, \ud{w}
\]
as an integral in $\We^2_v(\strip{\omega})$, again by Lemma
\ref{bsc.l.res}. It follows that 
\[   \Psi(f) = \frac{1}{2\upi \ui} \int_{\rand{\strip{\omega'}}} 
f(w) \Psi(r_w) \, \ud{w} 
= 
\frac{1}{2\upi \ui} \int_{\rand{\strip{\omega'}}} 
f(w) R(w,A) \, \ud{w} = f(A).
\]
This proves (ii) and $\Psi = \Phi_0$ on
$\calE[\cstrip{\omega}]$. From this it follows
that $\Psi= \Phi_{\We^2_v}$ as claimed.
\end{proof}

The next proposition, a generalization of \cite[Lemma 3.7]{KrWe2017},  
is similar to Proposition \ref{fcc.p.W2vreg-char}, 
and yields a characterization of group generators with a bounded Sobolev calculus.

\begin{prop}\label{fcc.p.bsobc-char}  
Let $-\ui A$ be the generator of a $C_0$-group $(U_s)_{s\in \R}$ on a
Banach space $X$, let $\omega \geq 0$ and let 
$v: \R \to [1, \infty)$ be strongly  admissible. Then the following statements are equivalent:
\begin{aufzii}
\item The mapping $\frac{\ue^{-\omega |\bfs|}}{v} U_\bfs$ has weakly square integrable orbits, i.e., for each $x\in X$ and $x' \in X'$,
\[
\frac{\ue^{-\omega |\bfs|}}{v} \langle U_\bfs x, x' \rangle\, \in \, \Ell{2}(\R ).
\]
\item $A$ has a bounded $\We_v^2 (\strip{\omega})$-calculus.
\end{aufzii} In that case, 
\[ 
\langle f(A)x, x' \rangle\, =\, \int_\R f^\vee (s) \langle U_s x, x'
\rangle\, \ud s \qquad (x\in X,\, x'\in X')
\]
for all $f\in \We_v^2(\strip{\omega})$.
\end{prop}

\begin{proof}
Conversely, suppose that (ii) holds. Then we obtain (i) as in the
proof of Proposition \ref{fcc.p.W2vreg-char} but 
with $h = \car$. 

\prfnoi
Suppose that (i) holds. By the closed graph theorem we find $c > 0$
with 
\[     \norm{v^{-1} \ue^{-\omega \abs{\bfs}} \dprod{U_\bfs x}{x'}
}_{\Ell{2}(\R)} 
\le c \norm{x} \norm{x'} \qquad (x \in X,\, x'\in X').
\]
Next, fix  $\theta > \omega$. Then, by Cauchy--Schwarz
\[   \int_\R \ue^{-\theta \abs{s}} \abs{\dprod{U_s x}{x'}} \, \ud{s}
\lesssim \norm{x} \norm{x'}  \qquad (x\in X, \, x'\in X').
\]
It now follows from \cite[Thm. 5.1]{Glu2015} that $A$ is of strip type $\le \omega$.

Fix $\theta > \theta(U)$. Then for $\psi \in
\Ho(\strip{\theta})$ one has
\[    \psi(A) = \int_\R \psi^\vee(s) U_s \, \ud{s}.
\]
Indeed, the integral converges strongly; and multiplying the identity
by  $(\lambda - A)^{-2}$ from the left reduces the matter to \eqref{fcc.id.UnbddSFC} for $h = (\lambda - \bfz)^{-2}$, 
which is already known to be true. Finally, as in the 
proof of Proposition \ref{fcc.p.W2vreg-char}, use Cauchy--Schwarz to estimate
\[   \abs{\dprod{\psi(A)x}{x'}} \le
\norm{\psi}_{\We^2_v(\strip{\omega})}
\, \norm{v^{-1} \ue^{-\omega \abs{\bfs}} \dprod{U_s x}{x'}
}_{\Ell{2}(\R)}
\le c \norm{\psi}_{\We^2_v(\strip{\omega})} \norm{x} \norm{x'}
\]
for $x\in X,\, x'\in X'$.  This gives (ii). 
The additional statement follows  by approximation. 
\end{proof}

\begin{exa}
Let $-\ui A$ be the generator of a $C_0$-group $(U_s)_{s\in \R}$ such
that $\norm{U_\bfs} \le  v_u \ue^{\omega \abs{\bfs}}$ for some
measurable function  $v_U: \R \to (0, \infty)$. If $v$ is admissible
with $v_U/v \in \Ell{2}(\R)$, then $A$ is
$\We^2_v(\strip{\omega})$-regular. This follows from Proposition
\ref{fcc.p.bsobc-char} (but can quite easily be seen directly as well). 
\end{exa}

In order to extend a bounded Sobolev calculus
to a bounded Hörmander calculus (which is the subject of
the next section), we need the following 
``convergence lemma''.

\begin{thm}[Convergence Lemma]\label{fcc.t.conlem}
Let $A$ be an operator on $X$ with 
a bounded $\We_v^2
(\strip{\omega})$-calculus, where
$\omega \geq 0$ and 
$v: \R \to [1, \infty)$ is strongly  admissible.
Suppose  that $f\in \Hr_v^2 (\strip{\omega})$ and that
$(f_n)_n$ is a sequence in $\Hr_v^2 (\strip{\omega})$ with the following properties:
\begin{aufziii}
\item There is an ample set 
$\calD \subseteq \Ho[\strip{\omega}]$ with 
\[
\norm{\psi (f_n -f )}_{\We_v^2 (\strip{\omega})} \, \xrightarrow{n\to \infty}\, 0 \qquad (\psi \in \calD);
\]
\item $\sup_{n\in \N} \norm{f_n(A)}  < \infty$.
\end{aufziii} Then  $f(A)\in\BL (X)$  
and $f_n (A) \to f(A)$ strongly on $X$.
\end{thm}

\begin{proof}
Observe that
\[
\norm{ (f_n(A) - f(A))\psi (A)}\, =\, \norm{ (\psi(f_n -f))(A)}
\, \lesssim_A\, \| \psi (f_n -f)\|_{\We_v^2 (\strip{\omega})}
\, \xrightarrow{n\to \infty}\, 0
\] 
for all $\psi \in \calD$. As $\calD$ is ample
and the sequence $(f_n (A))_n$ is uniformly bounded, the claim follows. 
\end{proof}

In the next lemma we show that for a given function 
$f\in \Hr_v^2 (\strip{\omega})$ one can 
always find a bounded sequence $(f_n)_n$ in $\Hr_v^2(\strip{\omega})$ such that condition {\rm 1)} of Theorem
\ref{fcc.t.conlem} is satisfied. 
Moreover, this sequence may be chosen in 
$\Ho(\strip{\theta})$ for any given  $\theta>\omega$.

\begin{lem}\label{bsc.l.appHo}
Let $\theta >\omega \geq 0$, and let $v: \R \to [1, \infty)$ be strongly admissible. 
Then there is a constant $K\geq 0$ with the following property: For each $f\in \Hr_v^2 (\strip{\omega})$ there is a sequence $(f_n)_n$ in $\Ho (\strip{\theta})$ with 
\begin{equation}\label{bsc.eq.appHo}
\sup_{n\in \N} \| f_n\|_{\Hr_v^2 (\strip{\omega})} \, \leq\,  K \| f\|_{\Hr_v^2 (\strip{\omega})}
\end{equation} and 
\[
\| \psi_1\psi_2 (f_n -f )\|_{\We_v^2 (\strip{\omega})}\, \xrightarrow{n\to \infty}\, 0
\] 
for all $\psi_1, \psi_2\in \Ho[\cstrip{\omega}]$.
\end{lem}

\begin{proof}
Define $\vphi_n := \ue^{-\frac{1}{n} \bfz^2} \in \Ho(\strip{\theta})$
  for $n\in \N$. Then 
\[    \sup_{n \in \N} \norm{\vphi_n}_{\Hr^2_v(\strip{\omega})}
  \lesssim \sup_{n \in \N} \norm{\vphi_n}_{\Ha(\strip{\theta})} <
  \infty.
\]
Furthermore, $\vphi_n \to \car$ uniformly on compacts.
Hence, for each $\psi \in \Ho(\strip{\theta})$
\[   \norm{\psi (\car - \vphi_n)}_{\We^2_v(\strip{\omega})}
\lesssim \norm{\psi (\car - \vphi_n)}_{\Har{2}(\strip{\theta})} \to 0.
\]
Now, fix $f\in \Hr_v^2 (\strip{\omega})$ and set
\[
g_n \coloneqq \vphi_n f \in \We^2_v(\strip{\omega})\qquad (n\in \N).
\] 
Then 
\[
\sup_{n\in \N} \|g_n\|_{\Hr_v^2 (\strip{\omega})}
\, \lesssim_v \, 
\big(\sup_{n\in \N} \|\vphi_n\|_{\Hr_v^2 (\strip{\omega})}\big) \| f\|_{\Hr_v^2 (\strip{\omega})}
\, <\, \infty,
\] 
and 
\[
\norm{\psi_1 \psi_2 (g_n - f)}_{\We_v^2 (\strip{\omega})}
\, \lesssim_v\, 
\norm{\psi_1 (\phi_n - \car)}_{\We_v^2 (\strip{\omega})} \, 
\norm{\psi_2 f}_{\We_v^2 (\strip{\omega})}
\, \xrightarrow{n\to \infty}\, 0
\]
 for all $\psi_1, \psi_2 \in \Ho[\strip{\omega}]$ by the
 considerations above. 
Finally, recall that $\Ho (\strip{\theta})$ is dense in
$\We_v^2(\strip{\omega})$, which is why one can choose 
$f_n \in \Ho (\strip{\theta})$ with
\[
\| g_n - f_n \|_{\We_v^2 (\strip{\omega})}\, \leq\, \frac{1}{n} \|
f\|_{\Hr_v^2 (\strip{\omega})} \qquad (n\in \N).
\] 
It follows that 
\begin{align*}
\| f_n\|_{\Hr_v^2 (\strip{\omega})}
 &= \, \sup_{t\in \R} \|\tau_t \gauss\cdot f_n\|_{\We_v^2 (\strip{\omega})}\\
 &\leq\, \sup_{t\in \R} \big( \|\tau_t \gauss\cdot (f_n - g_n)\|_{\We_v^2 (\strip{\omega})} + \|\tau_t \gauss\cdot g_n\|_{\We_v^2 (\strip{\omega})} \big)\\
 & \lesssim_v\, \tfrac{1}{n} \| \gauss\|_{\We_v^2 (\strip{\omega})}\| f\|_{\Hr_v^2 (\strip{\omega})} + \| g_n\|_{\Hr_v^2 (\strip{\omega})}\\[0.2cm]
& \lesssim_v \,\Big( \| \gauss\|_{\We_v^2 (\strip{\omega})} + 
\sup_{k\in \N}\| \vphi_k\|_{\Hr_v^2 (\strip{\omega})}\Big) \|f\|_{\Hr_v^2 (\strip{\omega})}.
\end{align*} 
Therefore, we may set 
\[
K \coloneqq \| \gauss\|_{\We_v^2 (\strip{\omega})} + 
\sup_{k\in \N}\| \phi_k\|_{\Hr_v^2 (\strip{\omega})}
\] (times a constant depending only on $v$) to establish (\ref{bsc.eq.appHo}). One readily verifies that
\[
\| \psi_1\psi_2 (f_n - f) \|_{\We_v^2 (\strip{\omega})}\, \xrightarrow{n\to \infty}\, 0
\] 
for all $\psi_1, \psi_2 \in \Ho[\cstrip{\omega}]$, and this concludes the proof.
\end{proof}

\section{Bounded Hörmander Calculus}\label{s.bhc}

As in the previous sections, $X$ always denotes a Banach space,
$\omega \in \R_{\ge 0}$ and $v: \R \to [1, \infty)$ is strongly
admissible.

\begin{defn}
An operator $A$ on $X$ is said to have a \emdf{bounded
  ($\gamma$-bounded, $R$-bounded)
    $\Hr^2_v(\strip{\omega})$-calculus} if 
there is a bounded ($\gamma$-bounded, $R$-bounded)  unital algebra homomorphism 
\[  \Psi: \Hr^2_v(\strip{\omega}) \to \BL(X)
\]
such that $\Psi(r_\lambda) = R(\lambda, A)$ for one/all $\lambda \in
\C\ohne \cstrip{\omega}$.
\end{defn}

Suppose that $A$ has a bounded $\Hr^2_v(\strip{\omega})$-calculus. Then,
as  $\We^2_v(\strip{\omega}) \subseteq \Hr^2_v$ continuously, $A$ also
has a bounded Sobolev calculus and hence is of strip type
$\omega$. Moreover,
the Sobolev calculus extends (topologically) the natural holomorphic calculus
(Proposition \ref{bsc.p.bsc-char}). 
In particular, $\Psi( \ue^{-\bfz^2} ) =
(\ue^{-\bfz^2})(A)$ is injective (Lemma \ref{fcc.l.ample}). 
As each $f\in \Hr^2_v(\strip{\omega})$ can be
be written as the quotient $f = (\ue^{-\bfz^2}
f)/ \ue^{-\bfz^2}$, we conclude that $\Psi$ coincides (on its domain) with the algebraic
extension of the $\We^2_v(\strip{\omega})$-calculus. In particular, 
if $A$ has a bounded $\Hr^2_v(\strip{\omega})$-calculus, then it is
unique.

\medskip

As a corollary of Theorem \ref{fcc.t.conlem} and Lemma 
\ref{bsc.l.appHo} we obtain a characterization for when an operator admits a bounded $\Hr_v^2 (\strip{\omega})$-calculus.

\begin{thm}\label{bhc.t.bhc-char} 
Let $A$ be a densely defined operator on $X$ and let $\theta >
\omega$.
Then the following statements are equivalent.
\begin{aufzii}
\item $A$ has a bounded ($\gamma$-bounded, $R$-bounded)
 $\Hr_v^2(\strip{\omega})$-calculus.
\item $A$ is of strip type $\omega$ and the set
\[  \{ \psi(A) \suchthat \psi \in \Ho (\strip{\theta}),\, 
\norm{\psi}_{\Hr_v^2 (\strip{\omega})} \le 1\} \subseteq \BL(X)
\]
is bounded ($\gamma$-bounded, $R$-bounded).
\end{aufzii}
\end{thm}

\begin{proof}
The implication (i)$\dann$(ii) is clear from the remarks preceding
this theorem.  For the converse, 
suppose that {\rm (ii)} holds. Then, in particular, $A$ has a bounded
$\We_v^2 (\strip{\omega})$-calculus. 
Choose $K$ as in Lemma \ref{bsc.l.appHo}.

To a  given $f\in \Hr^2_v(\strip{\omega})$ with 
$\norm{f}_{\Hr(\strip{\omega})} \le 1$ we can then 
find a sequence $(f_n)_n$ in 
$\Ho(\strip{\theta})$ with $\sup_n \norm{f_n}_{\Hr^2_v} \le K$ and
$\norm{\psi^2(f_n - f)}_{\We^2_v} \to 0$ for all $\psi \in
\Ho[\cstrip{\omega}]$. By (ii), $\sup_n \norm{f_n(A)} < \infty$. 
Since the set $\calD := \{ \psi^2 \suchthat 
\psi \in \Ho[\cstrip{\omega}]\}$ is ample, 
Theorem \ref{fcc.t.conlem} is applicable and yields $f(A) \in \BL(X)$ 
and $f_n(A) \to f(A)$ strongly. Hence, the operator
$f(A)$ lies in the strong closure of the 
bounded $(\gamma$-bounded, $R$-bounded) set
\[   \{ \psi(A) \suchthat \psi \in \Ho(\strip{\theta}),\,
\norm{\psi}_{\Hr^2_v(\strip{\omega})} \le K\}.
\]
This yields (i). 
\end{proof}

The following is our first result on bounded Hörmander calculus. It
extends \cite[Theorem 2.2]{GCMMST2001} to general Banach spaces. The
proof is based on the approach in the classical works of
Meda \cite{Me1990} and Cowling and Meda \cite{CoMe1993}, see also
\cite[p. 944 - 946]{CarDra2017}.

\begin{thm}\label{bhc.t.meda}
Let $-\ui A$ be a densely defined operator with 
a bounded $\Ha(\strip{\theta})$-calculus for some $\theta > 0$ on a
Banach space $X$. Suppose 
\[  \norm{ \ue^{-is A}} \le \tilde{v}(s) \ue^{\omega \abs{s}} \qquad (s\in \R)
\]
for some $\omega \in \R_{\ge 0}$ and some measurable function
$\tilde{v}: \R \to \R_{\ge 0}$. Let $v: \R \to [1, \infty)$ be strongly
admissible such that $\tilde{v}/v \in \Ell{1}(\R)$.  
Then $A$ has a bounded
$\Hr^2_v(\strip{\omega})$-calculus. Moreover, for 
each $\vphi \in \Ho[\cstrip{2\omega}]$ with $\vphi(0)=1$
one has
\[ f(A) = \int_\R F_s(A) \ue^{-\ui s
  A}\, \ud{s} \qquad (f\in \Hr^2_v(\strip{\omega})),
\]
where $F_s = (\tau_\bfz \vphi \cdot f)^\vee(s)$ for $s\in \R$.
 \end{thm}

\vanish{
Before we start the proof, let us comment on the
assumptions in the theorem. Fix $\omega'> \omega$.
Since admissible functions
grow at most polynomially, 
$v \ue^{\omega\abs{\bfs}} \lesssim
\ue^{\omega'\abs{\bfs}}$. Hence 
$\norm{\ue^{-\ui \bfs A} \ue^{-\omega'\abs{\bfs}}} \le v_U/v
  \in\Ell{1}(\R)$. By Datko's theorem \cite[Thm. 5.1.2]{ABHN}, 
it follows that $\ue^{-\ui \bfs A} \ue^{\omega' \abs{\bfs}}$ is
bounded. This shows that the {\em group type} of 
$\ue^{-\ui \bfs A}$ is $\theta( \ue^{-\ui \bfs A}) \le \omega$.
From this it follows by a result of Cowling, Doust, McIntosh and Yagi
\cite[Thm. 5.4.1]{Ha2006} that $A$ has bounded
$\Ha(\strip{\theta})$-calculus  for each $\theta > \omega$. 
}

\begin{proof}
We show first that $A$ is of strip type $\omega$. To this end, 
let $\Re \lambda \ge  \omega'> \omega$. 
Then, since $v$ is growing at most polynomially,  
\[   \norm{\ue^{\pm \ui sA} \ue^{-\lambda s}} \le 
\tilde{v}(s) \ue^{\omega s} \ue^{-\omega' s} 
= \frac{\tilde{v}(s)}{v(s)}   v(s) \ue^{-(\omega'- \omega)s}
\lesssim \frac{\tilde{v}(s)}{v(s)} 
\qquad (s\ge 0).
\]
Since $\tilde{v}/v\in \Ell{1}(\R)$, it follows that $\lambda \in
\resol(\pm \ui A)$ and  $\sup_{\re \lambda \ge \omega'}
\norm{R(\lambda, \pm \ui A)}< \infty$. Hence, 
$A$ is of strip type $\omega$.

\smallskip

Fix, without loss of generality,  $\theta  > \omega$ 
such that $A$ has a bounded $\Ha(\strip{\theta})$-calculus,
and fix  $\vphi \in \Ho(\strip{2\theta})$ with
$\vphi(0)=1$. Then by Theorem \ref{hrs.t.strongadm}
 we have 
the representation formula
\[
f(z)\, =\, \int_\R (\tau_z \vphi f)^\vee (s)\ue^{-\ui sz}\, \ud s
= \int_\R F_s(z) \ue^{-\ui sz}\, \ud s
\qquad (z\in \cstrip{\omega})
\] 
 for 
$f\in \Hr^2_v(\strip{\omega})$.
Note, however, that this formula holds even for $z\in \strip{\theta}$
whenever $f\in \Ho(\strip{\theta})$. We shall use this
fact in the following to establish the formula
\beq\label{bhc.eq.bhc}
 f(A) = \int_\R F_s(A) \ue^{-\ui s A}\, \ud{s}
\eeq
for such $f$ and  therefrom  the estimate
\[   \norm{f(A)} \lesssim \norm{f}_{\Hr^2_v(\strip{\omega})}
\qquad (f\in \Ho(\strip{\theta})).
\]
From Theorem \ref{bhc.t.bhc-char} 
it then follows that $A$ has bounded
$\Hr^2_v(\strip{\omega})$-calculus. Finally, 
we shall establish \eqref{bhc.eq.bhc} 
for all $f\in \Hr^2_v(\strip{\omega})$.

\smallskip
Fix $f\in \Ho(\strip{\theta})$ and define (as above) 
\[ F(z,s) := F_s(z) := (\tau_z \vphi \cdot f)^\vee(s)
= \frac{1}{2\upi} \int_\R \vphi(t - z) f(t)\ue^{\ui t s}\, \ud{t}
\qquad (z\in
\strip{\theta}, s\in \R).
\]
It is easy to see that $F$ is continuous, 
$F_s  \in \Ho(\strip{\theta})$
(uniformly in $s\in \R$) and 
\[   s\mapsto F_s(A) = \frac{1}{2\upi}
\int_\R \vphi(t- A)f(t) \ue^{\ui ts}\, \ud{t}
\]
is continuous in operator norm. 
Moreover, by Theorem \ref{hrs.t.indep}.b
\[
\norm{F_s}_{\Ha(\strip{\theta})} \lesssim_\delta  
 \frac{\ue^{-\delta |s|}}{v(s)} \norm{f}_{\Hr_v^2 (\strip{\delta})}
\qquad (s\in \R,\, \omega \le \delta \le \theta).
\]
Since the $\Ha(\strip{\theta})$-calculus is
bounded,
\[ \norm{F_s(A)} \lesssim   \frac{\ue^{-\omega |s|}}{v(s)} 
\norm{f}_{\Hr_v^2 (\strip{\omega})}
\qquad (s\in \R).
\]
This shows that the integral  
\[ \int_\R F_s(A) \ue^{-\ui s A} \, \ud s
\]
is absolutely convergent.  Fix $\lambda \in \C \ohne 
\cstrip{\theta}$ and $\delta \in (\omega , \theta)$. Then 
\begin{align*}
R(\lambda, A)^2 f(A)  =&\, \frac{1}{2\pi \ui}\int_{\partial \strip{\delta}} \frac{f(z)}{(\lambda - z)^2} R(z, A)\, \ud z\\
=&\, \frac{1}{2\pi \ui}\int_{\partial \strip{\delta}} \bigg(\int_\R (\tau_z \vphi \cdot f)^\vee (s) \ue^{-\ui sz}\, \ud s\bigg)\frac{1}{(\lambda - z)^2} R(z, A)\, \ud z\\
=&\, \frac{1}{2\pi \ui}\int_{\partial \strip{\delta}} \int_\R (\tau_z \vphi \cdot f)^\vee (s) \ue^{-\ui sz} \frac{1}{(\lambda - z)^2} R(z, A)\,\ud s\,  \ud z\\
=&\,  \int_\R  \bigg(\frac{1}{2\pi \ui} \int_{\partial \strip{\delta}} (\tau_z \vphi \cdot f)^\vee (s) \ue^{-\ui sz} \frac{1}{(\lambda - z)^2} R(z, A)\, \ud z\bigg)\, \ud s\\
  =&\,  \int_\R  F_s(A) \ue^{-\ui s A} R(\lambda, A)^2\, \ud s
= R(\lambda, A)^2 \Big(\int_\R F_s(A) \ue^{-\ui s A} \, \ud s\Big). 
\end{align*} As $R(\lambda, A)^2$ is injective, 
the identity (\ref{bhc.eq.bhc}) follows.
 Hence, one can estimate
\begin{align*}
\| f(A)\| \leq&\, \int_\R \| F_s(A)\|\, \| \ue^{-\ui s A} \|\, \ud s
\lesssim \, \int_\R \frac{\ue^{-\omega |s|}}{v(s)} \| f\|_{\Hr_v^2
                (\cstrip{\omega})} \, \| \ue^{-\ui sA }\|\, \ud s\\
\leq & \|\tilde{v}/ v \|_1\,\| f\|_{\Hr_v^2(\strip{\omega})}, 
\end{align*} 
As said before, this  by Theorem \ref{bhc.t.bhc-char}  
yields  that  $A$ has a bounded $\Hr^2_v(\strip{\omega})$-calculus. 

\smallskip
It remains to show that the formula \eqref{bhc.eq.bhc}
holds for all $\vphi \in \Ho[\cstrip{2\omega}]$ with $\vphi(0)=1$ and
all $f\in \Hr^2_v(\strip{\omega})$. Fix $\theta > \omega$ and  $\vphi\in \Ho(\strip{2\theta})$
with $\vphi(0)=1$. By what we have already shown, $A$
has a bounded $\Ha(\strip{\theta})$-calculus 
(simply because $\Ha(\strip{\theta}) \subseteq
\Hr^2_v(\strip{\omega})$ continuously) and hence  the
proof from above yields  \eqref{bhc.eq.bhc}
for all  $f\in  \Ho(\strip{\theta})$. 

Fix $f\in \Hr^2_v(\strip{\omega})$ and pick $K \ge 0$ and 
a sequence
$(f_n)_n$ in $\Ho(\strip{\theta})$ as in Lemma \ref{bhc.eq.bhc}.
Define 
\[ F_{s,n}(z) := (\tau_z \vphi\cdot f_n)^\vee (s) \qquad (z\in 
\strip{\theta},\, s\in \R).
\]
By Theorem \ref{fcc.t.conlem}, $f_n(A) \to f(A)$ strongly, and 
by the classical
convergence theorem, $F_{n,s}(A) \to F_s(A)$ strongly. 
Moreover,
\[  \sup_n \norm{F_{n,s}(A)} \lesssim
\frac{\ue^{-\omega\abs{s}}}{v(s)} K \norm{f}_{\Hr^2_v(\strip{\omega})}
\qquad (s\in \R). 
\]
Hence (dominated convergence!)
\[  f_n(A)x = \int_\R F_{s,n}(A) \ue^{-\ui s A}x\, \ud{s}
\to \int_\R F_s(A)\ue^{-\ui s A} x\, \ud{s} \qquad (x\in X).
\]
This concludes the proof.  
\end{proof}

Theorem \ref{bhc.t.meda} holds for any Banach space and uses the
bounded $\Ha$-calculus in a very coarse way. However, if $X$ has finite
cotype, bounded $\Ha$-calulus implies ``square function estimates''
and ``dual square function estimates'', i.e. an abstract
Littlewood--Paley theorem. Combining this with $\gamma$-boundedness
arguments (which guarantee boundedness with respect to the relevant
square function norms) then leads to a more refined statement about
a bounded Hörmander type functional calculus as follows.

\begin{thm}\label{bhc.t.main}
Let $X$ be a Banach space of finite cotype $q \in [2, \infty)$ and
type $p \in [1,2]$, and let $A$ be a densely defined operator on $X$ with a
bounded $\Ha(\strip{\theta})$-calculus for some $\theta > 0$. Suppose
that $\omega\ge 0$ and 
$\tilde{v} : \R \to (0, \infty)$ is measurable with 
\[ \norm{\ue^{-\ui s A}} \le \tilde{v}(s) \ue^{\omega \abs{s}}\qquad
(s\in \R).
\]
Let $r \in [1,2]$ with $\frac{1}{r } > \frac{1}{p} - \frac{1}{q}$ and  
let $v: \R \to [1, \infty)$ be strongly admissible such that $\tilde{v}/{v} \in
\Ell{r}(\R)$.  
Then, $A$ has a bounded $\Hr_v^2(\strip{\omega})$-calculus. If, in
addition, $X$ has Pisier's contraction property, then this calculus is $\gamma$-bounded.
\end{thm}

\begin{rem}
Observe that for a Banach space $X$ with type $p\in [1,2]$ and finite cotype $q <
\infty$ one has  $\frac{1}{p} - \frac{1}{q} < 1$. Hence Theorem
\ref{bhc.t.main} is stronger  than Theorem \ref{bhc.t.meda} for such
spaces $X$. 
\end{rem}

\bigskip

The proof requires several steps. In all what follows $X, p, q, A, \omega,
\tilde{v}, v $ are as in the theorem. As in the beginning of the proof
of Theorem \ref{bhc.t.meda}, 
we see that $A$ is of strip type $\omega$. Fix
$\theta> \omega$ such that $A$ has a bounded
$\Ha(\strip{\theta})$-calculus.\footnote{By a theorem of Cowling,
  Doust, McIntosh and Yagi \cite[Thm. 5.4.1]{Ha2006} 
or also by Theorem \ref{bhc.t.meda}, 
$\theta$ can be any number $> \omega$. However, this information
is inessential here.} 
In order to employ abstract Littlewood--Paley theory, we need
the following ``partition of unity''-result.

\begin{lem}\label{bhc.l.pou}
Let $\theta > 0$. Then there are $\psi, \psi_1, \psi_2\in
\Ho(\strip{\theta})$ with 
\[  \sum_{n\in \Z} \psi_1(z-n) \psi_2(z-n)\psi(z- n) = 1
\qquad (z\in \strip{\theta}).
\]
\end{lem}

\begin{proof}
Fix $0 < \alpha < \frac{\upi}{2 \theta}$ and $c > 0$ such that 
\[ \eta := \frac{c \ue^{\alpha \bfz}}{(1+ \ue^{\alpha \bfz})^2} = 
\frac{c}{(1 + \ue^{\alpha \bfz}) (1 + \ue^{-\alpha \bfz})} \in
\Ho(\strip{\omega})
\]
satisfies
\[   \int_\R \eta(s)\, \ud{s} = 1.
\]
Then $\eta$ maps the strip $\strip{\theta}$ into the sector
$\sector{\alpha \theta}$. Define
\[ \vphi(z) := \int_0^1 \eta(s- z)\, \ud{s} \qquad (z\in
\strip{\theta}).
\]
Then, still, $\vphi \in \Ho(\strip{\theta})$ and, by convexity,
$\vphi$ maps into $\cls{\sector{\alpha \theta}}$. Moreover,
\[ \sum_{n \in \Z} \vphi(t - n) = 
  \sum_{n \in \Z} \int_0^1 \eta(s -t + n) = \int_\R \eta(s)\,\ud{s}
= 1 \qquad (t\in \R).
\]
By uniqueness of holomorphic functions,
\[   \sum_{n\in \Z} \vphi(z- n) = 1 \quad \text{for all $z\in  \strip{\theta}$}.
\] 
Finally, since $\eta$ 
has no zeroes, $\re\eta(z) > 0$ and hence also $\re
\vphi(z) > 0$ for each $z\in \strip{\theta}$. This means that we can 
take 
\[   \psi(z) := \psi_1(z) := \psi_2(z) :=
\bigl(\vphi(z)\bigr)^\frac{1}{3} \qquad (z\in \strip{\theta}).\qedhere
\]
\end{proof}

\medskip

Back to the proof of Theorem \ref{bhc.t.main}, 
fix $\psi, \psi_1, \psi_2 \in \Ho(\strip{\theta})$ as in 
Lemma \ref{bhc.l.pou}. By Proposition \ref{apl.p.PL} we have 
the norm equivalences
\begin{align*}   
\norm{x} & \eqsim
           \bignorm{\bigl(\psi_1(A-n)x\bigr)_n }_{\gamma(\Z;X)} 
\\ & \eqsim
           \bignorm{\bigl(\psi(A-n)\psi_1(A-n)x \bigr)_n}_{\gamma(\Z;X)} 
           \qquad (x\in X).
\end{align*}
Then 
\[ 
\norm{f(A)x}  \eqsim
\bignorm{\bigl( (\tau_n\psi \cdot  f)(A)  \psi_1(A-n)x \bigr)_n}_{\gamma(\Z;X)} 
\qquad (x\in X,\ f\in \Ho(\strip{\theta})).
\]
Hence, in order 
to obtain the desired estimate
\[  \norm{f(A)x} \lesssim \norm{f}_{\Hr^2_v(\strip{\omega})} \norm{x}
\qquad (x\in X,\, f\in \Ho(\strip{\theta}))
\]
it is sufficient to establish 
\[   \norm{ \bigl( ( \tau_n\psi \cdot  f)(A) \bigr)_{n\in \Z}}_{\BL(\gamma(\Z;X))} \lesssim 
 \norm{f}_{\Hr^2_v(\strip{\omega})} \qquad (f\in \Ho(\strip{\theta})).
\]
These considerations lead to the following 
intermediate result, a generalization 
of \cite[Cor. 5.2]{KrWe2018}.

\begin{prop}\label{bhc.p.main-aux}
Let $\theta > \omega \ge 0$ and let
$A$ be a densely defined closed operator
with a bounded $\Ha(\strip{\theta})$ calculus
on a Banach space $X$.
Furthermore, let $0 \neq \psi \in \Ho(\strip{\theta})$ as above 
and $v: \R \to
[1, \infty)$ strongly admissible. 
\begin{aufzi}
\item If $X$ has finite cotype and
\[ \gammabound{ ( \tau_n\psi \cdot  f)(A) \suchthat n \in \Z } \lesssim
\norm{f}_{\Hr^2_v(\strip{\omega})} 
\qquad (f\in \Ho(\strip{\omega})),
\]
then $A$ has a bounded $\Hr^2_v(\strip{\omega})$-calculus.
\item If $X$ has Pisier's contraction property  and 
\[ \gammabound{  ( \tau_n\psi \cdot  f)(A) \suchthat n\in \N,\, f\in
  \Ho(\strip{\theta}),\, \norm{f}_{\Hr^2_v(\strip{\omega})}\le 1 } <\infty,
\] 
then $A$ has a $\gamma$-bounded $\Hr^2_v(\strip{\omega})$-calculus. 
\end{aufzi}
\end{prop}

\begin{proof}
a) By the considerations above and the identity
\[ \norm{ (T_n)_{n\in \Z}}_{\BL(\gamma(\Z; X))} 
\le \gammabound{ T_n \suchthat n \in \Z}
\]
for any sequence $(T_n)_n \in \BL(X)^\Z$ (Lemma \ref{alp.l.gbd})
we obtain an estimate 
\[  \norm{f(A)x}\lesssim \norm{f}_{\Hr^2_v(\strip{\omega})} \norm{x}
\qquad (x\in X,\, f\in \Ho(\strip{\theta})).
\]
Hence, the claim follows from Theorem \ref{alp.l.gbd}.

\prfnoi
b)\ Define
\[ C := \gammabound{  ( \tau_n\psi \cdot  f)(A) \suchthat n\in \N,\, f\in
  \Ho(\strip{\theta}),\, \norm{f}_{\Hr^2_v(\strip{\omega})}\le 1 } < \infty.
\]
Pick a finite sequence $(f_k)_k$ in $\Ho(\strip{\theta})$ with 
$\norm{f_k}_{\Hr^2_v(\strip{\omega})} \le 1$ for each $k$, and a
finite
sequence $(x_k)_k$ in $X$. 
Then 
\begin{align*}
\Exp \Bignorm{ \sum_k \gamma_k f_k(A)x_k}_X^2
& \eqsim 
\Exp \Bignorm{ \Bigl( \sum_k \gamma_k (\tau_n \psi \cdot f_k)(A)
\psi_1(A-n)x_k \Bigr)_n }_{\gamma(\Z;X)}^2
\\ & =
\Exp \Bignorm{ \sum_k \gamma_k \Bigl( (\tau_n \psi \cdot f_k)(A)
\psi_1(A-n)x_k \Bigr)_n }_{\gamma(\Z;X)}^2
\\ & =  \Bignorm{  \Bigl( (\tau_n \psi \cdot f_k)(A)
\psi_1(A-n)x_k \Bigr)_{n,k} }_{\gamma(\Z;\gamma(\Z;X))}^2
\\ & \eqsim 
\Bignorm{  \Bigl( (\tau_n \psi \cdot f_k)(A)
\psi_1(A-n)x_k \Bigr)_{n,k} }_{\gamma(\Z\times \Z;X)}^2
\\ & \le C^2  
\Bignorm{  \Bigl(\psi_1(A-n)x_k \Bigr)_{n,k} }_{\gamma(\Z\times
     \Z;X)}^2
\eqsim C^2 \Bignorm{ \sum_k \gamma_k x_k}_X^2.
\end{align*}
(Note that Pisier's contraction property  was employed in 
passing from $\gamma(\Z; \gamma(\Z;X))$ to $\gamma(\Z\times \Z;X)$
and back.) We conclude that 
\[    \gammabound{ f(A) \suchthat f\in\Ho(\strip{\theta}),\,
  \norm{f}_{\Hr^2_v(\strip{\omega})} \le 1} \le C < \infty
\]
and hence, by Theorem \ref{alp.l.gbd}, that $A$ has $\gamma$-bounded
$\Hr^2_v(\strip{\omega})$-calculus. 
\end{proof}

We now make the final step in the proof of Theorem \ref{bhc.t.main}.

\begin{proof}[Proof of Theorem \ref{bhc.t.main}]
Let, as before, $X, p, q, A, \omega,
\tilde{v}, v $ as in the formulation of the theorem, and $\psi,
\psi_1, \psi_2 \in \Ho(\strip{\theta})$ as in Lemma \ref{bhc.l.pou}. By
Proposition \ref{bhc.p.main-aux} 
it suffices to show
\[   \gammabound{  ( \tau_n\psi \cdot  f)(A) \suchthat n\in \N,\, f\in
  \Ho(\strip{\theta}),\, \norm{f}_{\Hr^2_v(\strip{\omega})}\le 1 } <\infty.
\]
To that end, fix $f\in \Ho (\strip{\theta})$, and note that
\begin{align*}
(\tau_n \psi \cdot f) (A)x  =&\, \int_\R (\tau_n \psi\cdot f)^\vee
                                (s) \ue^{-\ui sA} x\, \ud s\\
=&\, \int_\R \left( v(s) \ue^{\omega |s|} (\tau_n \psi\cdot  f)^\vee
   (s) \right) \cdot \left(\frac{\ue^{-\omega |s|}}{v (s)} \ue^{-\ui sA} x\right)\ud s
\end{align*} for each $n\in \Z$ and $x\in X$. By hypothesis,
\[
\frac{\ue^{-\omega |\bfs |}}{v} \|\ue^{-\ui \bfs A} x\|\, \lesssim \,
\frac{\tilde{v}}{ v} \norm{x}\quad \text{and}\quad \frac{\tilde{v}}{ v} \in \, \Ell{r} (\R ).
\] 
Hence, by Theorem \ref{alp.t.intmeans} 
it suffices to establish an estimate of the form
\[ 
\sup_{n\in \Z} \left(\int_\R  \big| v(s) \ue^{\omega |s|} (\tau_n
  \psi\cdot f)^\vee (s) \big|^{r'}\, \ud s \right)^{\frac{1}{r'}}\, \
\lesssim\, \| f\|_{\Hr_v^2(\strip{\omega})} \qquad(f\in \Ho(\strip{\theta})),
\]
where  $\tfrac{1}{r'} = 1 - \tfrac{1}{r}$ is the
dual exponent. For this, note that on one hand
\[
\sup_{n\in \Z} \left(\int_\R  \big| v(s) \ue^{\omega |s|} (\tau_n
  \psi\cdot f)^\vee (s) \big|^{2}\, \ud s \right)^{\frac{1}{2}}\, =\,
\sup_{n\in \Z } \| \tau_n \psi \cdot f\|_{\We_v^2(\strip{\omega})}\,
\lesssim\, \| f\|_{\Hr_v^2 (\strip{\omega})}. 
\] 
On the other hand,  by Theorem \ref{hrs.t.indep}.b),
\[
\sup_{n\in \Z } \sup_{s\in \R} \big| v(s) \ue^{\omega |s|} (\tau_n \psi\cdot f)^\vee (s) \big|\, \lesssim\, \| f\|_{\Hr_v^2 (\strip{\omega})}.
\] 
As $r'\in [2, \infty]$, the
desired estimate follows by interpolation.
\end{proof}

One can do with weaker geometric assumptions
if one assumes a stronger boundedness condition on the group. Recall
from Appendix \ref{s.alp}
the notion of a semi-$\gamma$-bounded family of operators.

\begin{thm}\label{bhc.t.main-semi}
Let $X$ be a Banach space 
with Pisier's contraction property
and let $A$ be a densely defined operator on $X$ with a
bounded $\Ha(\strip{\theta})$-calculus for some $\theta > 0$. Suppose
that $\omega\ge 0$ and 
$\tilde{v} : \R \to (0, \infty)$ is measurable such that the operator family
\[  \{     \tfrac{\ue^{-\omega \abs{s}}}{\tilde{v}(s)} \ue^{-\ui sA} \suchthat s\in \R\}
\]
is semi-$\gamma$-bounded. Let 
$v: \R \to [1, \infty)$ be strongly admissible such that
$\frac{\tilde{v}}{v} \in \Ell{2}(\R)$. 
Then, $A$ has a $\gamma$-bounded $\Hr_v^2(\strip{\omega})$-calculus.
\end{thm}

\begin{proof}
We proceed as in the proof of Theorem \ref{bhc.t.main} and pick
$\theta > \omega$
such that $A$ has bounded $\Ha(\strip{\theta}$-calculus, as well as
functions
$\psi, \psi_1, \psi_2 \in \Ho(\strip{\theta})$ as in Lemma
\ref{bhc.l.pou}. As before
it suffices to show
\[   \gammabound{  ( \tau_n\psi \cdot  f)(A) \suchthat n\in \N,\, f\in
  \Ho(\strip{\theta}),\, \norm{f}_{\Hr^2_v(\strip{\omega})}\le 1 } <\infty.
\]
To that end, we  again fix $f\in \Ho (\strip{\theta})$ and consider
the formula
\begin{align*}
(\tau_n \psi \cdot f) (A)x  =&\, \int_\R (\tau_n \psi\cdot f)^\vee
                                (s) \ue^{-\ui sA} x\, \ud s\\
=&\, \int_\R \left( v(s) \ue^{\omega |s|} (\tau_n \psi\cdot  f)^\vee
   (s) \right) \cdot \frac{\tilde{v}(s)}{v(s)} \left(\frac{\ue^{-\omega |s|}}{\tilde{v} (s)} \ue^{-\ui sA} x\right)\ud s
\end{align*} for each $n\in \Z$ and $x\in X$. Now an application of  Theorem
\ref{alp.t.semi} proves the claim.
\end{proof}

\begin{rems}\label{bhc.r.main}
1)\ Theorem \ref{bhc.t.main}   is particularly relevant when $X= \Ell{p}(\Omega)$
for some $1 < p < \infty$ and some measure space $\Omega$. Namely, in
this case $\type(X) = \min\{2, p\}$ and $\cotype(X) =
\max\{p,2\}$. This yields
\[  \frac{1}{\type(X)} - \frac{1}{\cotype{X}} = \bigabs{\frac{1}{2} -
  \frac{1}{p}} < \frac{1}{2}.
\]
Hence, {\em each} $r\in [1,2]$ is a possible choice, and one obtains the
largest possible variety of bounded Hörmander calculi. 
However, as integrability is
usually realized via growth conditions, one should think of
$\frac{\tilde{v}}{v}$ being bounded. In this situation one has the
implication $\frac{\tilde{v}}{v} \in
\Ell{r},\, r\in [1,2] \,\dann\, \frac{\tilde{v}}{v} \in \Ell{2}$ and, hence, 
$r=2$ is the optimal integrability exponent.

\prfnoi
2)\ Theorems \ref{bhc.t.main} and \ref{bhc.t.main-semi} extend \cite[Thm. 10.2]{KrWe2018} by Kriegler
and Weis (for the
case $r=2$) from strip type $\omega = 0$ to general $\omega \ge
0$, and from classical Hörmander spaces (polynomial weights) to
the finer scale of admissible functions. 

To illustrate the latter, consider the case 
$\tilde{v}  = (1+|\bfs|)^\alpha$ and $X$ being an $\Ell{p}$-space as
above. Then, restricting to classical
Sobolev spaces (=polynomial weights) just allows one
to infer a $\gamma$-bounded $\Hr^2_{v_\beta}(\strip{\omega})$-calculus
for $\beta > \alpha + \frac{1}{2}$. For $\omega=0$, this
is Kriegler and Weis' result.  However, with our finer
regularity scale we are able to infer
a $\gamma$-bounded $\Hr^2_v(\strip{\omega})$ calculus even for the
admissible functions 
\[
v\, =\, (1+|\bfs|)^{\alpha + \frac{1}{2}} \big( \ln (\ue + |\bfs|)\big)^{\beta}\qquad (\beta > 1)
\] 
or
\[
v\, =\, (1+|\bfs|)^{\alpha + \frac{1}{2}} \ln (\ue+|\bfs|) \big( \ln \big(\ue \ln (\ue+|\bfs|)\big) \big)^\beta \qquad (\beta >1).
\] 

\prfnoi
3)\ However, \cite[Thm. 10.2]{KrWe2018} is not completely
covered by Theorem \ref{bhc.t.main}, as  we exclusively work with $\Ell{2}$-Sobolev spaces, where
Kriegler and Weis also consider 
Hörmander spaces based on Sobolev spaces $\We^{\beta, p} (\R)$ for
$p\neq 2$. 
\end{rems}

\section{Sectorial Operators}\label{s.sec}

In this section we describe how the obtained
(function- and operator-theoretic) results 
transfer to the sectorial situation. This 
is the ``original'' set-up and certainly
the one of most interest from the point of view of
applications.

We start with the function theory. In all what follows,
$\omega \in [0, \upi)$ and $v: \R \to [1, \infty)$ is admissible. 
Recall that
\[  \sector{\omega} := \begin{cases}
\{ z\in \C \ohne \{0\} \suchthat 
\abs{\arg z} < \omega\} & \text{if $\omega > 0$}\\
\R_{> 0} & \text{if $\omega = 0$}
\end{cases}
\]
is the open sector of angle $2\omega$ symmetric about the positive
real axis. In addition, we let 
\[
\csect{\omega}\coloneqq \cl{\sector{\omega}}\ohne\{0\} \qquad (\omega\in [0, \upi]).
\]
The function $\ue^\bfz$ maps $\strip{\omega}$ biholomorphically 
onto $\sector{\omega}$, with inverse $\log \bfz$. With the help
of this ``change of coordinates'' we shall transfer concepts
from strips to sectors.

For $\omega > 0$, the  \emdf{Hardy space} of order $p\in [1,\infty]$
on $\sector{\omega}$ is \label{f.Hardy-sector}
\[    \Har{p}(\sector{\omega}) := \{f \in \Hol(\sector{\omega})
\suchthat 
f(\ue^\bfz) \in \Har{p}(\strip{\omega})\}
\]
endowed with its natural norm. Then the space of 
{\emdf elementary functions} on $\sector{\omega}$ is
\label{f.elem-sector}
\[ \calE(\sector{\omega}) = \Har{1}(\sector{\omega}) \cap
\Har{\infty}(\sector{\omega} = \{ f \suchthat f(\ue^\bfz) \in
\calE(\strip{\omega}\}.
\]
Similarly, \label{f.holSchwartz-sector}
\[ \Ho(\sector{\omega}) := \{ f\suchthat  f(\ue^\bfz) \in
\Ho(\strip{\omega})\}.
\]
Observe that this space is larger than the space
$\Ho(\strip{\omega})$ considered in  \cite[Sec. 2.2]{Ha2006}. 
However, this should not lead to confusion.

\medskip

The \emdf{generalized Sobolev space} \label{f.Sob-sector} 
on $\sector{\omega}$ associated with $v$ is 
\[
\We_v^2 (\sect{\omega}) \coloneqq \big\{ f : \sector{\omega}\to \C
\suchthat  f(\ue^{\bfz}) \in \We^2_v(\strip{\omega})\}.
\]
It is  equipped with the norm 
\[
\| f\|_{\We_v^2 (\sect{\omega})} \coloneqq \|f(\ue^\bfs) \|_{\We_v^2(\strip{\omega})}
\] 
Similarly, the  \emdf{generalized Hörmander space} 
 with respect to $v$ is \label{f.Hrse}
\[
\Hr_v^2 (\sect{\omega}) \coloneqq \big\{ f: \sect{\omega}  \to \C \, \big|\,  f(\ue^\bfz) \in \Hr_v^2 (\strip{\omega})  \big\}.
\] 
It is equipped with the norm 
\[
\| f\|_{\Hr_v^2 (\sect{\omega})} \coloneqq \|f(\ue^\bfz) \|_{\Hr_v^2(\strip{\omega})}.
\]
We abbreviate \label{f.Sob-sector-class}\label{f.Hrse-class}
\[
\We^{\alpha, 2} (\sect{\omega}) \coloneqq \We_{v_\alpha}^2 (\sect{\omega}), 
\qquad \text{and}\qquad
\Hr^{\alpha, 2} (\sect{\omega})\coloneqq \Hr_{v_\alpha}^2 (\sect{\omega}),
\] 
where $\alpha \ge  0$ and $v_\alpha = (1+ |\bfs|)^{\alpha}$. Recall
that $v_\alpha$ is strongly admissible if and only if $\alpha > \tfrac{1}{2}$.

\medskip

Key properties of the spaces $\We_v^2 (\sect{\omega})$ and 
$\Hr_v^2 (\sect{\omega})$ can easily be deduced from 
their strip counterparts. Here we just list a few of them
for the case that $v$ is strongly admissible. 

\begin{prop}\label{sec.p.sobhr}
Let $0\leq \omega <\theta < \pi$, and  $v: \R \to [1, \infty )$ be 
strongly admissible. Then the following statements hold:
\begin{aufzi}
\item $\calE(\sect{\theta}) \subseteq \We_v^2 (\sect{\omega})$ and
  the space
\[ \{ f \in \Hol(\sect{\theta}) \suchthat \exists\, a > 0 : \abs{f}
\lesssim \ue^{- a (\ln \abs{\bfz})^2}\}
\]
is dense in $\We_v^2 (\sect{\omega})$.
\item There are canonical (continuous) embeddings 
\[
\We_v^2(\sect{\omega})\hookrightarrow \Co(\csect{\omega}), \quad \Hr_v^2 (\sect{\omega}) \hookrightarrow \Cb (\csect{\omega}), \quad\text{and}\quad \Ha (\sect{\theta}) \hookrightarrow \Hr_v^2 (\sect{\omega}).
\] 
\item The spaces $\We_v^2(\sector{\omega})$ and  $\Hr_v^2(\sector{\omega})$
are Banach algebras with respect to pointwise multiplication.

\item In case  $\omega > 0$ one has  $f\in \Hr^2_v(\sect{\omega})$ if and
  only if $f\in \Ha ( \sect{\omega})\cap \Cb (\csect{\omega})$ and
  $f(\ue^{\pm \ui \omega}\bfs) \in \Hr^2_v(\R_{> 0})$.
\end{aufzi}
\end{prop}

\begin{proof}
a) follows from Lemma \ref{sob.l.den}.
b) and c) follow from Proposition \ref{sob.p.We}.d and Theorem
\ref{hrs.t.strongadm}, part a) and c). d) follows from 
Theorem
\ref{hrs.t.strongadm}.d.
\end{proof}

\begin{rem}[Connection to Classical Hörmander Spaces]\label{sec.r.cHr}
Let $\alpha > \tfrac{1}{2}$ and $\omega > 0$.  
The classical Hörmander condition of order $\alpha$ for 
$g\in \Ell{\infty}(\R_+)$ reads
\begin{equation}\label{sec.eq.cHr}
\sup_{t>0} \| \eta\cdot g(t\bfs)\|_{\We^{\alpha, 2} (\R)}\, <\, \infty,
\end{equation} 
where $0\neq \eta\in \Test(\R_{> 0})$ is arbitrary. Since test functions on $\R_{> 0}$
and test functions on $\R$ correspond to each other via the
$\exp$-$\log$-correpondence, \eqref{sec.eq.cHr} is equivalent to
$g(\ue^{\bfs}) \in \Hr^2_{v_\alpha}(\R)$. (Recall Theorem
\ref{hrr.t.indep} for the independence from the defining function $\eta$.)  
It follows that the Hörmander space $\Hr^{\alpha, 2}(\sect{0})$ as
defined above coincides precisely with the classical Hörmander space
$\Hr^{\alpha, 2}(\R_+)$.

Consequently, by Proposition \ref{sec.p.sobhr}.d, the Hörmander
space $\Hr^{\alpha, 2}(\sect{\omega})$ as defined above coincides
precisely with the space  $\Ha (\sect{\omega}; \alpha)$ considered in
\cite{GCMMST2001}, \cite{MMS2004}, \cite{Sasso2005}, and
\cite{CarDra2017}, cf. also \eqref{int.eq.hor-frac}.
\end{rem}

Let us now turn to the operator theory. An operator $A$ on a Banach
space
$X$ is \emdf{sectorial} of angle $\omega \in [0, \upi)$ if  $\spec(A)
\subseteq \cls{\sect{\omega}}$ and $\sup \{ \norm{\lambda R(\lambda,A)}
  \suchthat \lambda \in C\ohne \cls{\sect{\omega'}}\} < \infty$ for
  each $\omega < \omega'< \upi$. 

The theory of sectorial operators and their
  holomorphic functional calculus is presented at many places in the
  literature, e.g. in \cite{Ha2006} or \cite[Chap. 10]{HvNVW2017}, and
  we assume the reader to be  familiar with the basic facts.
  In particular, we recall that only for {\em injective} sectorial
  operators there is a satisfying notion of a  bounded
  $\Ha$-calculus on a sector. 
Given a  non-injective sectorial operator $A$ on a Banach space $X$ 
one either has to pass
to $X/\ker(A)$ or to $\cls{\ran}(A)$,
cf. \cite[Prop. 2.2.1.h]{Ha2006}.

If $A$ is injective and sectorial of angle $\omega$, 
the operator $\log(A)$   (defined via the functional
calculus) is of  strip-type  $\omega$,  and it is
densely defined
if and only if $A$ is densely defined and has dense range. 

The holomorphic functional calculi
for $A$ and $\log(A)$ are connected via the \emdf{composition rule}
\beq\label{sec.eq.cr}    f(A) = [f(\ue^\bfz)](\log(A)),\qquad  g(\log(A)) = [g(\log \bfz)](A),
\eeq
see \cite[Cor. 4.2.5]{Ha2006}.  Analogous to the strip
case we can define the notion of 
a bounded $\We^2_v(\sect{\omega})$-calculus
and a bounded ($\gamma$-bounded) $\Hr^2_v(\sect{\omega})$-calculus.
The following is straightforward to check.

\begin{prop}[Composition rule]\label{fcc.thm.CompositionRule} 
Let $A$
 be an injective sectorial operator on a Banach space $X$, let
$\omega \in [0, \upi)$ and $v: \R \to [1, \infty)$ strongly
admissible. Then $\log A$ has a bounded $\We_v^2
(\strip{\omega})$-calculus if and only if $A$ has a bounded $\We_v^2
(\sect{\omega})$-calculus, and the calculi are related
via the composition rule \eqref{sec.eq.cr}.
\end{prop}

As in the strip case, a bounded Sobolev calculus implies
an (at least) unbounded Hörmander calculus by algebraic
regularization. The connection via the composition rule
remains valid. Therefore, we can simply transfer 
theorems about bounded Hörmander calculus on strips
to theorems about bounded Hörmander calculus on sectors.
The first is the sectorial analogue of Theorem \ref{bhc.t.meda}.

\begin{thm}\label{sec.t.meda}
Let $A$ be a densely defined operator on a Banach space $X$, 
with dense range and  a bounded $\Ha(\sect{\theta})$-calculus for some $\theta > 0$. Suppose 
\[  \norm{ A^{-\ui s}} \le \tilde{v}(s) \ue^{\omega \abs{s}} \qquad (s\in \R)
\]
for some $\omega \in \R_{\ge 0}$ and some measurable function
$\tilde{v}: \R \to \R_{\ge 0}$. Let $v: \R \to [1, \infty)$ be strongly
admissible such that $\tilde{v}/v \in \Ell{1}(\R)$.  
Then $A$ has a bounded
$\Hr^2_v(\sect{\omega})$-calculus. 
 \end{thm}

\begin{proof}
The operator $\log(A)$ satisfies the hypotheses of Theorem
\ref{bhc.t.meda}. Hence, $\log(A)$ has a bounded 
$\Hr^2_v(\strip{\omega})$-calculus. Transferring this
back to the sector with the composition rule yieds that 
$A$ has a bounded $\Hr^2_v(\sect{\omega})$-calculus.
\end{proof}

And here is the sectorial analogue of Theorem \ref{bhc.t.main}.

\begin{thm}\label{sec.t.main}
Let $X$ be a Banach space of finite cotype $q \in [2, \infty)$ and
type $p \in [1,2]$, and let $A$ be a densely defined operator on $X$
with dense range and a
bounded $\Ha(\sect{\theta})$-calculus for some $\theta > 0$. Suppose
that $\omega \in [0, \upi)$ and 
$\tilde{v} : \R \to (0, \infty)$ is measurable with 
\[ \norm{A^{-\ui s}} \le \tilde{v}(s) \ue^{\omega \abs{s}}\qquad
(s\in \R).
\]
Let $r \in [1,2]$ with $\frac{1}{r } > \frac{1}{p} - \frac{1}{q}$ and  
let $v: \R \to [1, \infty)$ be strongly admissible such that $\tilde{v}/{v} \in
\Ell{r}(\R)$.  
Then, $A$ has a bounded $\Hr_v^2(\sect{\omega})$-calculus. 
If, in addition, $X$ has Pisier's contraction property, then this calculus is $\gamma$-bounded.
\end{thm}

\begin{proof}
Analogous to the proof of Theorem \ref{sec.t.meda}, now by virtue
of Theorem \ref{bhc.t.main}. 
\end{proof}

Finally, here is the analogue of 
Theorem \ref{bhc.t.main-semi}.

\begin{thm}\label{sec.t.main-semi}
Let $X$ be a Banach space with Pisier's contraction property, and let $A$ be a densely defined operator on $X$
with dense range and a
bounded $\Ha(\sect{\theta})$-calculus for some $\theta > 0$. Suppose
that $\omega \in [0, \upi)$ and 
$\tilde{v} : \R \to (0, \infty)$ is measurable such that the operator family
\[ \bigl\{ \tfrac{\ue^{-\omega \abs{s}}}{\tilde{v}(s)} A^{-\ui s}
\bigr\} 
\]
is semi-$\gamma$-bounded. Let
$v: \R \to [1, \infty)$ be strongly admissible such that $\tilde{v}/{v} \in
\Ell{2}(\R)$.  
Then, $A$ has a $\gamma$-bounded $\Hr_v^2(\sect{\omega})$-calculus. 
\end{thm}

Theorems \ref{sec.t.main} and 
\ref{bhc.t.main-semi} are (for $r=2$) 
a generalization of \cite[Thm. 6.1(2)]{KrWe2018}
to sectoriality angles $\omega >0$ and to more general smoothness
conditions, cf. also the 
Remarks \ref{bhc.r.main} above.

\begin{rem}[Non-injective Operators]\label{sec.r.injpart}
Let  $A$ be a densely defined sectorial operator of angle $\omega$ on
a Banach space $X$, and let $Y := \cls{\ran}(A)$. Then $Y$ is
invariant under the resolvent of $A$ and hence 
the  part $A_Y$ of $A$ in $Y$ is sectorial of angle
$\omega$, with $R(\lambda, A_Y) = R(\lambda,A)\res{Y}$ for $\lambda
\in\resol(A)$, 
cf. \cite[Sec. 2.6.2]{Ha2006}.
Moreover,  $A_Y$ is densely defined and has dense range. (This follows
from \cite[Prop. 2.1.1.c]{Ha2006}.)  We shall call $A_Y$ the
\emdf{injective part} of $A$. It is easy to see that $f(A_Y) =
f(A)\res{Y}$ for $f\in \calE[\sect{\omega}]$
\cite[Prop. 2.3.6]{Ha2006}. 

By abuse of notation, we shall say that $A$ \emdf{has a bounded
  $\calF$-calculus} (where $\calF$ may be $\Ha(\sect{\omega})$ or
$\We^2_v(\sect{\omega})$ or $\Hr^2_v(\sect{\omega})$), whenever its
injective part $A_Y$
does so.  With this terminology, Theorems \ref{sec.t.meda}---\ref{sec.t.main-semi}
remain valid even without the assumption that $A$ has dense range. 

Note that if $X$ is
reflexive then it decomposes as a direct sum
\[ X = \cls{\ran}(A) \oplus \ker(A).
\]
Clearly $f(A)\res{\ker(A)} = 0$  for each elementary function
$f$. Hence, one may formulate characterizations of a bounded
$\calF$-calculus for $A$ (i.e., in effect, for $A_Y$) similar to
Theorems \ref{bsc.p.bsc-char} and \ref{bhc.t.bhc-char}. We shall not
elaborate on this.

An alternative way to deal with a non-injective oprator $A$ is to perturb
it into $A+ \veps$ and then require functional calculus bounds that are
independent of $\veps > 0$. This has, e.g., been done in
\cite{GCMMST2001}. Again, this leads to equivalent concepts, but also
here we do not go further into detail.
\end{rem}

\section{Applications}\label{s.app}

\subsection*{Groups on Hilbert Spaces}

Let $H$ be a Hilbert space and $-\ui A$ the generator
of a $C_0$-group $(U_s)_{s\in \R}$ on $H$ with 
\[   \norm{U_s} \le \tilde{v}(s) \ue^{\omega \abs{s}}
\]
where $\omega \ge 0$ and $\tilde{v}$ is some
polynomially bounded function. By the
Boyadzhiev--de Laubenfels theorem \cite[Thm. 7.2.1]{Ha2006} 
$A$ has a bounded $\Ha(\strip{\theta})$-calculus for
each $\theta > \omega$. 

As $H$ has type and cotype $2$, Theorem \ref{bhc.t.main} implies
that $A$ actually has a bounded $\Hr^2_v(\strip{\omega})$-calculus
whenever $v$ is strongly admissible and
$\frac{\tilde{v}}{v} \in \Ell{r}(\R)$ for some $r\in [1,2]$.
As discussed in Remarks \ref{bhc.r.main},  integrability is usually 
realized through growth conditions, and 
one should think of $\frac{\tilde{v}}{v}$ being bounded. In this
respect, $r=2$ represents the optimal exponent.

For example, if $\tilde{v} \equiv c$ is constant, then $A$ has
bounded $\Hr^2_v(\strip{\omega})$-calculus for each strongly
admissible function $v$. In particular, $A$ has a (classical) $\Hr^{\alpha,2}(\strip{\omega})$-calculus
for each $\alpha > \frac{1}{2}$.

\subsection*{Sectorial Operators  on $\Ell{p}$-Spaces}

Let $(\Omega, \mu)$ be a measure space, let $1 < p < \infty$  and
$A$ a sectorial operator on $X = \Ell{p}(\Omega)$ with a bounded
$\Ha(\sect{\theta})$-calculus for some angle $\theta > 0$. By abuse of
notation we write $A$ also for its injective part (see Remark
\ref{sec.r.injpart}). Suppose one has an estimate
\[ 
\norm{ A^{-\ui s}} \le \tilde{v}(s) \ue^{\omega \abs{s}} \qquad
(s\in \R)
\]
where $\tilde{v}$ is some positive measurable function and $\omega
 \ge 0$. As discussed in Remarks \ref{bhc.r.main}, we have
\[   \frac{1}{\type(X)} - \frac{1}{\cotype(X)} = \bigabs{ \frac{1}{2} -
  \frac{1}{p} } <  \frac{1}{2}.
\]
Hence,  $r=2$ is a possible choice in Theorem \ref{sec.t.main} (and
the optimal one if, in addition, $\frac{\tilde{v}}{v}$ is bounded). It
follows that  $A$ has a $\gamma$-bounded $\Hr^2_v(\sect{\omega})$-calculus
whenever $\frac{\tilde{v}}{v} \in \Ell{2}(\R)$. In particular,
if $\tilde{v} \lesssim (1 + \abs{\bfs})^\alpha$ for some $\alpha > 0$
is constant, then $A$ has
$\gamma$-bounded $\Hr^{\beta,2 }(\sector{\omega})$-calculus for each $\beta > \alpha + 
\frac{1}{2}$.

\subsection*{Symmetric Contraction Semigroups}

Let $(\Omega, \mu)$ denote a $\sigma$-finite measure space. 
We shall abbreviate $\Ell{1}\coloneqq \Ell{1}(\Omega, \mu)$ and
$\Ell{\infty}\coloneqq \Ell{\infty}(\Omega, \mu)$.
A family of operators
\[
T_t : \Ell{1}\cap \Ell{\infty} \to \Ell{1}+ \Ell{\infty} \qquad (t\geq 0)
\] 
is called a \emdf{symmetric contraction semigroup} if it has the following properties:
\begin{aufziii}
\item $T_t T_s = T_{t+s}$ for all $s, t\geq 0$;
\item $\| T_t f\|_1 \leq \| f\|_1$ and $\| T_t f\|_\infty \leq \| f\|_\infty$ for all $f\in \Ell{1}\cap \Ell{\infty}$ and all $t\geq 0$;
\item $\|T_t f - f\|_p \xrightarrow{t\searrow 0} 0$ for all $f\in \Ell{1}\cap \Ell{\infty}$ and $p\in [1, \infty )$;
\item $\dps \int_\Omega T_t f \cdot \cl{g}\, \ud \mu = \int_\Omega  f \cdot \cl{T_t g}\, \ud \mu$ for each $f, g\in \Ell{1}\cap \Ell{\infty}$ and $t\geq 0$.
\end{aufziii}
It is an easy consequence of the Riesz-Thorin interpolation theorem
that each symmetric contraction 
semigroup extends to a contractive $C_0$-semigroup 
on $\Ell{p}(\Omega , \mu)$ for each $1\leq p < \infty$. We denote the respective generator by $-A_p$. 
Cowling \cite{Co1983} has shown that $A_p$ has a bounded
$\Ha(\sect{\theta})$-calculus for each angle $\theta > \upi
\abs{\frac{1}{p} - \frac{1}{2}}$.  In \cite[Thm. 2.1, Cor. 2.2]{CoMe1993}, Cowling and
Meda combined this result with ideas from \cite{Me1990} to obtain
sharper functional calculus bounds from additional hypotheses about the
growth behaviour of the group of imaginary powers $(A_p^{-\ui
  s})_{s\in \R}$. Moreover, a bounded Hörmander-type calculus on
$\sect{0}$ was inferred when the imaginary powers grow at most
polynomially \cite[Cor. 2.3]{CoMe1993}.
As described in the Introducion, in  \cite{GCMMST2001} García-Cuerva et al. extended the latter result to the
sectorial case, see Theorem \ref{int.t.GCetal}.

For a long time it remained an open problem to
determine the  optimal angle for the calculus. 
The angle in Cowling's result from above  was obtained by complex
interpolation and turned out to be  not optimal. 
Finally, in \cite{CarDra2017}, Carbonaro and Dragičević showed that 
this angle is
\[
\omega_p \coloneqq \arcsin \left|1 - \frac{2}{p}\right|.
\] 
In fact, they showed more \cite[Proposition 11]{CarDra2017}:

\begin{thm}[Carbonaro--Dragičević]\label{app.t.CD}
Let $-A$ be the generator of a symmetric contraction semigroup
over some measure space $(\Omega,\mu)$. Then for each $1 < p < \infty$
there is a constant $c_p \geq 0$ such that
\[
\| A_p^{-\ui s} f\|_p \, \leq\, c_p (1+ |s|)^{\frac{1}{2}} \ue^{\omega_p
  |s|} \| f\|_p \qquad (s\in \R,\, f\in\cl{\ran}(A_p)).
\] 
\end{thm}

As said in the introduction, one can pull out of this a bounded
$\Hr^{\beta, 2}(\sect{\omega_p})$-calculus for $A_p$ for all $\beta >
\frac{3}{2}$ with the help of the result by
García-Cuerva et al., see \cite[Theorem 1]{CarDra2017}.
However, with our refined methods, we can say more:

\begin{thm}\label{app.t.CD-improved}
Let $1<p<\infty$  and let $-A_p$ be the $\Ell{p}$-generator of a symmetric
contraction semigroup. Then  (the injective part of) $A_p$  
has a $\gamma$-bounded
$\Hr^{2}_v(\sect{\omega_p})$-calculus for each admissible function
$v:\R \to [1, \infty)$ such that 
\[
\frac{(1+ |\bfs|)^{\frac{1}{2}}}{v}\, \in \, \Ell{2}(\R).
\]
In particular, $A_p$ has a $\gamma$-bounded $\Hr^{\beta, 2}(\sect{\omega_p})$-calculus for each $\beta > 1$.  
\end{thm}

\subsection*{The Ornstein--Uhlenbeck Semigroup}

The Carbonaro--Dragičević result is optimal 
at least with respect to the angle $\omega_p$. 
In fact, it had been shown before that a limiting
case is provided by the so-called
{\em Ornstein--Uhlenbeck semigroup} on $(\R^d, \gamma_d)$, where
$\gamma_d$ is the standard Gaussian measure on $\R^d$, i.e.,
\[
\gamma_d (\ud x) \, =\, \frac{1}{(2\pi)^{\nicefrac{d}{2}}}\ue^{-\frac{|x|^2}{2}}\, \ud x.
\] 
The $\Ell{p}$-generator is given by $-\calL_p$, where
\[
\mathcal{L}_p \coloneqq - \frac{1}{2}\Delta + \mathbf{x}\cdot \nabla,
\]
and  where $\Delta$ is the Laplacian on $\R^d$, 
$\bfx$ is the mapping $x\mapsto x$, $\nabla$ is the gradient operator on
$\R^d$. 

It has been noted in \cite[Thm. 2]{GCMMST2001} that $\calL_p$ does not
have a bounded $\Ha(\sect{\omega})$-calculus for any $\omega <
\omega_p$. Actually, $\omega_p$ is the precise sectoriality angle
of $\calL_p$, see \cite{ChFaMePa2005}. 
Hence, the angle $\omega_p$ in Theorem
\ref{app.t.CD-improved} is optimal.

Nevetheless, $\calL_p$ itself allows for a stronger statement. 
In \cite[Thm. 4.3]{MMS2004}, Mauceri, Meda and Sjögren 
proved that for $1 < p < \infty$ there is $c_p \ge 0$ such that 
\[
\| \mathcal{L}_p^{-\ui s} f\|_p \, \leq\, c\ue^{\omega_p |s|} \|f\|_p 
\qquad (f\in  \cls{\ran} \mathcal{L}_p,\, s\in \R). 
\] 
where by abuse of notation we again write $\mathcal{L}_p$ for the injective
part of $\mathcal{L}_p$.  Hence, we may apply our findings 
from above and obtain:

\begin{thm}\label{app.t.OU}
Let $1<p<\infty$. Then the Ornstein--Uhlenbeck operator $\calL_p$ 
on $\Ell{p}(\R^d, \gamma_p)$
has a $\gamma$-bounded
$\Hr^{2}_v(\sect{\omega_p})$-calculus for each admissible function
$v:\R \to [1, \infty)$. 
In particular, $A_p$ has a $\gamma$-bounded $\Hr^{\beta, 2}(\sect{\omega_p})$-calculus for each $\beta > \frac{1}{2}$.  
\end{thm}

\vanish{
\begin{rem}
Theorem \ref{app.t.OU} is  optimal in that neither the
angle $\omega_p$ can be improved nor the smoothness condition on $v$. 
The latter is trivial, and the former follows from the estimate
\[
\| \mathcal{L}^{-\ui s}\|\, \geq\, \ue^{\varphi_p |s|}  \qquad (s\in \R).
\]
which is implicit in \cite{GCMMST2001}, 
\end{rem}
}

\appendix

\section{Auxiliary Lemmas}\label{s.aux}

The following is related to Fubini's theorem. 

\begin{lem}\label{aux.l.L1}
Let $f,g \in \Ell{1}(\R)$. Then 
\begin{aufzi}
\item $\tau_tf \cdot g \in \Ell{1}(\R)$ for almost all $t\in \R$;
\item $\tau_\bft f\cdot g\in \Ell{1}(\R; \Ell{1}(\R))$;
\item $\dps \int_\R \tau_tf \cdot g \, \ud{t} = \Bigl(\int_\R f(t)\, \ud{t} \Bigr)\,
  g$ as an integral in $\Ell{1}(\R)$.
\end{aufzi}
\end{lem}

\begin{proof}
  First let $f,g \in \Cc(\R)$. Then $\tau_tf \cdot
  g \in \Ell{1}(\R)$ for {\em all} $t\in \R$ and $\tau_\bft f
  \cdot g \in \Cc(\R; \Ell{1}(\R) \cap \Co(\R))$. As point evaluations
are continuous on $\Co(\R)$, 
\[  \int_\R \tau_t f\cdot g\, \ud{t}(x) = \int_\R f(x-t)g(x)\, \ud{t}= 
\int_\R f(t)\, \ud{t} \cdot g(x) \qquad (x\in \R),
\]
and hence c) holds.  It follows that 
\[   \Bigl(\int_\R \abs{f} \Bigr) \abs{g} = \int_\R \tau_t
\abs{f}\cdot \abs{g}\, \ud{t}= 
\int_\R \abs{ \tau_t f \cdot g}\, \ud{t}
\]
and integrating (which amounts to  applying a  bounded linear functional on
$\Ell{1}(\R)$) yields
\[  \norm{f}_1 \norm{g}_1 = \int_\R \norm{\tau_t f\, \cdot g}_1\,
\ud{t}.
\]
Thus we obtain a bounded bilinear mapping
\[ \Phi: \Ell{1}(\R) \times \Ell{1}(\R) \to \Ell{1}(\R ;
\Ell{1}(\R)),\qquad
\Phi(f,g) := \tau_\bft f\cdot g\quad \text{if $f,g\in \Cc(\R)$}.
\]
Next, let $f,g$ be arbitrary and take $f_n, g_n
\in \Cc(\R)$ with $f_n \to f$ and $g_n \to g$ in $\Ell{1}(\R)$. 
Without loss of generality $f_n \to f$ and $g_n \to g$ almost
everywhere
and $\tau_\bft f_n \cdot g_n = \Phi(f_n, g_n) \to \Phi(f,g)$ almost everywhere.
On the other hand, 
$\tau_t f_n\cdot g_n \to \tau_t f \cdot
g$ almost everywhere for {\em each} $t\in \R$, hence $\Phi(f,g)(t) = \tau_t f \cdot
g$ in $\Ell{1}(\R)$ for almost all $t\in \R$. This proves a) and b). And c) follows since
integration is a bounded linear operator/functional. 
\end{proof}

\begin{lem}\label{aux.l.Lp}
Let $f\in \Ell{1}(\R) \cap \Ell{\infty}(\R)$ and  $g\in
\Ell{p}(\R)$, $1\le p < \infty$. Then $\tau_\bft f \cdot g \in \Cb(\R; \Ell{p}(\R))$ and
\[  \int_\R \tau_t f \cdot g\, \ud{t} = \Bigl(\int_\R f\, \Bigr) \, g
\]
as a weak integral. 
\end{lem}

\begin{proof}
Write $\tau_t f \cdot g = \tau_t( f \cdot \tau_{-t}g)$ to see that
the function $\tau_\bft f \cdot g$ is continuous with values in
$\Ell{p}(\R)$. To see that it is weakly integrable, let $h\in
\Ell{p'}(\R)$ and compute
\begin{align*}
\int_\R & \dprod{ \tau_t f \cdot g}{h}\, \ud{t} = 
\int_\R \dprod{ \tau_t f \cdot (gh)}{\car}\, \ud{t} \stackrel{(*)}{=}
\dprod{ \int_\R \tau_t f \cdot (gh)\, \ud{t}}{\car} 
\\ & =
\Bigl( \int_\R f \Bigr) \,   \dprod{ gh }{\car} =
\Bigl( \int_\R f \Bigr) \,   \dprod{ g }{h},
\end{align*}
where we used Lemma \ref{aux.l.L1} and that $gh \in \Ell{1}(\R)$ at ($*$).
\end{proof}

\begin{lem}\label{aux.l.exp}
Let $1\le p < \infty$ and 
$f,g $ be measurable such that there are $\alpha >1$ and $\beta >
\frac{1}{p}$ 
with 
\[ \abs{f(s)} \lesssim \frac{1}{(1 + \abs{s})^\alpha},\quad  
\abs{g(s)} \lesssim \frac{1}{(1 + \abs{s})^{\alpha +\beta}} \qquad (s\in \R).
\]
Then the function 
$\tau_\bft f \cdot g$ is in $ \Cb(\R; \Ell{p}(\R))$ and 
\[  \int_\R \norm{\tau_t f\cdot g}_{\Ell{p}} \, \ud{t}< \infty
\quad \text{and}\quad 
\int_\R \tau_t f\cdot g \, \ud{t} = \Bigl(\int_\R f\Bigr)\, g.
\]
\end{lem}

\begin{proof}
We note that $f\in \Ell{1}\cap \Ell{\infty}$ and hence Lemma \ref{aux.l.Lp}
is applicable. We only need to show the strong integrability. 
From $1 + \abs{x + y} \le (1 + \abs{x}) (1 + \abs{y})$ it follows 
(with $x= t-s$ and $y=s$) that  
\[   \abs{\tau_t f \cdot g}(s) 
\lesssim  \frac{1}{(1 + \abs{t-s})^{\alpha}}   \frac{1}{(1 + \abs{s})^{(\alpha + \beta)}}
\le \frac{1}{(1 + \abs{t})^{\alpha}} \frac{1}{(1 + \abs{s})^{\beta}}
\]
and hence 
\[ \norm{\tau_t f\cdot g}_p \lesssim  \frac{1}{(1+ \abs{t})^{\alpha}} \qquad
(t\in \R).\qedhere
\]
\end{proof}

\vanish{
\begin{lem}\label{app.l.exp}
Let $f,g $ be measurable such that there is $c >0$ with 
\[ \abs{f(s)}, \abs{g(s)} \lesssim \ue^{-c\abs{s}} \qquad (s\in \R).
\]
Then for each $1\le p < \infty$ the function 
$\tau_\bft f \cdot g$ is in $ \Cb(\R; \Ell{p}(\R))$ and 
\[  \int_\R \norm{\tau_t f\cdot g}_{\Ell{p}} \, \ud{t}< \infty
\quad \text{and}\quad 
\int_\R \tau_t f\cdot g \, \ud{t} = \Bigl(\int_\R f\Bigr)\, g.
\]
\end{lem}

\begin{proof}
We note that $f\in \Ell{1}\cap \Ell{\infty}$ and hence Lemma \ref{aux.l.Lp}
is applicable. We only need to show the strong integrability. 
Fix $\theta \in (0,1)$. Then 
\[   \abs{\tau_t f \cdot g}(x) 
\lesssim \ue^{-c\abs{t-x}} \ue^{-c\abs{x}} 
\le  \ue^{-\theta c \abs{t-x}} \ue^{-\theta c\abs{x}}\ue^{-(1{-}\theta)c\abs{x}} 
\le \ue^{-\theta c \abs{t}}\ue^{-(1{-}\theta)c\abs{x}}
\]  
and hence 
\[ \norm{\tau_t f\cdot g}_p \lesssim \ue^{-\theta c\abs{t}} \cdot
\Bigl( \frac{2}{(1{-}\theta) p c}\Bigr)^\frac{1}{p} \qquad (t\in \R).\qedhere
\]
\end{proof}
}

\section{Abstract Littlewood--Paley Theory}\label{s.alp}

In this appendix we provide some definitions and results from
the theory of abstract square function estimates and their connection with
functional calculus as promoted by Kalton and Weis in their seminal work \cite{KalWei2016}.
A more detailed account can be found in
\cite{HaHa2013},  an exhaustive treatment in \cite{HvNVW2017}.

\medskip

In the following, $X,\,Y$ are Banach spaces and 
$(\gamma_j)_{j\in J}$ denotes, generically, a family (over a
context-depending index set $J$)  of
independent $\C$-valued normalized Gaussian random variables.
We assume that the reader is familiar with the notion 
of (Rademacher) type and cotype as developed in \cite{DJT}, see
also \cite[Chap. 7]{HvNVW2017}. 
\vanish{
Let $x_1, \dots, x_{n} \in X$. Then, since
$\gamma_n$ is symmetric, 
\begin{align*}
    \norm{\sum_{j=1}^{n{-}1} \gamma_j x_j}_{\Ell{2}(\Omega;X)}
& \le \frac{1}{2} \norm{\sum_{j=1}^{n{-}1} \gamma_j x_j + \gamma_n x_n }_{\Ell{2}(\Omega;X)} 
+ \frac{1}{2} \norm{\sum_{j=1}^{n{-}1} \gamma_j x_j - \gamma_n x_n
}_{\Ell{2}(\Omega;X)} 
\\ & = \norm{\sum_{j=1}^{n} \gamma_j x_j}_{\Ell{2}(\Omega;X)}.
\end{align*}
}

Let $(x_j)_{j \in J}$ be any family in $X$. It is easy to see (by the
symmetry of a standard Gaussian) that the net
\[   \bigl(\Exp \bignorm{   \sum_{j \in F} \gamma_j x_j}_X^2
\bigr)_{F\subseteq J\, \text{finite}}
\]
is increasing. The family $(x_j)_j$ is called \emdf{$\gamma$-summing}
if
this net is bounded, i.e. if 
\[   \norm{ (x_j)_j }_{\gamma(J; X)} := \sup_F  \Bigl(\Exp \bignorm{
  \sum_{j\in F}
\gamma_n x_j }_X^2\Bigr)^\frac{1}{2} < \infty,
\]
where the supremum is taken over all finite subsets $F \subseteq J$.
The space
of all $\gamma$-summing families $(x_j)_j \in X^J$, endowed with the
norm $\norm{\cdot}_{\gamma(J;X)}$, is denoted by $\gamma_\infty(J; X)$.
The space of 
\emdf{$\gamma$-radonifying} families is the closure
\[  \gamma(J;X) := \cls{\coo(J;X)}
\]
of the set of finite sequences within $\gamma_\infty(J;X)$. By
a classical result, 
$\gamma_\infty(J;X)=
\gamma(J;X)$ if $X$ has finite cotype, see \cite[Thm. 6.4.10]{HvNVW2017}.

\medskip

A Banach space $X$ has \emdf{(Pisier's) contraction
  property}\footnote{This property is also called ``property $(\alpha)$'',
  but we follow the suggestion in \cite{HvNVW2017} to change that terminology.}  
if the canonical mapping
\[   \coo(\Z^2;X) \to \coo(\Z; \gamma(\Z;X)),\qquad 
(x_{ij})_{i,j} \mapsto \bigl((x_{ij})_{i} \bigr)_{j}
\]
extends to an isomorphism
\[     \gamma(\Z^2; X )  \cong \gamma(\Z; \gamma(\Z; X)).
\]%
\vanish{
\[  \norm{ (x_{ij})_{i,j} }_{\gamma(\Z \times \Z;X)}
\eqsim 
\bignorm{ \bigl((x_{ij})_{i} \bigr)_{j} }_{\gamma(\Z ; \gamma(\Z;X))} 
\qquad ((x_{ij})_{i,j} \in \coo(\Z\times \Z; X)).
\]
there is $C > 0$ such that
\[  \frac{1}{C}   
\norm{ (x_{ij})_{(i,j) \in I\times J} }_{\gamma(I \times J;X)}
\le \norm{ \bigl((x_{ij})_{i\in I} \bigr)_{j  \in J} }_{\gamma(J ; 
\gamma(I;X))}
\le C \norm{ (x_{ij})_{(i,j) \in I\times J} }_{\gamma(I \times J;X)}
\]
for all finite sets  $I, J$ and all families $(x_{ij})_{i,j} \in X^{I
  \times J}$.}%
If $X$ has Pisier's contraction property, then $X$ has finite
cotype. For further information, see  \cite[Sec. 7.5]{HvNVW2017}.

\begin{rems}\label{alp.r.geom}
1)\ If $X$ has finite cotype, Gaussian and Rademacher random sums
are equivalent in the sense that the space $\gamma_\infty(J;X)$
does not change (up to an equivalent norm) if in its definition
one employs Rademacher variables instead of Gaussians,
see \cite[Cor. 7.2.10]{HvNVW2017}. 

Similarly, Pisier's contraction property  can equivalently be defined by employing
Rademacher variables instead of Gaussian variables (because either
formulation implies that $X$ has finite cotype).

\medskip
\noindent
2)\ Let $X\coloneqq \Ell{p}(\Omega', \mu)$ for some measure space
$(\Omega', \mu)$ and $1\le p \le \infty$. Then, 
\begin{enumerate}
\item $X$ has Pisier's contraction property  if $1 < p < \infty$;
\item $X$ has type $\min \{  2, \, p\}$, and
\item $X$ has cotype $\max \{  2, \, p\}$. In
  particular, 
\item $X$ is of finite cotype if $1\le p < \infty$.
\end{enumerate}
\end{rems}

\medskip

A set $\calT \subseteq \BL(X;Y)$ is called \emdf{$\gamma$-bounded} if
there is a constant $C\ge 0$ such that 
\[   \norm{ (T_j x_j)_{j\in J} }_{\gamma(J;Y)} \le C \norm{ (x_j)_j}_{\gamma(J;X)}
\]
for all finite families $(T_j, x_j)_j$ in $\calT \times X$. In this
case, the least
number $C$ with this property is denoted by $\gammabound{\calT}$ and 
called the \emdf{$\gamma$-bound} of $\calT$. Clearly, if
$\calT'\subseteq \calT$ then $\gammabound{\calT'}\le \gammabound{\calT}$.
The proof of the following lemma is straightfoward.

\begin{lem}\label{alp.l.gbd}
Let $X,Y$ be Banach space and $\calT \subseteq \BL(X;Y)$ $\gamma$-bounded.
Then any (arbitrary)  family 
 $(T_j)_{j \in J}$ in $\calT$ induces a bounded ``diagonal operator''
\[        T := (T_j)_j: \gamma_\infty(J;X) \to \gamma_\infty(J;Y),\qquad (x_j)_j
\mapsto (T_j x_j)_j. 
\]
Moreover, $\norm{T}_{\gamma_\infty(J;X)\to \gamma_\infty(J;Y)} \le \gammabound{\calT}$.
\end{lem}

\vanish{
\begin{defn}\label{prel.def.GammaTwo}
An algebra homomorphism $\Phi : \mathcal{A}\to \BL (X)$, where $\mathcal{A}$ is some normed algebra, is called $\gamma$-\emdf{bounded} if the set
\[
\left\{ \Phi (a)\colon \|a\|_{\mathcal{A}}\leq 1 \right\}
\] is a $\gamma$-bounded subset of $\BL (X)$. 
\end{defn}
}

\begin{rem}
Exchanging Gaussian variables for Rademachers yields the notion of an \emdf{$R$-bounded}
  set of operators. However, as long as $X, Y$ have finite cotype,
$\gamma$-boundedness is equivalent to $R$-boundedness. 
\end{rem}

\vanish{
We note the following two results for easy reference. Both statements show how $\gamma$-bounded sets can be obtained in certain situations.

\begin{thm}\label{prel.thm.ClosureStillGammaBounded}
Let $X$ be a Banach space, and let $\calT\subseteq \BL (X)$ be a $\gamma$-bounded set. Then, both the closure $\overline{\calT}^{\, \mathrm{st}}$ in the strong operator topology and the closure $\overline{\calT}^{\,\mathrm{we}}$ in the weak operator topology are $\gamma$-bounded with
\[
\gammabound{ \overline{\calT}^{\, \mathrm{st}}}\, =\, \gammabound{\overline{\calT}^{\, \mathrm{we}}}\, =\, \gammabound{\calT}.
\]
\end{thm}
}

\begin{thm}\label{alp.t.intmeans}
Let $X$ be a Banach space, let
$(\Omega, \mu)$ be a measure space, and let $r\in [1, \infty)$ with
\[
\frac{1}{r}\, >\, \frac{1}{\mathrm{type}(X)} - \frac{1}{\mathrm{cotype}(X)}.
\] 
Furthermore,  let $F: \Omega \to \BL (X)$ be strongly
$\mu$-measurable with 
\[
\sup_{\| x\|\leq 1} \| F(\bfs)x\|_{\Ell{r}(\Omega;X)}\, <\, \infty.
\] 
Consider for each $\varphi \in \Ell{r'}(\Omega , \mu)$ the operator
\[
I_{F, \varphi} : X\to X, \quad x\mapsto \int_\Omega \varphi (s) \cdot F(s)x\, \ud \mu (s).
\]  
Then, the set $\calT_{F, r'} \coloneqq \{ I_{F, \varphi}\, |\,  \|\varphi\|_{r'}\leq 1 \}$ is $\gamma$-bounded with
\[
\gammabound{\calT_{F, r'}}\, \lesssim  \sup_{\|x\|\leq 1} \| F(\bfs)x\|_{\Ell{r}(\Omega;X)}.
\]
\end{thm}

\begin{proof}
The hypotheses imply that $X$ has finite cotype and, thus,  Rademacher
and Gaussian sums are equivalent. Hence, 
the result is contained in \cite[Theorem 8.5.12]{HvNVW2017}.
\end{proof}

\begin{rem}

Theorem \ref{alp.t.intmeans} also holds in the case $p=1$ and $q=
\infty$, see \cite[Thm. 8.5.4]{HvNVW2017} and exchange Gaussians
for Rademachers in the proof. 
\end{rem}

A subset $\calT \subseteq \BL(X;Y)$ is called \emdf{semi-$\gamma$-bounded}
if 
\[  \sup_{\norm{x}\le 1} \gammabound{ \C \ni \lambda \mapsto \lambda Tx
  \suchthat T \in \calT} < \infty.
\]
If $X$ has finite cotype, semi-$\gamma$-boundedness is the same as
semi-$R$-boundedness, a notion that has been introduced and studied in
\cite {VeWe2010}.  
We need the following auxiliary result.

\begin{thm}\label{alp.t.semi}
Let $X$ be a Banach space with Pisier's property $(\alpha)$ and let
$(T_s)_{s\in \R}$ be a strongly measurable and semi-$\gamma$-bounded
family.
Then 
\[    \gammabound{  x\mapsto \int_\R f(s)g(s) T_sx \, \ud{s} \suchthat 
f\in \Ell{2}(\R),\, \norm{f}_2\le 1} \lesssim \norm{g}_2 \qquad (g\in
\Ell{2}(\R)).
\]
\end{thm}

\begin{proof}
Let $H := \Ell{2}(\R)$ and  let $\gamma(H;X)$ be the set of operators
$T: H \to X$ such that $(Te_n)_{n \in \N} \in \gamma(\N;X)$ for
one/each orthonormal basis $(e_n)_n \in H$.
These are the so-called {\em $\gamma$-radonifying opertors}, see
\cite[Sec. 9.1.b]{HvNVW2017}. Identifying $H '= \Ell{2}(\R)'$ with
$\Ell{2}(\R) = H$ via the canonical duality we obtain $\gamma(H;\C) =
H$. Fix $g\in H$ and $x\in X$ and consider the operator
\[   S_{g,x} : H \to X,\qquad f \mapsto \int_\R f(s)g(s)T_s x\, \ud{s}.
\]
Since the family $(T_s)_{s\in \R}$ is semi-$\gamma$-bouded, the
so-called {\em $\gamma$-multiplier theorem}
\cite[Thm. 9.5.1]{HvNVW2017} implies
that $S_{g,x}\in \gamma(H;X)$ and 
\[   \norm{T_{g,x}}_\gamma \lesssim \norm{g}_H \norm{x} \qquad (x\in
X,\, g\in H).
\]
Since $X$ has property $(\alpha)$,
\[  \gammabound{ \gamma(H;X) \ni S \mapsto Sf \in X\suchthat f\in H, \norm{f}_H \le 1} <
\infty.
\]
Composing this with the bounded operator
$x \mapsto T_{g,x}$ yields
\[  \gammabound{ x \mapsto T_{g,x}f 
\suchthat \norm{f}_H \le 1} 
\lesssim \norm{x \mapsto T_{g,x}}_{X\to \gamma(H;X)} \lesssim
\norm{g}_H \qquad (g\in H)
\]
as claimed. 
\end{proof}

The following is an abstract Littlewood--Paley
theorem.

\begin{prop}\label{apl.p.PL}
Let $A$ be a densely defined operator on a Banach space
$X$ with finite cotype. Suppose that $A$ has a
 bounded $\Ha (\strip{\theta})$-calculus 
for some $\theta > 0$ and $\psi, \eta 
\in \calE[ \cstrip{\theta}]$ are such that 
\[  \sum_{n \in \Z} \psi(t-n)\eta(t-n) = 1 \quad \text{for all $t\in \R$}. 
\]
Then  one has
the norm equivalence
\[   \norm{x} \eqsim \norm{ (\psi(A-n)x )_{n\in \Z} }_{\gamma(\Z; X)}
\qquad (x\in X). 
\]
\end{prop}

\begin{proof}
This is similar to \cite[Thm.s 10.4.4 and 10.4.8]{HvNVW2017}.
Note that 
\[   \sup_{z\in \strip{\theta}} \norm{ (\psi(z-n))_{n\in
    \Z}}_{\ell^1(\Z)} < \infty,
\]
which in the terminology of \cite{HaHa2013} means that 
the mapping 
\[  \strip{\theta} \to \ell^2(\Z), \quad z \mapsto(\psi(z-n))_n
\]
has $\ell^1$-frame bounded image. This implies that the associated
square function and dual square function
are bounded \cite[Thm. 4.11]{HaHa2013}.  The additional hypothesis implies
(by holomorphy)
\[  \sum_{n\in \Z} \psi(z-n)\eta(z-n) = 1 \qquad (z\in \strip{\theta}),
\]
and hence the norm equivalence holds by 
\cite[Cor. 4.7]{HaHa2013}.  
\end{proof}

\bigskip

\noindent
{\bf Acknowledgements. } Florian Pannasch would like to thank Christoph
Kriegler 
and Lutz Weis for inspiring discussions on Hörmander-type 
functional calculus during the 21st Internet Seminar on ``Functional
Calculus'' (Wuppertal 2018). Markus Haase is grateful to Lukas
Hagedorn  
for giving helpful comments on earlier versions of the paper.

\nocite{GCMST1999}
\nocite{Stein1970}
\nocite{Str1967}


\begin{thebibliography}{GCMM{\etalchar{+}}01}

\bibitem[Blu03]{Blunck2003}
Blunck, S.:
\newblock A {H}\"{o}rmander-type spectral multiplier theorem for operators
  without heat kernel.
\newblock {\em Ann. Sc. Norm. Super. Pisa Cl. Sci. (5)}, 2(3):449--459, 2003.

\bibitem[CD17]{CarDra2017}
Carbonaro, A. and Dragi\v{c}evi\'{c}, O.:
\newblock Functional calculus for generators of symmetric contraction
  semigroups.
\newblock {\em Duke Math. J.}, 166(5):937--974, 2017.


\bibitem[CFMP05]{ChFaMePa2005}
Chill, R. and Fa\v{s}angov\'{a}, E.  and Metafune, G.  and Pallara,
D.:
\newblock The sector of analyticity of the Ornstein-Uhlenbeck
semigroup on $\Ell{p}$-spaces with respect to invariant measure.
\newblock {\em  J. London Math. Soc. (2)}  71(3): 703--722, 2005.

		


\bibitem[Chr91]{Christ1991}
Christ, M.:
\newblock {$L^p$} bounds for spectral multipliers on nilpotent groups.
\newblock {\em Trans. Amer. Math. Soc.}, 328(1):73--81, 1991.

\bibitem[CM93]{CoMe1993}
Cowling, M. and Meda, S.
\newblock Harmonic analysis and ultracontractivity.
\newblock {\em Trans. Amer. Math. Soc.}, 340(2):733--752, 1993.

\bibitem[Cow83]{Co1983}
Cowling, M.:
\newblock Harmonic analysis on semigroups.
\newblock {\em Ann. of Math. (2)}, 117(2):267--283, 1983.


\bibitem[DJT]{DJT}
Diestel, J. and  Jarchow, H. and Tonge, A.:
\newblock {\em Absolutely Summing Operators.}
\newblock Cambridge University Press, 1995.



\bibitem[GCM{\etalchar{+}}01]{GCMMST2001}
Garc\'{i}a-Cuerva, J. and Mauceri, G. and Meda, S. 
and Sj\"{o}gren, P. and Torrea, J. L.:
\newblock Functional calculus for the {O}rnstein-{U}hlenbeck operator.
\newblock {\em J. Funct. Anal.}, 183(2):413--450, 2001.

\bibitem[GCM{\etalchar{+}}99]{GCMST1999}
Garc\'{\i}a-Cuerva, J. and  Mauceri, G. and Sj\"{o}gren, P. and
  Torrea, J. L.:
\newblock Spectral multipliers for the {O}rnstein-{U}hlenbeck semigroup.
\newblock {\em J. Anal. Math.}, 78:281--305, 1999.


\bibitem[Glu15]{Glu2015}
Glück, J.:
\newblock On weak decay rates and uniform stability of bounded linear
operators.
\newblock  {\em Arch. Math}, 104:347--356, 2015.


\bibitem[Gra09]{Gra2009}
Grafakos, L.:
\newblock {\em Modern {F}ourier analysis}, volume 250 of {\em Graduate Texts in
  Mathematics}.
\newblock Springer, New York, second edition, 2009.

\bibitem[HH13]{HaHa2013}
Haak, B. H. and Haase, M:.
\newblock Square function estimates and functional calculi. 
\newblock Unpublished preprint, 2013. 
Arxiv:  \url{https://arxiv.org/abs/1311.0453}.
\newblock Extended version to appear as
{\em The {S}quare {F}unctional {C}alculus}.



\bibitem[Haa05]{Ha2005}
Haase, M.:
\newblock A general framework for holomorphic functional calculi.
\newblock {\em Proc. Edinb. Math. Soc. (2)}, 48(2):423--444, 2005.

\bibitem[Haa06]{Ha2006}
Haase, M.:
\newblock {\em The functional calculus for sectorial operators}, volume 169 of
  {\em Operator Theory: Advances and Applications}.
\newblock Birkh\"{a}user Verlag, Basel, 2006.

\bibitem[Haa18]{Ha2018}
Haase, M.:
\newblock {\em Functional calculus.}
\newblock {\em Lecture notes of the 21st Internet Seminar}, 2018.

\bibitem[Haa20]{Ha2020}
Haase, M.:
\newblock On the fundamental principles of unbounded functional   calculi.
\newblock submitted manuscript, 2020. 
\url{https://doi.org/10.48550/arXiv.2104.04865}

\bibitem[Heb90]{Heb1990}
Hebisch, W.:
\newblock A multiplier theorem for {S}chr\"{o}dinger operators.
\newblock {\em Colloq. Math.}, 60/61(2):659--664, 1990.



\bibitem[Hoe60]{Hoe1960}
H\"{o}rmander, L.:
\newblock Estimates for translation invariant operators in {$L^{p}$}spaces.
\newblock {\em Acta Math.}, 104:93--140, 1960.


\bibitem[HvNVW17]{HvNVW2017}
Hyt\"{o}nen, T. and van Neerven, J. and Veraar, M. and Weis, L.:
\newblock {\em Analysis in {B}anach spaces. {V}ol. {II}}. 
\newblock Springer, Cham, 2017.


\bibitem[KW16]{KalWei2016}
Kalton, N. and Weis, L.:
\newblock { \em The $H^\infty$-calculus and square function estimates.}
\newblock Pages 716--764 in: Nigel J. Kalton Selecta,  Volume 1,
Birkhäuser, Cham 2016. 

\bibitem[Kat]{Katznelson2004}
Katznelson, Y.:
\newblock {\em An Introduction to Harmonic Analysis}.
\newblock Third corrected edition, Cambridge UP, 2004

\bibitem[Kri09]{Kr2009}
Kriegler, C.:
\newblock {\em Spectral {M}ultipliers, {$R$}-bounded {H}omomorphisms, and {A}nalytic
  {D}iffusion {S}emigroups.}
\newblock Ph.D. Thesis, 2009.
\newblock URL: \url{http:
//digbib.ubka.uni-karlsruhe.de/volltexte/1000015866}.

\bibitem[Kri14]{Kr2014Eins}
Kriegler, C.:
\newblock H\"{o}rmander functional calculus for {P}oisson estimates.
\newblock {\em Integral Equations Operator Theory}, 80(3):379--413, 2014.

\bibitem[KW17]{KrWe2017}
Kriegler, C. and Weis, L.:
\newblock Spectral multiplier theorems and averaged {$R$}-boundedness.
\newblock {\em Semigroup Forum}, 94(2):260--296, 2017.

\bibitem[KW18]{KrWe2018}
Kriegler, C. and Weis, L.:
\newblock Spectral multiplier theorems via {$H^\infty$} calculus and
  {$R$}-bounds.
\newblock {\em Math. Z.}, 289(1-2):405--444, 2018.

\bibitem[Mue98]{Muller1998}
M\"{u}ller, D.:
\newblock Functional calculus on {L}ie groups and wave propagation.
\newblock In {\em Proceedings of the {I}nternational {C}ongress of
  {M}athematicians, {V}ol. {II} ({B}erlin, 1998)}, pages
  679--689, 1998.

\bibitem[Med90]{Me1990}
Meda, S.:
\newblock A general multiplier theorem.
\newblock {\em Proc. Amer. Math. Soc.}, 110(3):639--647, 1990.

\bibitem[MM90]{MaMe1990}
Mauceri, G. and Meda, S.:
\newblock Vector-valued multipliers on stratified groups.
\newblock {\em Rev. Mat. Iberoamericana}, 6(3-4):141--154, 1990.

\bibitem[MMS04]{MMS2004}
Mauceri, G., Meda, S., and Sj\"{o}gren, P.:
\newblock Sharp estimates for the {O}rnstein-{U}hlenbeck operator.
\newblock {\em Ann. Sc. Norm. Super. Pisa Cl. Sci. (5)}, 3(3):447--480, 2004.


\bibitem[Pa19]{Pan2019}
Pannasch, F.:
\newblock{ \em The Holomorphic Hörmander Functional Calculus}.
\newblock Ph.D. Thesis, Kiel University, 2019.
\newblock URL: \url{https://macau.uni-kiel.de/receive/macau_mods_00000194}.

\bibitem[Rud]{RudinFA}
Rudin, W.:
\newblock {\em Functional Analysis}, 
\newblock 2nd edition, McGraw Hill 1991.


\bibitem[Sas05]{Sasso2005}
Sasso, E.:
\newblock Functional calculus for the {L}aguerre operator.
\newblock {\em Math. Z.}, 249(3):683--711, 2005.

\bibitem[Ste70]{Stein1970}
Stein, E. M.:
\newblock {\em Topics in harmonic analysis related to the {L}ittlewood-{P}aley
  theory}.
\newblock Annals of Mathematics Studies, No. 63. Princeton University
Press, 1970.

\bibitem[Str67]{Str1967}
Strichartz, R. S.:
\newblock Multipliers on fractional {S}obolev spaces.
\newblock {\em J. Math. Mech.}, 16:1031--1060, 1967.



\bibitem[Str71]{Str1971}
Strichartz, R. S.:
\newblock A note on {S}obolev algebras.
\newblock {\em Proc. Amer. Math. Soc.}, 29:205--207, 1971.


\bibitem[Str72]{Str1972}
Stromberg, K.:
\newblock An elementary proof of {S}teinhaus' theorem.
\newblock {\em Proc. Amer. Math. Soc.}, 36(1):308, 1972.


\bibitem[VW10]{VeWe2010}
Veraar, M. and Weis, L.:
\newblock On semi-$T$-boundedness and its applications.
\newblock {\em 
J. Math. Anal. Appl.}, 363: 431–443, 2010. 
\end{thebibliography}
\bibliographystyle{alpha}

\newcommand{\etalchar}[1]{$^{#1}$}


\end{document}